\documentclass[reqno]{amsart}
\usepackage{amsfonts}
\usepackage{graphicx}
\usepackage{amsmath,amssymb,amsfonts,amsthm,mathrsfs}
\usepackage[english]{babel}
\usepackage{hyperref}
\usepackage{marginnote}
\usepackage{color}
\usepackage{enumerate}

\newcommand{\dis}{\displaystyle}

\newcommand{\R}{\mathbb{R}}

\theoremstyle{plain}

\newtheorem{theorem}{Theorem}[section]
\newtheorem{lemma}[theorem]{Lemma}

\newtheorem{corollary}[theorem]{Corollary}

\newtheorem{definition}[theorem]{Definition}
\newtheorem{remark}[theorem]{Remark}
\theoremstyle{remark}
\numberwithin{equation}{section}

\def\v{\varepsilon}

\def\t{\theta}

\def\a{\alpha}

\def\g{\gamma}
\def\d{\delta}

\def\s{\sigma}

\def\f{\frac}

\def\sp{\shortparallel}

\usepackage{tikz}

\newcommand{\dd}{{\rm d}}

\newcommand{\Mr}{\mathcal{R}}

\newcommand{\FL}{\mathbf{L}}

\newcommand{\Fi}{\mathbf{1}}
\newcommand{\la}{\lambda}
\newcommand{\pa}{\partial}
\newcommand{\vep}{{\varepsilon}}

\begin{document}

\title[Boundary layer for diffusive boundary conditions]{Boundary layer solution of the Boltzmann equation for diffusive reflection boundary conditions in half-space}

\author[F.M. Huang]{Feimin Huang}
\address[F.M. Huang]{Academy of Mathematics and Systems Science, Chinese Academy of Sciences, Beijing 100190, China; School of Mathematical Sciences, University of Chinese Academy of Sciences, Beijing 100049, China}
\email{fhuang@amt.ac.cn}

\author[Y. Wang]{Yong Wang}
\address[Y. Wang]{Academy of Mathematics and Systems Science, Chinese Academy of Sciences, Beijing 100190, China; School of Mathematical Sciences, University of Chinese Academy of Sciences, Beijing 100049, China}
\email{yongwang@amss.ac.cn}

\begin{abstract}
We study steady Boltzmann equation in half-space, which arises in the Knudsen boundary layer problem, with diffusive reflection boundary conditions.
Under certain admissible conditions and the source term decaying exponentially,  we establish the existence of  boundary layer solution for both linear and nonlinear Boltzmann equation in half-space  with diffusive reflection boundary condition in $L^\infty_{x,v}$  when the far-field Mach number of the Maxwellian is zero.  The  continuity and the spacial decay  of the solution are  obtained. The uniqueness is established under some constraint conditions.
\end{abstract}

\date{\today}
%\thanks{}
\maketitle

\tableofcontents

\thispagestyle{empty}

%%%%%%%%%%%%%%%%%%%%%%%%%%%%%%%%%%%%%%%%%%%%%%%%%%%%%%%%%%%%%%%%%%

\section{Introduction}

%\subsection{Boltzmann equation}
In the present paper, we consider the steady Boltzmann equation
\begin{equation}\label{1.1}
v_3\cdot \partial_x F=Q(F,F)+\mathcal{S},\quad \   (x,v)\in \mathbb{R}_+\times\mathbb{R}^3,
\end{equation}
with $\R_+=(0,+\infty)$.
The Boltzmann collision term $Q(F_1,F_2)$ on the right is defined in terms of the following bilinear form
\begin{align}\label{1.2}
Q(F_1,F_2)&\equiv\int_{\mathbb{R}^3}\int_{\mathbb{S}^2} B(v-u,\t)F_1(u')F_2(v')\,{d\omega du}
-\int_{\mathbb{R}^3}\int_{\mathbb{S}^2} B(v-u,\t)F_1(u)F_2(v)\,{d\omega du}\nonumber\\
&:=Q_+(F_1,F_2)-Q_-(F_1,F_2),
\end{align}
where the relationship between the post-collision velocity $(v',u')$ of two particles with the pre-collision velocity $(v,u)$ is given by
\begin{equation*}%\label{1.3}
u'=u+[(v-u)\cdot\omega]\omega,\quad v'=v-[(v-u)\cdot\omega]\omega,
\end{equation*}
for $\omega\in \mathbb{S}^2$, which can be determined by conservation laws of momentum and energy
\begin{equation*}%\label{1.3-1}
u'+v'=u+v,\quad |u'|^2+|v'|^2=|u|^2+|v|^2.
\end{equation*}
The Boltzmann collision kernel $B=B(v-u,\theta)$ in \eqref{1.2} depends only on $|v-u|$ and $\theta$ with $\cos\theta=(v-u)\cdot \omega/|v-u|$.   Throughout this paper,  we  consider the hard sphere model, i.e.,
\begin{equation*}%\label{1.4}
B(v-u,\t)=|(v-u)\cdot \omega|
\end{equation*}

There have been many studies on the half-space problem of the steady Boltzmann equation in the literature.  The existence, uniqueness and properties of asymptotic behavior were proved in \cite{BCN} for the linearized Boltzmann equation of a hard sphere gas for the  Dirichlet type boundary condition, see \cite{CGS} for a classification of well-posed kinetic boundary layer problem. Later, the existence of nonlinear boundary layers with small magnitudes and Dirichlet boundary condition for the hard sphere model was established in \cite{UYU}, it was shown that the existence of a solution depends on the Mach number of the far field Maxwellian, and an implicit solvability conditions yielding the co-dimensions of the boundary data (see also \cite{CLY}), and we refer \cite{UYU1,DWY} for the time-asymptotic stability of such boundary layer solution. Wu \cite{Wu} established the existence of unique solution for a modified boundary layer solution with Dirichlet boundary condition in $L^\infty_{x,v}$ space which is used to prove the Hilbert expansion (diffusive scaling) in a disk.  Recently,
Bernhoff-Golse \cite{BG1} offered the existence and uniqueness of a uniformly decaying boundary layer type solution in the situation that gas is in contact with its condensed phase.
For  the specular reflection condition and the solution tends to a global Maxwellian with zero Mach number at the far field, Golse-Perthame-Sulem \cite{GPS-1988} and Huang-Jiang-Wang \cite{HJW} proved the existence, uniqueness and asymptotic behavior  in different  functional space, respectively. %This problem has been extensively studied, first by Sone, Aoki and their collaborators, by means of careful numerical simulations.
%To prove the Hilbert expansion of Boltzmann equation for half-space problem with specular boundary condition, the continuity and uniform estimate in $L^\infty_{x,v}$ are needed,  so in the present paper, we will focus on the existence steady solution of \eqref{1.1}  in the functional space $L^{2}_{x,v}\cap L^\infty_{x,v}$.

\smallskip

For the diffusive reflection boundary condition,  the existence of steady nonlinear Boltzmann solution is proved \cite{DHWZ,EGKM2013} in weighted  $L^\infty_{x,v}$ in smooth bounded domain, and the time-asymptotic stability of such steady solutions is also obtained. For half-space problem, Coron-Golse-Sulem \cite{CGS} proved the existence of solution for linear Boltzmann equation in half-space with zero source term in the functional space $L^\infty(e^{\g x} dx; L^2(|v_3|dv))$. To the best of our knowledge,  there is no result on the existence of solution to the nonlinear boundary layer problem \eqref{1.1} under diffusive reflection boundary condition in half-space with continuity. The purpose of the present paper is to establish the existence and continuity of solution of \eqref{1.1} in half-space with diffusive reflection boundary condition in $L^\infty_{x,v}$ when the Mach number of far-field Maxwellian is zero.  In fact, the continuity of boundary layer solution is very important to close the Hilbert expansion of Boltzmann equation in the initial boundary value problem, see \cite{Guo-Huang-Wang} for instance.

\smallskip

We supplement the Boltzmann equation \eqref{1.1} with the perturbed diffusive reflection boundary conditions
\begin{align}\label{B.C}
F(0,v)|_{v_3>0}=\sqrt{2\pi} \mu(v) \int_{u_3<0} F(0,u) |u_3| du +R(v)
\end{align}
and
\begin{equation}\label{Farfield}
\lim_{x\rightarrow\infty}F(x,v)=p_E^0 \mu(v),
\end{equation}
where $p_E^0>0$ is some given positive constant, and $\mu(v)$ is the normalized global Maxwellian
\begin{equation}\label{1.5}
	\mu(v)=\f{1}{(2\pi)^{\f32}} \, e^{-\f{|v|^2}{2}}.
\end{equation}
%Here the far-field function $F^\infty(v)$ will be determined by the solution itself, and can not be imposed in advance. That means $F^\infty(v)$ is a unknown function.
The boundary function $R(v)$ is defined only for $v_3>0$, and one can extend its definition to $\R^3$ by taking $R(v)=0$ when $v_3\leq 0$. The positive constant $p_E^0$  appear in general when we derive the Knudsen boundary layers in the process of Hilbert expansion of  Boltzmann equation (compressible Euler scaling) with physical boundary. In fact, the $p_E^0$ involves with the boundary pressure of Euler solution.

\smallskip

In general, for given source term $\mathcal{S}(x,v)$ and boundary perturbed term $R(v)$, the boundary value problem (BVP) \eqref{1.1} \eqref{B.C}-\eqref{Farfield} may be not solvable. In fact, for given $\mathcal{S}$ and $R$, we have to replace the far-field condition \eqref{Farfield} by a new one. Hence the main difficulty is to construct the new far-field state and obtain uniform estimates for both the new far-field state and the solution. To achieve our goal, we shall use the $L^2_{x,v}$-$L^\infty_{x,v}$-framework initiated by Guo \cite{Guo2}.

\medskip

We search the solution of \eqref{1.1} in the following form
\begin{equation}\label{1.6}
F(x,v)=p_E^0 \mu(v) + \sqrt{\mu(v)} \, f(x,v).
\end{equation}
 Then \eqref{1.1}, \eqref{B.C}, \eqref{Farfield} are rewritten as
\begin{align}\label{1.7}
\begin{cases}
&\dis v_3\cdot \pa_x f + p_E^0 \FL f= \Gamma(f,f) + \mathfrak{S},\\
&\dis f(0,v)|_{v_3>0}=P_\g f(0,v) + \mathfrak{R}(v),\\
&\dis \lim_{x\to\infty} f(x,v)=0,
\end{cases}
\end{align}
where we have used the notations
\begin{align}
\mathfrak{S}&:=\f{\mathcal{S}}{\sqrt{\mu(v)}},\qquad \mathfrak{R}:=\f{R}{\sqrt{\mu(v)}},\label{2.0}\\
P_\g f(v)&:=\sqrt{2\pi\mu(v)} \int_{u_3<0} |u_3| \sqrt{\mu(u)} f(u) du,\label{2.1}
\end{align}
and the linearized operator  $\FL$  is given by  $\FL=\nu(v)-K$ with  $K=K_2-K_1$ and
\begin{align}\notag
	(K_1f)(v)&=\int_{\mathbb{R}^3}\int_{\mathbb{S}^2} |(v-u)\cdot \omega|\sqrt{\mu(v)\mu(u)}f(u)\,d\omega du,%\label{1.11}
	\\
	(K_2f)(v)&=\int_{\mathbb{R}^3}\int_{\mathbb{S}^2} |(v-u)\cdot \omega|\sqrt{\mu(u)\mu(u')}f(v')\,d\omega du,\nonumber\\
	&\quad+\int_{\mathbb{R}^3}\int_{\mathbb{S}^2}B(v-u,\t)\sqrt{\mu(u)\mu(v')}f(u')\,d\omega du,%\label{1.12}
	\notag
	\\
	\nu(v)&=\int_{\mathbb{R}^3}\int_{\mathbb{S}^2}B(v-u,\t)\mu(u)\,d\omega du\cong 1+|v|,
	\notag\\%\label{1.12-1}
	\Gamma(f,f)&=\frac{1}{\sqrt{\mu(v)}} \,  Q(\sqrt{\mu}f, \, \sqrt{\mu}f).\notag
\end{align}
The null space of the operator $\FL$  is the 5-dimensional space of collision invariants:
$$
\mathcal{N}=\mbox{Ker}\,\FL=span \Big\{\sqrt{\mu},\quad  v\sqrt{\mu},\quad \f12 (|v|^2-3)\sqrt{\mu} \Big\}.
$$
And let $\bf{P}$ denote the projection operator from $L^2(\mathbb{R}^3)$ to $ \mathcal{N}$. We list some useful properties of $\FL$ and $K, \Gamma$ in Appendix \ref{APP}.

\smallskip

For later use we define the velocity weight function
\begin{equation}\label{wt}
w(v)=(1+|v|^2)^{\f{\beta}{2}} e^{\varpi |v|^2},\quad 0\leq \varpi<\f14.
\end{equation}

\medskip

In general, for given $\mathfrak{S}$ and $\mathfrak{R}$, the boundary value problem \eqref{1.7} may be not solvable. In fact, our main result is

\begin{theorem}\label{thm1.1}
Recall the weight function $w(v)$ in \eqref{wt}, and let $\beta\geq 3$ and $0\leq \varpi<\f18$. Let $\s_0\in(0,1)$ be suitably small. Assume $\mathfrak{S}\in \mathcal{N}^{\perp}$ and $\int_{v_3>0} v_3 \sqrt{\mu(v)} \mathfrak{R}(v) dv=0$.
There exist a small $\delta_0>0$ such that if
\begin{equation}\label{2.5}
	\|e^{\s_0 x}\nu^{-1}w\mathfrak{S}\|_{L^\infty_{x,v}} +  |w\mathfrak{R}|_{L^\infty_v(\R^3_+)}:=\delta \leq \delta_0,
\end{equation}
then there is a unique function $\mathfrak{f}^{\infty}:=\mathbb{G}(\mathfrak{S}, \mathfrak{R})$ with
\begin{align}\label{2.6-3}
\begin{split}
&\mathfrak{f}^\infty(v)=\Big\{ \mathfrak{b}_{1}^\infty\, v_1+   \mathfrak{b}_{2}^\infty \,v_2 + \mathfrak{c}^\infty \,  (\f{|v|^2}{2}-\f32)\Big\}\sqrt{\mu(v)}\in\mathcal{N},\\
&|(\mathfrak{b}_1^\infty, \mathfrak{b}_2^\infty, \mathfrak{c}^\infty)|\leq C\Big\{\|e^{\s_0 x}\nu^{-1}w\mathfrak{S}\|_{L^\infty_{x,v}} +  |w\mathfrak{R}|_{L^\infty_v(\R^3_+)}  \Big\},
\end{split}
\end{align}
such that the following nonlinear boundary layer problem of Boltzmann equation
\begin{align}\label{1.7-5}
\begin{cases}
&\dis v_3\cdot \pa_x \mathfrak{f} + p_E^0 \FL \mathfrak{f}= \Gamma(\mathfrak{f},\mathfrak{f}) + \mathfrak{S},\,\, (x,v)\in \R_+\times \R^3,\\
&\dis \mathfrak{f}(0,v)|_{v_3>0}=P_\g \mathfrak{f}(0,v) -(I-P_\g)\mathfrak{f}^{\infty}(v) + \mathfrak{R}(v),\\
&\dis \lim_{x\to\infty} \mathfrak{f}(x,v)=0,
\end{cases}
\end{align}
has a unique mild solution $\mathfrak{f}(x,v)$ satisfying
\begin{align}\label{1.14}
\|e^{\s x}w\mathfrak{f}\|_{L^\infty_{x,v}}+| w\mathfrak{f}(0)|_{L^\infty_v(\R^3)}
\leq \frac{C}{\s_0-\s}\Big\{\|e^{\s_0 x}\nu^{-1}w\mathfrak{S}\|_{L^\infty_{x,v}} +  |w\mathfrak{R}|_{L^\infty_v(\R^3_+)}  \Big\},
\end{align}
where $\sigma>0$ is a constant such $\sigma\in (0,\sigma_0)$, and $C>0$ is a constant independent of $\sigma$.
Moreover, if
$\mathfrak{S}$ is continuous in $\R_+\times \mathbb{R}^3$ and $\mathfrak{R}(v)$ is continuous in $\{v\in\R^3_+\}$, then $\mathfrak{f}(x,v)$ is continuous away from the  grazing set $\big\{(0,v)\ : \ v\in\R^3,\  v_3=0\big\}$.

\smallskip

Furthermore, let $\mathfrak{S}_i\in\mathcal{N}^{\perp}$ and $\int_{v_3>0} v_3 \sqrt{\mu(v)} \mathfrak{R}_i(v) dv=0$, $i=1,2$ satisfying \eqref{2.5}. Let $\mathfrak{f}_i, \mathfrak{f}^\infty_i$ be the solutions obtain above by replacing $\mathfrak{S},\, \mathfrak{R}$ by $\mathfrak{S}_i,\, \mathfrak{R}_i$, and  denote
\begin{equation}
\mathfrak{f}_i^\infty(v)= \mathbb{G}(\mathfrak{S}_i,\mathfrak{R}_i)=\Big\{ \mathfrak{b}_{i,1}^\infty\, v_1+   \mathfrak{b}_{i,2}^\infty \,v_2 + \mathfrak{c}_i^\infty \,  (\f{|v|^2}{2}-\f32)\Big\}\sqrt{\mu(v)}\,\,\, i=1,2,
\end{equation}
then it holds that $\mathbb{G}(0,0)=0$ and
\begin{align}\label{1.30}
&\|e^{\s x}w(\mathfrak{f}_1-\mathfrak{f}_2)\|_{L^\infty_{x,v}}+|w(\mathfrak{f}_1-\mathfrak{f}_2)(0)|_{L^\infty_v(\R^3)} +|(\mathfrak{b}_{1,1}^\infty-\mathfrak{b}_{2,1}^\infty, \, \mathfrak{b}_{1,2}^\infty-\mathfrak{b}_{2,2}^\infty, \, \mathfrak{c}_{1}^\infty-\mathfrak{c}_{2}^\infty)| \nonumber\\
&\leq \frac{C}{\s_0-\s}\Big\{\|e^{\s_0 x}\nu^{-1}w(\mathfrak{S}_1-\mathfrak{S}_2)\|_{L^\infty_{x,v}} +  |w(\mathfrak{R}_1-\mathfrak{R}_2)|_{L^\infty_v(\R^3_+)}  \Big\}.
\end{align}
That means the solution of $\mathfrak{f}, \mathfrak{f}^\infty$ of \eqref{1.7-5} depend continuously on $\mathfrak{S},\, \mathfrak{R}$ in the sense of \eqref{1.30}.
\end{theorem}

\begin{remark}
We note that $(I-P_\g)(\mathfrak{a}^\infty\sqrt{\mu}+\mathfrak{f}^\infty)\equiv (I-P_\g)(\mathfrak{f}^\infty)$ for all $\mathfrak{a}^\infty \in \R$, so the uniqueness of $\mathfrak{f}^\infty$ is in the sense that we normalize $\mathfrak{a}^\infty=0$. In fact, the uniqueness of $\mathfrak{f},\mathfrak{f}^{\infty}$ is under the constraints conditions \eqref{2.6-3} and \eqref{1.14}.
\end{remark}

\begin{corollary}\label{cor1.3}
Recall the weight function $w(v)$ in \eqref{wt}, and let $\beta\geq 3$ and $0\leq \varpi<\f18$. Let $\s_0\in(0,1)$ be suitably small. Assume $\mathfrak{S}\in \mathcal{N}^{\perp}$ and $\int_{v_3>0} v_3 \sqrt{\mu(v)} \mathfrak{R}(v) dv=0$. If there is a function $\mathfrak{r}(v)$ with $w\mathfrak{r} \in L^\infty_v(\R^3_+)$ such that
 \begin{equation}
 \mathfrak{R}(v)=-(I-P_\g)\mathbb{G}(\mathfrak{S},\, \mathfrak{r})(v) + \mathfrak{r}(v),
 \end{equation}
and $ \|e^{\s_0 x}\nu^{-1} w \mathfrak{S}\|_{L^\infty_{x,v}}+\| w\mathfrak{r}\|_{L^\infty_v(\R^3_+)}\leq \delta_0$, then there is a unique mild solution $f(x,v)$ to the nonlinear boundary layer problem   \eqref{1.7}, and satisfies
\begin{align}
\begin{split}
 |\mathbb{G}(\mathfrak{S},\, \mathfrak{r})|_{L^\infty_v(\R^3)} \leq C \Big\{\|e^{\s_0 x}\nu^{-1}w\mathfrak{S}\|_{L^\infty_{x,v}} +  |w \mathfrak{r}|_{L^\infty_v(\R^3_+)}  \Big\},\\
 \|e^{\s x}wf\|_{L^\infty_{x,v}}+|wf(0)|_{L^\infty_v(\R^3)} \leq \frac{C}{\s_0-\s}\Big\{\|e^{\s_0 x}\nu^{-1}w\mathfrak{S}\|_{L^\infty_{x,v}} +  |w \mathfrak{r}|_{L^\infty_v(\R^3_+)}  \Big\}.
\end{split}
\end{align}
Moreover, if
$\mathfrak{S}$ is continuous in $\R_+\times \mathbb{R}^3$ and $\mathfrak{R}(v)$ is continuous in $\{v\in\R^3_+\}$, then $\mathfrak{f}(x,v)$ is continuous away from the  grazing set $\big\{(0,v)\ : \ v\in\R^3,\  v_3=0\big\}$.
\end{corollary}

\begin{remark}
Corollary \ref{cor1.3} means that the nonlinear boundary layer problem \eqref{1.7} are solvable if the boundary perturbation $\mathfrak{R}$ is on the locally continuous manifold $\mathbb{M}:=\{\mathfrak{R}\, : \, \mathfrak{R}=-(I-P_\g)\mathbb{G}(\mathfrak{S},\mathfrak{r})+ \mathfrak{r},\,\, \mbox{with}\,\, \|e^{\s_0 x}\nu^{-1} w \mathfrak{S}\|_{L^\infty_{x,v}}+\| w\mathfrak{r}\|_{L^\infty_v(\R^3_+)}\leq \delta_0\, \}$.
\end{remark}

\

To study the well-posedness of  nonlinear boundary layer problem, we need first to consider the existence of solution for the following linearized boundary layer problem with a source term
\begin{align}\label{1.7-2}
	\begin{cases}
		\dis v_3 \partial_x f+p_E^0 {\bf L} f=g,\quad 	(x,v)\in \R_+\times \R^3,\\[2mm]
		\dis f(0,v)|_{v_3>0}=P_\g f(0,v)+r(v),\\[2mm]
		\dis \lim_{x\rightarrow\infty}f(x,v)=0.
	\end{cases}
\end{align}

\begin{theorem}\label{thm3.1}
Recall the weight function $w(v)$ in \eqref{wt}, and let $\beta\geq 3$ and  $0\leq \varpi<\f18$. Let $\s_0\in(0,1)$ be suitably small. Assume $g\in \mathcal{N}^{\perp}$ and $\int_{v_3>0} v_3 \sqrt{\mu(v)} r(v) dv=0$ with
\begin{equation}\label{2.20}
\|e^{\s_0 x}\nu^{-1} wg\|_{L^\infty_{x,v}} +  |wr|_{L^\infty_v(\R^3_+)}<+\infty.
\end{equation}
Then there exists a unique function $f^\infty(v):=\mathcal{G}(g,r)(v)$ with
\begin{align}\label{2.21}
	\begin{split}
		&f^\infty(v)=\Big\{ b_{1}^\infty\, v_1+   b_{2}^\infty \,v_2 + c^\infty \,  (\f{|v|^2}{2}-\f32)\Big\}\sqrt{\mu(v)},\\
		&|(b_1^\infty,b_2^\infty, c^\infty)|\leq C\Big\{ \|e^{\s_0 x}\nu^{-1} wg\|_{L^\infty_{x,v}} +  |wr|_{L^\infty_v(\R^3_+)}  \Big\},
	\end{split}
\end{align}
such that the following linearized boundary layer problem of Boltzmann equation
\begin{align}\label{1.7-3}
\begin{cases}
\dis v_3 \partial_x f+p_E^0 {\bf L} f=g,\quad 	(x,v)\in \R_+\times \R^3,\\[2mm]
\dis f(0,v)|_{v_3>0}=P_\g f(0,v)-(I-P_\g)f^{\infty}+r(v),\\[2mm]
\dis \lim_{x\to \infty}f(x,v)=0,
\end{cases}
\end{align}
has a unique solution $f(x,v)$ satisfying
\begin{equation}\label{1.8-0}
\|e^{\s x}wf\|_{L^\infty_{x,v}}+|wf(0)|_{L^\infty_v(\R^3)}
\leq C\Big\{\frac{1}{\s_0-\s}\|e^{\s_0 x}\nu^{-1}wg\|_{L^\infty_{x,v}} +  |wr|_{L^\infty_v(\R^3_+)}  \Big\}.
\end{equation}
The  uniqueness of $f^\infty(v)$ is under the constraints \eqref{2.21}. Moreover, if
$g$ is continuous in $\R_+\times \mathbb{R}^3$ and $r(v)$ is continuous in $\{v\in\R^3_+\}$, then $f(x,v)$ is continuous away from the  grazing set $\big\{(0,v)\ : \ v\in\R^3,\  v_3=0\big\}$.
\end{theorem}

\begin{remark}
Let $g_i\in\mathcal{N}^{\perp}$ and $\int_{v_3>0} v_3 \sqrt{\mu(v)} r_i(v) dv=0$, $i=1,2$ satisfying \eqref{2.20}. Let $f_i, f^\infty_i$ be the solutions obtain in \eqref{2.21}-\eqref{1.8-0} by replacing $g,\, r$ by $g_i,\, r_i$, and  denote
\begin{equation}
	f_i^\infty(v)= \mathcal{G}(g_i,r_i)=\Big\{ b_{i,1}^\infty\, v_1+   b_{i,2}^\infty \,v_2 + c_i^\infty \,  (\f{|v|^2}{2}-\f32)\Big\}\sqrt{\mu(v)}\,\,\, i=1,2,
\end{equation}
then it follows from the uniqueness, $\eqref{2.21}_2$ and \eqref{1.8-0} that
\begin{equation}\label{1.21}
	\mathcal{G}(g_1,\, r_1)+\mathcal{G}(g_2,\, r_2)=\mathcal{G}(g_1+g_2,\, r_1+r_2),\qquad \mathcal{G}(0,0)=0
\end{equation}
and
\begin{align}\label{1.22}
\begin{split}
&|(b_{1,1}^\infty-b_{2,1}^\infty, \, b_{1,2}^\infty-b_{2,2}^\infty, \, c_{1}^\infty-c_{2}^\infty)| \leq C\Big\{\|e^{\s_0 x}\nu^{-1}w(g_1-g_2)\|_{L^\infty_{x,v}} +  |w(r_1-r_2)|_{L^\infty_v(\R^3_+)}  \Big\},\\[2mm]
&|w(f_1-f_2)(0)|_{L^\infty_v(\R^3)}+
\|e^{\s x}w(f_1-f_2)\|_{L^\infty_{x,v}} \\
&\hspace{3.5cm}\leq C\Big\{\frac{1}{\s_0-\s}\|e^{\s_0 x}\nu^{-1}w(g_1-g_2)\|_{L^\infty_{x,v}} +  |w(r_1-r_2)|_{L^\infty_v(\R^3_+)}  \Big\}.
\end{split}
\end{align}
That means  $f, f^\infty$ depend continuously on $g,\, r$ in the sense of \eqref{1.22}.
\end{remark}

\begin{remark}
Theorem \ref{thm3.1} can be used to determine the boundary conditions of higher order viscous boundary layers and Knudsen boundary layers in the Hilbert expansion of Boltzmann equation (compressible Euler scaling) with diffusive reflection boundary conditions.
\end{remark}

\smallskip

\begin{corollary}\label{cor1.4}
 Assume $g, r$ satisfy the conditions in Theorem \ref{thm3.1}, if there exists a function $\mathfrak{r}(v)$ with $w\mathfrak{r} \in L^\infty_v(\R^3_+)$ such that
 \begin{align}\label{1.23}
 r(v)&=-(I-P_\g) \mathcal{G}(g,\mathfrak{r})(v)+\mathfrak{r}(v),\quad \forall \,\, v\in \R^3_+,
 \end{align}
then there is a unique mild solution $f(x,v)$ to the linearized steady boundary layer problem \eqref{1.7-2}, and satisfies
\begin{equation}\label{1.24}
\begin{split}
 |\mathcal{G}(g,\mathfrak{r})|_{L^\infty_v(\R^3)} \leq C\Big\{\|e^{\s_0 x}\nu^{-1}wg\|_{L^\infty_{x,v}} +  |w\mathfrak{r}|_{L^\infty_v(\R^3_+)}  \Big\},\\
\|e^{\s x}wf\|_{L^\infty_{x,v}}+|wf(0)|_{L^\infty_v(\R^3)}
\leq C\Big\{\frac{1}{\s_0-\s}\|e^{\s_0 x}\nu^{-1}wg\|_{L^\infty_{x,v}} +  |w\mathfrak{r}|_{L^\infty_v(\R^3_+)}  \Big\}.
\end{split}
\end{equation}
Moreover, if
$g$ is continuous in $\R_+\times \mathbb{R}^3$ and $r(v)$ is continuous in $\{v\in\R^3_+\}$, then $f(x,v)$ is continuous away from the  grazing set $\big\{(0,v)\ : \ v\in\R^3,\  v_3=0\big\}$.
\end{corollary}

\begin{remark}
Corollary \ref{cor1.4} means that the steady boundary value problem \eqref{1.7-2} are solvable if the boundary perturbation $r$ is on the manifold $\mathcal{M}:=\{r\, : \, r=-(I-P_\g)G(g,\mathfrak{r})+ \mathfrak{r} \,\, \mbox{with}\,\, e^{\s_0 x}\nu^{-1} w g \in L^\infty_{x,v},\, w\mathfrak{r} \in L^\infty_v(\R^3_+)\, \}$.
\end{remark}

We now briefly comment on the analysis of the present paper. The main part of this paper is to prove Theorem \ref{thm3.1}, that is to prove the existence of solution for the linear boundary layer problem of Boltzmann equation with diffusive reflection boundary conditions. We start with the construction of approximate solutions of the truncated problem \eqref{S3.3} with penalization. Firstly we establish {\it a priori} uniform $L^\infty_{x,v}$-estimate which is independent of both truncation parameter $d\geq $ and penalization parameter $\v\in[0,1]$, see Lemma \ref{lemS3.3} for details. With the help of the uniform $L^\infty_{x,v}$-estimate, and by similar arguments as in \cite{DHWZ}, we can obtain the solution $f$ of  approximate problem \eqref{S3.106} with
\begin{align}\label{1.29}
\|wf\|_{L^\infty_{x,v}} +|wf|_{L^\infty(\gamma)} \leq C_d \{\|\nu^{-1}wg\|_{L^\infty_{x,v}} + |wr|_{L^\infty(\g_-)}\},
\end{align}
see Lemmas \ref{lemS3.4}-\ref{lemS3.7} for details.

We note that the bound on the right hand side (RHS) of \eqref{1.29} depends on $d$. However, to construct the boundary layer solution of Boltzmann equation in half-space, we need some  uniform   estimates independent of $d$ so that we can take the limit $d\to+\infty$. We define
\begin{equation*}
	\bar{f}(x,v):=f(x,v)-\sqrt{2\pi \mu(v)}\, z_{\g_+}(f),
\end{equation*}
where  $\dis z_{\g_+}(f):=\int_{v_3<0} |v_3| \sqrt{\mu(v)} f(0,v) dv$. Clearly, $z_{\g_{+}}(\bar{f})=\int_{v_3<0} |v_3| \sqrt{\mu(v)} \bar{f}(0,v) dv=0$ and the equation of $\bar{f}$ has the same form as $f$, see \eqref{S3.123-5}.
Then, by energy estimate, we can prove
\begin{align}\label{1.31}
	&|(I-P_\g)\bar{f}(0)|^2_{L^2(\g_+)} + \int_0^d e^{2\sigma_1 x} \|({\bf I-P})\bar{f}(x)\|_{\nu}^2dx\nonumber\\
	&\leq C\Big\{|r|^2_{L^2(\g_-)} + \int_0^d e^{2\sigma_1 x} \|g(x)\|^2_{L^2_{x,v}} dx \Big\}\quad\mbox{for}\quad \sigma_1\in[0,\s_0].
\end{align}
It is very important that the term on RHS of \eqref{1.31} is independent of $d$. To obtain the uniform in $d$ estimate for macroscopic term, motivated by \cite{GPS-1988}, we define
\begin{equation*}
\tilde{f}(x,v):=\bar{f}(x,v)+\Phi(v)
\end{equation*}
with
\begin{equation*}
\Phi(v):=\big[\phi_0+\phi_1 v_1+\phi_2v_2+\phi_3(\frac{|v|^2}{2}-\f32)\big] \sqrt{\mu(v)},
\end{equation*}
where $(\phi_0, \phi_1, \phi_2,\phi_3)(d)$ are four constants determined in Lemma \ref{lem2.12}. And $\tilde{f}$ satisfies
\begin{align*}
	\begin{cases}
		v_3 \partial_x\tilde{f}+p_E^0 \FL \tilde{f}=g,\quad (x,v)\in \Omega_d\times \R^3, \\
		\tilde{f}(0,v)|_{v_3>0}=P_\g \tilde{f}(0,v) + (I-P_\g)\Phi + r(v),\\
		\tilde{f}(d,v)|_{v_3<0}=\tilde{f}(x,\Mr v).
	\end{cases}
\end{align*}
Then we can get the uniform $L^2_{x,v}$-estimate for $\tilde{f}$, i.e.,
\begin{equation}\label{3.141}
	\|e^{\s x}\tilde{f}\|_{L^2_{x,v}}\leq C\left\{|r|_{L^2(\g_-)}+\frac{1}{\s_1-\s}\|e^{\s_1 x} g\|_{L^2_{x,v}}\right\},
\end{equation}
with $0<\s<\s_1\leq \s_0$, and the constant $C>0$ is independent of $d$, see Lemma \ref{lem2.13}. Fortunately, the right hand side of \eqref{3.141} is independent of $\Phi$. Applying the $L^2_{x,v}$-$L^\infty_{x,v}$ estimate to $e^{\s x} \tilde{f}$,  we can prove
\begin{align}\label{3.166}
	&\|e^{\s x}w\tilde{f}\|_{L^\infty_{x,v}}+|e^{\s x}w\tilde{f}|_{L^\infty(\g)}\nonumber\\
	&\leq C\Big\{\frac{1}{\s_0-\s}\|e^{\s_0 x}\nu^{-1}wg\|_{L^\infty_{x,v}} +  |w\big((I-P_\g)\Phi + r\big)|_{L^\infty(\g_-)}  \Big\},
\end{align}
see Lemma \ref{lem2.16} for details.
From the proof of Lemma \ref{lem2.12}, we know that the constants $(\phi_0,\phi_1,\phi_2, \phi_3)(d)$ should depend on $d$, hence the RHS of \eqref{3.166} still depends on $d$. Hence, to obtain the uniform estimate for $\tilde{f}$, we have to establish the uniform in $d$ estimate  for $(\phi_0,\phi_1,\phi_2, \phi_3)(d)$. This is the key part of the present paper. In fact, we express the macroscopic part  $\mathbf{P} \bar{f}(x,v)=\big[\bar{a}(x)+\bar{b}(x)\cdot v+\bar{c}(x) (\f{|v|^2}{2}-\f32)\big] \sqrt{\mu}$ by using the boundary value (see Lemma \ref{lem2.14}) and get
\begin{align*}%\label{3.182}
|(\bar{a},\bar{b}_{1},\bar{b}_{2},\bar{c})(d)|
&\leq C \Big\{ |(\mathbf{I-P})\tilde{f}_d(d)|_{L^\infty_v} + \|g\|_{L^2_{x,v}} + |(I-P_\g)\bar{f}(0)|_{L^2(\g_+)} + |r|_{L^2(\g_-)}\Big\},
%&\leq C \Big\{ |(\mathbf{I-P})\tilde{f}_d(d)|_{L^\infty_v} + \|g\|_{L^2_{x,v}} + |r|_{L^2(\g_-)}\Big\}.
\end{align*}
which, together with \eqref{S3.140-1}, \eqref{1.31} and \eqref{3.166}, yields that
\begin{align*}%\label{3.183}
	&|(\phi_0,\phi_1,\phi_2,\phi_3)(d)|
	%\leq C |(\bar{a},\bar{b}_{1},\bar{b}_{2},\bar{c})(d)|+ C |(\mathbf{I-P})\tilde{f}(d)|_{L^\infty_v}\nonumber\\
	%&\leq C|(\mathbf{I-P})\tilde{f}_d(d)|_{L^\infty_v} + C \Big\{\|g\|_{L^2_{x,v}}  + |r|_{L^2(\g_-)}\Big\}\nonumber\\
	\leq Ce^{-\sigma d} |e^{\sigma d} w\tilde{f}_d(d)|_{L^\infty_{v}} + C \Big\{\|g\|_{L^2_{x,v}}  + |r|_{L^2(\g_-)}\Big\}\nonumber\\
	&\leq Ce^{-\sigma d} |(\phi_0,\phi_1,\phi_2,\phi_3)(d)| + C\Big\{\frac{1}{\s_0-\s}\|e^{\s_0 x}\nu^{-1}wg\|_{L^\infty_{x,v}} +  |wr|_{L^\infty(\g_-)}  \Big\}.
\end{align*}
Then we have
\begin{align}\label{3.180}
	|(\phi_0,\phi_1,\phi_2,\phi_3)(d)|\leq C\Big\{\|e^{\s_0 x}\nu^{-1}wg\|_{L^\infty_{x,v}} + |wr|_{L^\infty(\g_-)}  \Big\},
\end{align}
and hence
\begin{align}\label{3.181}
	\|e^{\s x}w\tilde{f}\|_{L^\infty_{x,v}}+|e^{\s x}w\tilde{f}|_{L^\infty(\g)}
	\leq C\Big\{\frac{1}{\s_0-\s}\|e^{\s_0 x}\nu^{-1}wg\|_{L^\infty_{x,v}} +  |wr|_{L^\infty(\g_-)}  \Big\},
\end{align}
where the constants $C>0$ are independent of $d$, see Lemma \ref{lem2.17} for details. Therefore we have established the uniform in $d$ estimates for both $\tilde{f}$ and $(\phi_0,\phi_1,\phi_2,\phi_3)(d)$.

To take the limit $d\to +\infty$, we still need to obtain the asymptotic behavior of $(\phi_0,\phi_1,\phi_2,\phi_3)(d)$. Using \eqref{3.181} and energy estimate, we can prove
	\begin{align}\label{3.158}
	|({I-P_\g})(\bar{f}_{d_2}-\bar{f}_{d_1})(0)|_{L^2(\g_+)}
	\leq  Ce^{-\sigma_1 d_1 }\Big\{|r|^2_{L^2(\g_-)} + \int_0^{d_1} e^{2\sigma_1 x} \|g(x)\|^2_{L^2_{x,v}} dx \Big\}^{\f12}.
\end{align}
where we have denoted $\bar{f}$ to be  $\bar{f}_{d}$ to emphasize the dependent on $d$, see Lemma \ref{lem2.15} for details. Then using \eqref{S3.140-1},  \eqref{3.158}, \eqref{3.181} and  Lemma \ref{lem2.14}, we can show that there exist constants  $(\phi_0^\infty,\phi_1^\infty,\phi_2^\infty, \phi_3^\infty)$ such that $\lim_{d\to +\infty}(\phi_0,\phi_1,\phi_2, \phi_3)(d)=(\phi_0^\infty,\phi_1^\infty,\phi_2^\infty, \phi_3^\infty)$  and
\begin{equation*}
	|(\phi_0^\infty,\phi_1^\infty,\phi_2^\infty, \phi_3^\infty)|\leq C\Big\{\|e^{\s_0 x} \nu^{-1}wg\|_{L^\infty_{x,v}} +  |wr|_{L^\infty(\g_-)}  \Big\}.
\end{equation*}
With above uniform estimates, we can finally  prove Theorem \ref{thm3.1}.

\smallskip

The paper is organized as follows: In section 2, we prove Theorem \ref{thm3.1} by a series of approximations and  some uniform estimates.
In section 3, we prove Theorem \ref{thm1.1}. Some useful known results are presented in Appendix \ref{APP}.

\smallskip

We list some notations that will be used in this paper.
Throughout this paper, $C$ denotes a generic positive constant which may vary from line to line.
And $C_a$, $C_b, \cdots$ denote the generic positive constants depending on $a, b, \cdots$, respectively, which also may vary from line to line. We denote  $\Omega=\R_+\ \mbox{or} \  (0,d)$, and the phase boundary  $\gamma=\partial \Omega \times  \mathbb{R}_v^3 $. We denote  $\| \cdot \|_{L^p}$  the standard $L^p(\Omega \times \mathbb{R}_v^3)$-norm or $L^p(\mathbb{R}_v^3)$-norm, and $\|\cdot \|_{\nu}=\|\sqrt{\nu}\cdot \|_{L^2}$. When the norms need to be distinguished from each other, we write $\|\cdot\|_{L^p_v}$,
$\|\cdot\|_{L^p_x}$ and $\|\cdot\|_{L^p_{x,v}}$, respectively.     For the phase boundary integration, we define $|\cdot|_{L^\infty(\gamma)}$ denotes the $L^\infty(\gamma)$-norm,
$|\cdot|_{L^2(\gamma)}$ denotes the $L^2(\gamma, |v_3|dv)$-norm.

%%%%%%%%%%%%%%%%%%%%%%%%%%%%%%%%%%%%%%%%%%%%%%%%%%%%%%%%%%%%%%%%%%%%%%%%%%%%%%%%%%%%%%%%%%%%%%%%%%%%%%%%%%%%%%%%%%%%%%%%%%%%%%%%%%%%%%%%%%%%%%%%%%%%%%%%%%%

\section{Existence for the Linearized problem }\label{L}

 We denote $\Omega_d:=(0,d)$ with $d\geq1$ and denote the phase boundary of  $\Omega_d \times \mathbb{R}^{3}$ as
$\gamma =\partial \Omega_d \times \mathbb{R}^{3}$. We split $\gamma$ into three disjoint parts, outgoing
boundary $\gamma _{+},$ the incoming boundary $\gamma _{-},$ and the
singular boundary $\gamma _{0}$ for grazing velocities$:$
\begin{align}
\gamma _{+} &=\{(x,v)\in \partial \Omega_d \times \mathbb{R}^{3}:
n(x)\cdot v>0\}, \nonumber\\
\gamma _{-} &=\{(x,v)\in \partial \Omega_d \times \mathbb{R}^{3}:
n(x)\cdot v<0\}, \nonumber\\
\gamma _{0} &=\{(x,v)\in \partial \Omega_d \times \mathbb{R}^{3}:
n(x)\cdot v=0\}.\nonumber
\end{align}
where $\vec{n}(x)$ is the outward unit normal. It is direct to know that $\partial\Omega_d=\{0,d\}$, $\vec{n}(0)=(0,0,-1)$ and $\vec{n}(d)=(0,0,1)$. %For the phase boundary integration, we define $|\cdot|_{L^\infty(\gamma)}$ denotes the $L^\infty(\gamma)$-norm, $|\cdot|_{L^2(\gamma)}$ denotes the $L^2(\gamma, |v_3|dv)$-norm, where $\gamma=\partial \Omega \times  \mathbb{R}_v^3 $ with $\Omega=\R_+\ \mbox{or} \  \Omega_d$.

\smallskip

%For later use we define the notation
%\begin{align}
%(P_\gamma f)(x,v):=\sqrt{2\pi \mu(v)} \int_{u\cdot \vec{n}>0} f(x,u)\sqrt{\mu} |u_3| du.
%\end{align}
To construct a solution of the linearized boundary layer problem in half-space,  we first consider the truncated  problem with penalized term
\begin{align}\label{S3.3}
\begin{cases}
\dis \v f^\v+ v_3 \cdot \partial_x f^\v+p_E^0 \FL f^\v=g,\quad (x,v)\in \Omega_d\times \R^3,\\[2mm]
\dis f^\v(0,v)|_{v_3>0}=(P_{\g} f)(0,v)+r(v),\\[2mm]
\dis f^{\v}(d,v)|_{v_3<0}=f^{\v}(d,\mathcal{R}v),
\end{cases}
\end{align}
where $\v\in(0,1]$ and $\Mr v:=(v_{\sp}, -v_3)$ with $v_{\sp}:=(v_1,v_2)$. We define
\begin{equation*}%\label{S3.7}
h^\v(x,v):=w(v) f^\v(x,v), %.
\end{equation*}
then %Then
\eqref{S3.3} can be rewritten as
\begin{align}\label{S3.8}
\begin{cases}
\dis \v h^\v+v_3\cdot \partial_x h^\v+p_E^0\nu(v) h^\v=p_E^0K_{w} h^\v+wg,\\[2mm]
\dis h^\v(0,v)|_{v_3>0}=\f{1}{\tilde{w}(v)} \int_{u_3<0} h^\v(0,u) \tilde{w}(u) d\sigma+(wr)(v),\\[2mm]
h^\v(d,v)|_{v_3<0}=h^\v(d, \Mr v),
\end{cases}
\end{align}
where
\begin{equation}\label{2.7-1}
\tilde{w}(v)=\frac{1}{\sqrt{2\pi \mu(v)} w(v)}\quad \mbox{and} \quad d\sigma=\sqrt{2\pi} \mu(v) |v_3| dv.
\end{equation}
It is easy to check that
\begin{equation*}
\int_{v_3<0} d\sigma=\int_{v_3<0}  \sqrt{2\pi} \mu(v) |v_3| dv=1.
\end{equation*}
Hereafter $K_wh=wK(\f{h}{w})$ and
\begin{align}\label{3.8}
K_wh(v)=\int_{\R^3} k_w(v,u) h(u) du\quad\mbox{with}\quad
k_w(v,u)=w(v) k(v,u) w(u)^{-1}.
\end{align}

%\begin{remark}
%Via the approach developed in \cite{Guo2},  our proof is organized as follows: we first establish uniform estimates in $\v, n$ for the solution to the truncated problem and pass to the limit as $\v \to 0, n \to \infty$. Then, we prove uniform estimates in $d$ and pass to the limit as $d \to \infty$.  The estimates \eqref{3.132-0}, \eqref{3.132-01} and the uniqueness follows consequently.
%\end{remark}

\subsection{A priori $L^\infty_{x,v}$ estimate}
For the approximate problem  \eqref{S3.8},  the most difficult part is to obtain the {\it a priori} $L^\infty_{x,v}$-bound uniform in $\v\in[0,1]$ and $ d\in[1,\infty)$.

\begin{definition}
Given $(t,x,v),$ let $[X(s),V(s)]$ %=[X(s;t,x,v),V(s;t,x,v)]=[x+(s-t)v,v]$
be the backward characteristics for \eqref{S3.8}, which is determined by
\begin{align*}%\label{ode}
\begin{cases}
\dis \frac{{\dd}X(s)}{{\dd}s}=V_3(s),
\quad \frac{{\dd}V(s)}{{\dd}s}=0,\\[2mm]
[X(t),V(t)]=[x,v].
\end{cases}
\end{align*}
The solution is then given by
\begin{equation*}
%\label{def.bic}
[X(s),V(s)]=[X(s;t,x,v),V(s;t,x,v)]=[x-(t-s)v_3,v].
\end{equation*}
\end{definition}
Now for each $(x,v)$ with $x\in \bar{\Omega}_d$ and $v_3\neq 0,$ we define its {backward exit time} $t_{\mathbf{b}}(x,v)\geq 0$ to be the last moment at which the
back-time straight line $[X({-\tau};0,x,v),V({-\tau};0,x,v)]$ remains in $\bar{\Omega}_d$:
\begin{equation*}%\label{exit}
t_{\mathbf{b}}(x,v)={\sup\{s \geq 0:x-\tau v_3\in\bar{\Omega}_d\text{ for }0\leq \tau\leq s\}.}
\end{equation*}
We also define
\begin{equation*}%\label{xb}
x_{\mathbf{b}}(x,v)=x(t_{\mathbf{b}})=x-t_{\mathbf{b}}(x,v)\, v_3\in \partial \Omega_d .
\end{equation*}
We point out that $X(s)$, $t_{\mathbf{b}}(x,v)$ and $x_{\mathbf{b}}(x,v)$ are independent of  the horizontal velocity.

\vspace{1mm}

Let $x\in \bar{\Omega}_d$, $(x,v)\notin \gamma _{0}\cup \g_{-}$ and denote
$
(t_{0},x_{0},v_{0}):=(t,x,v)$ hereafter. We firstly define
\begin{equation}
v_{k}\in \mathcal{V}_{k}:=\{ v_{k}\in \R^3 \,:\, v_{k}\cdot \vec{n}(0) =-v_{k,3}>0\},\quad\mbox{for}\,\, k\geq1,
\end{equation}
which is defined only at the lower boundary $x=0$.
Then the back-time cycle is defined as
\begin{align}\label{2.7}
\begin{cases}
\dis X_{cl}(s;t,x,v)&=\sum_{k}\Fi_{[t_{k+1},t_{k})}(s)\{x_{k}-v_{k,3}\cdot(t_{k}-s)\},\\[2mm]
\dis V_{cl}(s;t,x,v)&=\sum_{k}\Fi_{[t_{k+1},t_{k})}(s)v_{k},
\end{cases}
\end{align}
with
\begin{equation*}
(t_{k+1},x_{k+1},v_{k+1})=(t_{k}-t_{\mathbf{b}}(x_{k},v_{k}),x_{\mathbf{b}}(x_{k},v_{k}), v_{k+1} ),
\end{equation*}
where
\begin{equation}
\begin{split}
v_{k+1}=\Mr v_{k}=(v_{\sp}, -v_{k,3})\quad \mbox{if}\quad x_{k+1}=d,\\
v_{k+1}\in \mathcal{V}_{k+1}\quad \mbox{if} \quad x_{k+1}=0.
\end{split}
\end{equation}

%It is direct to check that the velocity space $\mathcal{V}_l$ implicitly depends on $\{x, v,  v_1,v_2,\cdots, v_{l-1}\}$.
 Clearly, for $k\geq 1$ and $(x,v)\notin \g_0\cup \g_-$, it holds that
\begin{equation}\label{2.8}
\begin{split}
&x_k=\frac{1-(-1)^k}{2} x_1+\frac{1+(-1)^k}{2} x_2,\\
&t_k-t_{k+1}=\frac{d}{|v_{k,3}|}>0.
\end{split}
\end{equation}

\

We also define the iterated integral
\begin{equation}\label{2.8-1}
\int_{\Pi_{j=1}^{k-1}\mathcal{V}_j}  \Pi_{j=1}^{k-1} d\sigma_j=\int_{\mathcal{V}_1} d\sigma_1\int_{\mathcal{V}_2} d\sigma_2\cdots\int_{\mathcal{V}_{k-1}} d\sigma_{k-1},
\end{equation}
where
\begin{equation*}
\dis d\sigma_j=\sqrt{2\pi} \mu(v_j) |v_{j,3}| dv_j,\quad j=1,2,\cdots, k-1.
\end{equation*}
It is direct to know that $d\sigma_j, \, j=1,2,\cdots, k-1$ are probability measure.

\smallskip

The following lemma will be used in the proof of Lemma \ref{lemS3.3}, and the proof  is very similar to the one in Guo \cite{Guo2}, so we omit the proof here.
\begin{lemma}\label{lem2.1}
%For any given $\epsilon>0$, there exists $k_0(\epsilon,T_0)\ll1 $ such that for $k\geq k_0$,
Let $(t,x,v)$ with $t\in[0,T_0]$, $x\in\bar{\Omega}_d,\, v\in \R^3$ and $v_3>0$, %it holds that
%\begin{equation}
%\int_{\Pi_{l=1}^{k-1}\mathcal{V}_l} I_{t_k(t,x,v,v_1,\cdots, v_{k-1})>0} \Pi_{l=1}^{k-1} d\sigma_l \leq \epsilon.
%\end{equation}
there exist constants $C_1, C_2>0$, independent of $T_0$, such that  for $ k\geq k_0:=C_1 T_0^{\f43}\gg1$, it holds
\begin{equation}
\int_{\Pi_{l=1}^{k-1}\mathcal{V}_{2l-1}} \mathbf{1}_{\{t_{2(k-1)+1}(t,x,v,v_1,\cdots, v_{2(k-1)-1})>0\}} \Pi_{l=1}^{k-1} d\sigma_{2l-1}  \leq \left(\f12\right)^{C_2 T_0^{\f43}}.
\end{equation}
We point out that all the constants above are independent of $d\geq 1$.
\end{lemma}

We can represent the solution of \eqref{S3.8} in a mild formulation which enables us to get the $L^\infty_{x,v}$ bound of solutions.  Indeed, for later use, we consider the following iterative linear problems involving a parameter $\la\in [0,1]$:
\begin{equation}\label{2.9}
\begin{cases}
\v h^{i+1}+v_3\partial_x h^{i+1}+p_E^0 \nu(v) h^{i+1}=\lambda p_E^0 K_w h^i +wg,\\[3mm]
\dis h^{i+1}(0,v)|_{v_3>0}=\f{1}{\tilde{w}(v)} \int_{v'_3<0} h^i (0,v') \tilde{w}(v') d\sigma'+(wr)(v),\\[3mm]
\dis h^{i+1}(d,v)|_{v_3<0}=h^i(d,\Mr v) +(w\tilde{r})(v)
\end{cases}
\end{equation}
for $i=0,1, 2,\cdots$, where $h^0\equiv0$ and $(wr)(v), \, (w\tilde{r})(v)\in L^\infty(\gamma)$ are some given functions. For the mild formulation of \eqref{2.9}, we have the following lemma whose proof can be given by induction arguments, but omitted for brevity here.

\begin{lemma}\label{lem2.2}
Let $\v\in[0,1]$ and $\lambda\in[0,1]$. For each $(x,v)\in \bar{\Omega}_d\times \mathbb{R}^3\setminus (\gamma_0\cup \gamma_-)$,
we have:

\smallskip
\noindent{\it Case 1.} If $v_{0,3}>0$, it holds
\begin{align}\label{2.10}
h^{i+1}(x,v)&=\sum_{n=1}^3 J_n+\sum_{n=4}^{16} I_{\{t_1>0\}} J_n,%\quad\mbox{if} \quad v_3>0,
\end{align}
with $\nu_{\v}(v):=\v+p_E^0 \nu(v)$, and
\begin{align*}
J_1&=I_{\{t_1\leq 0\}} \, e^{-\nu_{\v}(v) t}\, h^{i+1}(x-v_{0,3}t,v),\\
J_2+J_3&= \int_{\max\{t_1,0\}}^t e^{-\nu_{\v}(v) (t-s)} \Big(\lambda K_w h^i+wg\Big)(x-v_{0,3}(t-s),v) ds,
\end{align*}
\begin{align*}
J_4&=e^{-\nu_{\v}(v) (t-t_1)} (wr)(x_1,v),\\
J_5&=\frac{e^{-\nu_{\v}(v) (t-t_1)}}{\tilde{w}(v)} \Big\{\int_{\Pi_{j=1}^{k-1}\mathcal{V}_{2j-1}}  \sum_{l=1}^{k-2} \Fi_{\{t_{2l+1}>0\}} (wr)(v_{2l-1}) d\Sigma^{k-1}_{l}(t_{2l+1})\\
&\qquad\qquad\qquad+\int_{\Pi_{j=1}^{k-1}\mathcal{V}_{2j-1}}  \sum_{l=1}^{k-1} \Fi_{\{t_{2l}>0\}} (w\tilde{r})(v_{2l-1}) d\Sigma^{k-1}_{l}(t_{2l})\Big\},
\end{align*}
and
\begin{align*}
J_6&=\frac{e^{-\nu_{\v}(v) (t-t_1)}}{\tilde{w}(v)}  \int_{\Pi_{j=1}^{k-1}\mathcal{V}_{2j-1}}  \sum_{l=1}^{k-1} \Fi_{\{t_{2l}\leq 0<t_{2l-1}\}} h^{i+2-2l}(x_{2l-1}-v_{2l-1,3} t_{2l-1},v_{2l-1}) d\Sigma^{k-1}_{l}(0),\\
J_7+J_8&=\frac{e^{-\nu_{\v}(v) (t-t_1)}}{\tilde{w}(v)} \int_{\Pi_{j=1}^{k-1}\mathcal{V}_{2j-1}}  \sum_{l=1}^{k-1} \Fi_{\{t_{2l}\leq 0<t_{2l-1}\}} \int_0^{t_{2l-1}} \nonumber\\
&\qquad\qquad\qquad\quad\times (\lambda K_w h^{i+1-2l}+wg)(x_{2l-1}-v_{2l-1,3} (t_{2l-1}-s),v_{2l-1}) ds d\Sigma^{k-1}_{l}(s),\\
J_9+J_{10}&=\frac{e^{-\nu_{\v}(v) (t-t_1)}}{\tilde{w}(v)} \int_{\Pi_{j=1}^{k-1}\mathcal{V}_{2j-1}}  \sum_{l=1}^{k-1} \Fi_{\{t_{2l}>0\}} \int_{t_{2l}}^{t_{2l-1}}\nonumber\\
&\qquad\qquad\qquad\quad\times (\lambda K_w h^{i+1-2l}+wg)(x_{2l-1}-v_{2l-1,3} (t_{2l-1}-s),v_{2l-1}) ds d\Sigma^{k-1}_{l}(s),\\
J_{11}&=\frac{e^{-\nu_{\v}(v) (t-t_1)}}{\tilde{w}(v)}  \int_{\Pi_{j=1}^{k-1}\mathcal{V}_{2j-1}}  \sum_{l=1}^{k-1} \Fi_{\{t_{2l+1}\leq 0<t_{2l}\}} h^{i+1-2l}(x_{2l}+v_{2l-1,3} t_{2l},\Mr v_{2l-1}) d\Sigma^{k-1}_{l}(0),\\
J_{12}+J_{13}&=\frac{e^{-\nu_{\v}(v) (t-t_1)}}{\tilde{w}(v)} \int_{\Pi_{j=1}^{k-1}\mathcal{V}_{2j-1}}  \sum_{l=1}^{k-1} \Fi_{\{t_{2l+1}\leq 0<t_{2l}\}} \int_{0}^{t_{2l}}\nonumber\\
&\qquad\qquad\qquad\quad\times (\lambda K_w h^{i-2l}+wg)(x_{2l}+v_{2l-1,3} (t_{2l}-s),\Mr v_{2l-1}) ds d\Sigma^{k-1}_{l}(s),\\
J_{14}+J_{15}&=\frac{e^{-\nu_{\v}(v) (t-t_1)}}{\tilde{w}(v)} \int_{\Pi_{j=1}^{k-1}\mathcal{V}_{2j-1}}  \sum_{l=1}^{k-1} \Fi_{\{t_{2l+1}>0\}} \int_{t_{2l+1}}^{t_{2l}}\nonumber\\
&\qquad\qquad\qquad\quad\times (\lambda K_w h^{i-2l}+wg)(x_{2l}+v_{2l-1,3} (t_{2l}-s),\Mr v_{2l-1}) ds d\Sigma^{k-1}_{l}(s),\\
J_{16}&=\frac{e^{-\nu_{\v}(v) (t-t_1)}}{\tilde{w}(v)} \int_{\Pi_{j=1}^{k-1}\mathcal{V}_{2j-1}} \Fi_{\{t_{2(k-1)+1}>0\}} h^{i+1-2(k-1)}(x_{2(k-1)+1},\Mr v_{2(k-1)-1}) d\Sigma^{k-1}_{k-1} (t_{2(k-1)+1}),
\end{align*}
where we have used the notation
\begin{align}\label{2.11}
d\Sigma_l^{k-1}(s):= \{\Pi_{j=l+1}^{k-1} d\sigma_{2j-1}\}&\cdot \{\tilde{w}(v_{2l-1}) e^{-\nu_{\v}(v_{2l-1}) (t_{2l-1}-s)} d\sigma_{2l-1}\}\nonumber\\
&\times \{\Pi_{j=1}^{l-1} e^{-\nu_{\v}(v_{2j-1}) (t_{2j-1}-t_{2j+1})} d\sigma_{2j-1}\}.
\end{align}

\smallskip

\noindent{\it Case 2.} If $v_{0,3}<0$, it holds
\begin{align}\label{2.12}
h^{i+1}(x,v)&=\Fi_{\{t_1\leq 0\}}\cdot e^{-\nu_{\v}(v) t}\cdot h^{i+1}(x-v_{0,3}t,v)\nonumber\\
&\quad+ \int_{\max\{t_1,0\}}^t e^{-\nu_{\v}(v) (t-s)}\Big(\lambda K_w h^i +w g\Big)(x-v_{0,3}(t-s),v) ds\nonumber\\
&\quad + \Fi_{\{t_1> 0\}}\cdot e^{-\nu_{\v}(v) (t-t_1)}\cdot (w\tilde{r})(v_1) + \Fi_{\{t_1> 0\}}\cdot e^{-\nu_{\v}(v) (t-t_1)}\cdot h^{i}(d,v_1).
\end{align}
where $v_1=\Mr v=(v_{0,\sp},-v_{0,3})$ and $v_{1,3}=-v_{0,3}>0$. Therefore we can express the last term $h^{i}(d,v_1)$ of \eqref{2.12}  by using the formula in {\it Case 1}.

\end{lemma}

%\begin{proof} The proof is very similar as Guo, the details are omitted here for simplicity.
%\end{proof}

\begin{lemma}\label{lemS3.3}
Let $h^i$, $i=0,1,2,\cdots$, be the solutions to \eqref{2.9}, satisfying
$$
\|h^i\|_{L^\infty_{x,v}}+|h^i|_{L^\infty(\gamma_+)}<\infty.
$$ %for $i=0,1,2,\cdots$.
Then there exists $T_0>1$ large enough  such that for $i\geq k_0:=C_1T_0^{\f43}$, it holds %, for $i\geq k$,
that
\begin{align}\label{3.13}
\|h^{i+1}\|_{L^\infty_{x,v}}+|h^{i+1}|_{L^\infty(\g_+)}
&\leq \frac18 \sup_{0\leq l\leq 2k_0} \{\|h^{i-l}\|_{L^\infty_{x,v}}+|h^{i-l}|_{L^\infty(\g_+)}\}\nonumber\\
&\hspace{-28mm} +C\{\|\nu^{-1}wg\|_{L^\infty_{x,v}}+|wr|_{L^\infty(\gamma_-)} +|w\tilde{r}|_{L^\infty(\gamma_-)}\}+C \sup_{0\leq l\leq 2k_0}\left\{\left\|\frac{h^{i-l}}{w}\right\|_{L^2_{x,v}}\right\}.
\end{align}
Moreover, if $h^i\equiv h$ for $i=1,2,\cdots$, i.e., $h$ is a solution, then \eqref{3.13} is reduced to the following estimate
\begin{align}\label{3.14}
\|h\|_{L^\infty_{x,v}}{+|h|_{L^\infty(\g)}}&\leq C\{\|\nu^{-1}wg\|_{L^\infty_{x,v}}+|wr|_{L^\infty(\gamma_-)}+|w\tilde{r}|_{L^\infty(\gamma_-)}\}+C\left\|\frac{h}{w}\right\|_{L^2_{x,v}}.
\end{align}
We emphasize that the positive constant  $C>0$  depends on $k_0$, and is  independent of $d$, $\lambda$ and $\v>0$.
\end{lemma}

\noindent{\bf Proof.} By the definition of $\nu_{\v}(v)$, we first note that
\begin{equation}\label{S3.25}
\nu_{\v}(v)\geq p_E^0 \nu(v)\geq \nu_E^0>0,
\end{equation}
where $\nu_E^0$ is a positive constant independent of $\v$ and $v\in\mathbb{R}^3$. Since the proof is complicate, we divide the proof into four steps.

\smallskip

\noindent{\it Step 1.} In this step, we assume $v_{0,3}>0$. Hence, $h^{i+1}(x,v)$ can be expressed by \eqref{2.10}. Then, for $J_1$, it follows from \eqref{S3.25} that
\begin{equation}\label{S3.26}
|J_1|\leq e^{-\nu_E^0 t} \|h^{i+1}\|_{L^\infty_{x,v}}.
\end{equation}

For   terms involving the source $g$, we notice that
\begin{align}\label{S3.27}
%\tilde{w}(v)=%\frac{1}{\sqrt{2\pi}}
%\sqrt{2\pi}\frac{e^{(\frac14-\varpi)|v|^2}}{(1+|v|^2)^{\frac{\beta}{2}}};\quad
\frac{1}{\tilde{w}(v)}\cong %\sqrt{2\pi}
%\frac{1}{\sqrt{2\pi}}(1+|v|^2)^{\frac{\beta}{2}} e^{-(\frac14-\varpi)|v|^2}\leq C e^{-\frac18|v|^2},
w(v)e^{-\f{|v|^2}{4}}\leq Ce^{-\f{|v|^2}{8}},
\qquad \tilde{w}(v) \mu(v)|v_3|\leq C\frac{\sqrt{\mu(v)}}{w(v)}\leq Ce^{-\f{|v|^2}{8}}
\end{align}
which immediately yields  that
\begin{equation}\label{S3.28}
%\left\{\begin{aligned}
\int_{\Pi _{j=1}^{k-1}\mathcal{V}_{2j-1}} \tilde{w}(v_l)  \  \Pi_{j=1}^{k-1} \dd\sigma_{2j-1}\leq C<\infty,\quad\mbox{for}\quad 1\leq l\leq k-1,%\\[1.5mm]
%&\int_{\Pi _{j=1}^{k-1}\hat{\mathcal{V}}_{j}} \sum_{l=1}^{k-1} \Fi_{\{t_{l+1}\leq0<t_l\}} \tilde{w}(v_l) \, \Pi_{j=1}^{k-1} \dd\sigma_j\leq Ck.
%\end{aligned}\right.
\end{equation}
Then it follows from \eqref{S3.27} and \eqref{S3.28} that
\begin{align}\label{S3.29}
|J_3|+|J_{8}|+|J_{10}|+|J_{13}|+|J_{15}|\leq Ck \| \nu^{-1}wg\|_{L^\infty_{x,v}},
\end{align}
\begin{align}\label{S3.30}
|J_4|+|J_5|\leq Ck  \{|wr|_{L^\infty{(\g_-)}}+|w\tilde{r}|_{L^\infty(\gamma_-)}\},
\end{align}
and
\begin{align}\label{S3.31}
|J_6|&\leq C e^{-\frac18|v|^2} e^{-\nu_E^0 (t-t_1)} \int_{\Pi _{j=1}^{k-1}\mathcal{V}_{2j-1}} \sum_{l=1}^{k-1} \Fi_{\{t_{2l}\leq0<t_{2l-2}\}} \|h^{i+2-2l}\|_{L^\infty_{x,v}} \dd\Sigma_{l}^{k-1}(0)\nonumber\\
&\leq Ck  e^{-\nu_E^0 t} e^{-\frac18|v|^2} \cdot \sup_{1\leq l\leq k-1}\{\|h^{i+2-2l}\|_{L^\infty_{x,v}}\}.
\end{align}
Similarly we have
\begin{align}\label{S3.31-1}
|J_{11}|\leq Ck e^{-\nu_E^0 t} e^{-\frac18|v|^2} \cdot \sup_{1\leq l\leq k-1}\{\|h^{i+1-2l}\|_{L^\infty_{x,v}}\}.
\end{align}

From the boundary condition   $\eqref{2.9}_2$, it holds that
\begin{align}\label{S3.32-1}
|h^{i+1-2(k-1)}(0)|_{L^\infty(\g_-)}\leq C|h^{i-2(k-1)}(0)|_{L^\infty(\g_+)}+|wr|_{L^\infty(\g_-)}.
\end{align}
Hence, for  $J_{16}$, it follows from \eqref{S3.27}, \eqref{S3.32-1} and  Lemma \ref{lem2.1} that
\begin{align}\label{S3.32}
|J_{16}| &\leq C e^{-\frac18|v|^2} \left(\frac12\right)^{C_2 T_0^{\frac43}}\cdot {|h^{i+1-2(k-1)}(0)|_{L^\infty(\g_-)}}\nonumber\\
&\leq C e^{-\frac18|v|^2} \left(\frac12\right)^{C_2 T_0^{\frac43}}\cdot \Big\{ |h^{i-2(k-1)}(0)|_{L^\infty(\g_+)}+|wr|_{L^\infty(\g_-)} \Big\}
\end{align}
where we have taken $k=C_1 T_0^{\frac43}$ and $T_0>1$ is a large constant to be  chosen later.

{For $J_7$}, it holds that %we have
\begin{align}\label{S3.34}
&|J_7| \leq  C e^{-\f18|v|^2} \sum_{l=1}^{k-1} \int_{\Pi _{j=1}^{l-1}\mathcal{V}_{2j-1}} \Pi_{j=1}^{l-1} \dd\sigma_{2j-1} \int_0^{t_{2l-1}} e^{-\nu_E^0(t-s)} \dd s \int_{\mathcal{V}_{2l-1}}\int_{\mathbb{R}^3}\nonumber\\
&\qquad\times  \Fi_{\{t_{2l}\leq0<t_{2l-1}\}} \tilde{w}(v_{2l-1}) |k_w(v_{2l-1},v') h^{i+1-2l}(x_{2l-1}-v_{2l-1,3}(t_{2l-1}-s),v')|\dd v' \dd\sigma_{2l-1} \nonumber\\
&=C e^{-\f18|v|^2} \sum_{l=1}^{k-1}\int_{\Pi _{j=1}^{l-1}\mathcal{V}_{2j-1}} \Pi_{j=1}^{l-1} \dd\sigma_{2j-1} \int_0^{t_{2l-1}} e^{-\nu_E^0(t-s)} \dd s \int_{\mathcal{V}_{2l-1}\cap \{|v_{2l-1}|\geq N\}}\int_{\mathbb{R}^3_{v'}} (\cdots)\dd v' \dd\sigma_{2l-1} \nonumber\\
&\,+C e^{-\f18|v|^2} \sum_{l=1}^{k-1}\int_{\Pi _{j=1}^{l-1}\mathcal{V}_{2j-1}} \Pi_{j=1}^{l-1} \dd\sigma_{2j-1} \int_0^{t_{2l-1}} e^{-\nu_E^0(t-s)} \dd s \int_{\mathcal{V}_{2l-1}\cap \{|v_{2l-1}|\leq N\}}\int_{\mathbb{R}^3_{v'}} (\cdots)\dd v' \dd\sigma_{2l-1}\nonumber\\
&:=\sum_{l=1}^{k-1} \left(J_{71l}+J_{72l}\right).
\end{align}
%\Red{Delete: \{For $J_{91l}$,\}}
We shall estimate the right-hand terms of \eqref{S3.34} as follows. A direct calculation shows
\begin{align}\label{S3.35}
\sum_{l=1}^{k-1} J_{71l}&\leq C e^{-\f18|v|^2} \sum_{l=1}^{k-1} \int_{\Pi _{j=1}^{l-1}\mathcal{V}_{2j-1}} \Pi_{j=1}^{l-1} \dd\sigma_{2j-1} \int_0^{t_{2l-1}} e^{-\nu_E^0(t-s)} \dd s\nonumber\\
&\qquad\qquad\qquad\times\int_{\mathcal{V}_{2l-1}\cap \{|v_{2l-1}|\geq N\}} e^{-\f18|v_{2l-1}|^2} dv_{2l-1} \cdot \sup_{1\leq l\leq k-1}\{\|h^{i+1-2l}\|_{L^\infty_{x,v}}\}\nonumber\\
&\leq Ck e^{-\frac1{16}N^2} \cdot e^{-\f18|v|^2}  \sup_{1\leq l\leq k-1}\{\|h^{i+1-2l}\|_{L^\infty_{x,v}}\}.
\end{align}
For each term $J_{72l}$, we have
\begin{align}\label{S3.36}
J_{72l}&\leq C e^{-\f18|v|^2} \int_{\Pi _{j=1}^{l-1}\mathcal{V}_{2j-1}} \Pi_{j=1}^{l-1} \dd\sigma_{2j-1} \int_{t_{2l-1}-\frac1N}^{t_{2l-1}} e^{-\nu_E^0(t-s)} \dd s\int_{\mathcal{V}_{2l-1}\cap \{|v_{2l-1}|\leq N\}}\int_{\mathbb{R}^3} (\cdots)\dd v' \dd\sigma_{2l-1}\nonumber\\
&\,+C e^{-\f18|v|^2} \int_{\Pi _{j=1}^{l-1}\mathcal{V}_{2j-1}} \Pi_{j=1}^{l-1} \dd\sigma_{2j-1} \int_0^{t_{2l-1}-\frac1N} e^{-\nu_E^0(t-s)} \dd s\int_{\mathcal{V}_{2l-1}\cap \{|v_{2l-1}|\leq N\}} e^{-\f18|v_{2l-1}|^2} \dd v_{2l-1}\nonumber\\
&\quad\times \int_{\{|v'|\geq 2N\}}  |k_w(v_{2l-1},v')| e^{\frac{|v_{2l-1}-v'|^2}{64}} \dd v' e^{-\frac{N^2}{64}}\cdot \sup_{1\leq l\leq k-1}\{\|h^{i+1-2l}\|_{L^\infty_{x,v}}\}\nonumber\\
&\,+C e^{-\f18|v|^2} \int_{\Pi _{j=1}^{l-1}\mathcal{V}_{2j-1}} \Pi_{j=1}^{l-1} \dd\sigma_{2j-1} \int_0^{t_{2l-1}-\frac1N} e^{-\nu_E^0(t-s)} \dd s \int_{\mathcal{V}_{2l-1}\cap \{|v_{2l-1}|\leq N\}} \int_{\{|v'|\leq 2N\}} \nonumber\\
&\quad\times  \Fi_{\{t_{2l}\leq0<t_{2l-1}\}} e^{-\frac18|v_{2l-1}|^2} |k_w(v_{2l-1},v') h^{i+1-2l}(x_{2l-1}-v_{2l-1,3}(t_{2l-1}-s),v')|\dd v' \dd v_{2l-1} \nonumber\\
&\leq C e^{-\f18|v|^2} \int_{\Pi _{j=1}^{l-1}\mathcal{V}_{2j-1}} \Pi_{j=1}^{l-1} \dd\sigma_{2j-1} \int_0^{t_{2l-1}-\frac1N} e^{-\nu_E^0(t-s)} \dd s \int_{\mathcal{V}_{2l-1}\cap \{|v_{2l-1}|\leq N\}} \int_{\{|v'|\leq 2N\}} \nonumber\\
&\quad\times  \Fi_{\{t_{2l}\leq0<t_{2l-1}\}} e^{-\frac18|v_{2l-1}|^2} |k_w(v_{2l-1},v') h^{i+1-2l}(x_{2l-1}-v_{2l-1,3}(t_{2l-1}-s),v')|\dd v' \dd v_{2l-1} \nonumber\\
&\quad+\frac{C}{N}e^{-\f18|v|^2}\cdot \|h^{i+1-2l}\|_{L^\infty_{x,v}}.
\end{align}
To estimate the first term on the right-hand side of  \eqref{S3.36},  a direct calculation shows
\begin{align}\label{S3.38}
&\int_{\mathcal{V}_{2l-1}\cap \{|v_{2l-1}|\leq N\}} \int_{\{|v'|\leq 2N\}}\Fi_{\{t_{2l}\leq0<t_{2l-1}\}} e^{-\frac18|v_{2l-1}|^2}\nonumber\\ &\qquad\qquad\times|k_w(v_{2l-1},v') h^{i+1-2l}(x_{2l-1}-v_{2l-1,3}(t_{2l-1}-s),v')|\dd v' \dd v_{2l-1}\nonumber\\
&\leq C_N \left\{\int_{\mathcal{V}_{2l-1}\cap \{|v_{2l-1}|\leq N\}} \int_{\{|v'|\leq 2N\}} e^{-\frac18|v_{2l-1}|^2} |k_w(v_{2l-1},v')|^2 \dd v' \dd v_{2l-1} \right\}^{\frac12}\nonumber\\
&\,\times \left\{\int_{\mathcal{V}_{2l-1}\cap \{|v_{2l-1}|\leq N\}} \int_{\{|v'|\leq 2N\}} \Fi_{\{t_{2l}\leq0<t_{2l-1}\}} \left|\frac{h^{i+1-2l}(y_{2l-1},v')}{w(v')}\right|^2 \dd v' \dd v_{2l-1} \right\}^{\frac12}\nonumber\\
&\leq C_N \left\{\int_{\mathcal{V}_{2l-1}\cap \{|v_{2l-1}|\leq N\}} \int_{\{|v'|\leq 2N\}} \Fi_{\{t_{2l}\leq0<t_{2l-1}\}} \left|\frac{h^{i+1-2l}(y_{2l-1},v')}{w(v')}\right|^2 \dd v' \dd v_{2l-1} \right\}^{\frac12}.
\end{align}
where we have used the notation $y_{2l-1}=x_{2l-1}-v_{2l-1,3}(t_{2l-1}-s)$. It is direct to know that $y_{2l-1}\in \Omega_d$ for $s\in [0,t_{2l-1}-\frac1N]$.  A direct computation
shows that
\begin{align}\label{S3.39}
\left|\f{\pa y_{2l-1}}{\pa v_{2l-1,3}}\right|=|t_{2l-1}-s|\geq \f{1}{N},\quad\mbox{for} \quad s\in [0,t_{2l-1}-\frac1N].
\end{align}
Thus, by making change of variable $v_{2l-1,3}\rightarrow y_{2l-1}$ and using \eqref{S3.39},  one obtains that
\begin{align}
&\left\{\int_{\mathcal{V}_{2l-1}\cap \{|v_{2l-1}|\leq N\}} \int_{\{|v'|\leq 2N\}} \Fi_{\{t_{2l}\leq0<t_{2l-1}\}} \left|\frac{h^{i+1-2l}(y_{2l-1},v')}{w(v')}\right|^2 \dd v' \dd v_{2l-1} \right\}^{\frac12}\nonumber\\
&\leq C_{N}\left\{\int_{\Omega_d}\int_{|v'|\leq 2N}\left|\frac{ h^{i+1-2l}(y_{2l-1},v')}{w(v')}\right|^2 \dd v' \dd y_{2l-1} \right\}^{\frac12}\leq C_N \left\|\frac{h^{i+1-2l}}{w}\right\|_{L^2_{x,v}},\nonumber
\end{align}
which together with \eqref{S3.38} and \eqref{S3.36} yield that
\begin{align}\label{S3.40}
J_{72l}\leq \frac{C}{N}e^{-\f18|v|^2}\cdot  \|h^{i+1-2l}\|_{L^\infty_{x,v}}+C_N  e^{-\f18|v|^2}\left\|\frac{h^{i+1-2l}}{w}\right\|_{L^2_{x,v}}.
\end{align}
Thus it follows from \eqref{S3.40}, \eqref{S3.35} and \eqref{S3.34} that
\begin{align}\label{S3.41}
|J_7| &\leq  \frac{Ck}{N}e^{-\f18|v|^2}\cdot \sup_{1\leq l\leq k-1}\{\|h^{i+1-2l}\|_{L^\infty_{x,v}}\}+C_N k  e^{-\f18|v|^2}\cdot \sup_{1\leq l\leq k-1}\left\{\left\|\frac{h^{i+1-2l}}{w}\right\|_{L^2_{x,v}}\right\}.
\end{align}
By similar arguments as in \eqref{S3.34}-\eqref{S3.41}, one can obtain
\begin{align}\label{S3.42}
|J_{9}|+|J_{12}|+|J_{14}|&\leq \frac{Ck}{N}e^{-\f18|v|^2}\cdot \sup_{1\leq l\leq k-1}\{\|h^{i+1-2l}\|_{L^\infty_{x,v}}+\|h^{i-2l}\|_{L^\infty_{x,v}}\}\nonumber\\
&\qquad+C_N k  e^{-\f18|v|^2}\cdot \sup_{1\leq l\leq k-1}\left\{\left\|\frac{h^{i+1-2l}}{w}\right\|_{L^2_{x,v}}+\left\|\frac{h^{i-2l}}{w}\right\|_{L^2_{x,v}}\right\}.
\end{align}
Now substituting \eqref{S3.41}-\eqref{S3.42}, \eqref{S3.32}, \eqref{S3.29}-\eqref{S3.31-1} and \eqref{S3.26} into \eqref{2.10}, we get, for $(t,x)\in[0,T_0]\times\Omega_d$ and $v_{0,3}>0$,  that
\begin{align}\label{S3.44}
|h^{i+1}(x,v)| & \leq \int_{\max\{t_1,0\}}^t e^{-\nu_E^0(t-s)}  \int_{\mathbb{R}^3}|k_w(v,v')h^{i}(x-v_{0,3}(t-s),v')|\dd v'\dd s
 + A_{1,i}(t,v),
\end{align}
where we have denoted
\begin{align*}%\label{S3.45}
A_{1,i}(t,v)&:=C e^{-\f1{32}|v|^2}\Big\{ k e^{-\nu_E^0 t}+ \left(\frac12\right)^{C_2T_0^{\frac43}}+\frac{k}N\Big\}\cdot \sup_{0\leq l\leq 2(k-1)}\{\|h^{i-l}\|_{L^\infty_{x,v}}+|h^{i-l}|_{L^\infty(\g_+)}\}\nonumber\\
&\quad+e^{-\nu_E^0 t} \|h^{i+1}\|_{L^\infty_{x,v}}+Ck \Big\{\|\nu^{-1}wg\|_{L^\infty_{x,v}}+|wr|_{{L^\infty(\g_-)}} + |w\tilde{r}|_{L^\infty(\gamma_-)}\Big\}\nonumber\\
&\quad+C_{N,k}\  e^{-\f18|v|^2}\cdot \sup_{1\leq l\leq 2(k-1)}\left\{\left\|\frac{ h^{i-l}}{w}\right\|_{L^2_{x,v}}\right\}.
\end{align*}

\medskip

\noindent{\it Step 2.} We consider the case $v_{0,3}<0$. In fact, it follows from \eqref{2.12} that
\begin{align}\label{S3.45}
|h^{i+1}(x,v)|&\leq  e^{-\nu_{E}^0 t} \|h^{i+1}\|_{L^\infty_{x,v}} + C\|\nu^{-1}wg\|_{L^\infty_{x,v}}\nonumber\\
&\quad+ \int_{\max\{t_1,0\}}^t e^{-\nu_E^0(t-s)}  \int_{\mathbb{R}^3}|k_w(v,v')h^{i}(x-v_{0,3}(t-s),v')|\dd v'\dd s \nonumber\\
&\quad + \Fi_{\{t_1> 0\}}\cdot e^{-\nu_{\v}(v) (t-t_1)}\cdot |h^{i}(d,v_1)|.
\end{align}
Noting $v_{1,3}=-v_{0,3}>0$ and $|v_1|=|\Mr v|=|v|$, then we apply \eqref{S3.44} to $h^{i}(d,v_1)$ to obtain
\begin{align}\label{S3.45-1}
&\Fi_{\{t_1> 0\}}\cdot e^{-\nu_{\v}(v) (t-t_1)}\cdot|h^{i}(d,v_1)| \nonumber\\
& \leq \Fi_{\{t_1> 0\}} \int_{\max\{t_2,0\}}^{t_1} e^{-\nu_E^0(t-s)}  \int_{\mathbb{R}^3}|k_w(v_1,v')h^{i}(d-v_{1,3}(t_1-s),v')|\dd v'\dd s
+ A_{2,i}(t,v)\nonumber\\
&\leq \Fi_{\{t_1> 0\}} \int_{\max\{t_2,0\}}^{t_1} e^{-\nu_E^0(t-s)}  \int_{\mathbb{R}^3}|k_w(\Mr v,v')h^{i}(2d-x+v_{0,3}(t-s),v')|\dd v'\dd s
+ A_{2,i}(t,v),
\end{align}
where we have used the fact $d-v_{1,3}(t_1-s)=2d-x+v_{0,3}(t-s)$ and denoted
\begin{align}
A_{2,i}(t,v)&:=C e^{-\f1{32}|v|^2}\Big\{ k e^{-\nu_E^0 t}+ \left(\frac12\right)^{C_2T_0^{\frac43}}+\frac{k}N\Big\}\cdot \sup_{1\leq l\leq 2(k-1)+1}\{\|h^{i-l}\|_{L^\infty_{x,v}}{+|h^{i-l}|_{L^\infty(\g_+)}}\}\nonumber\\
&\quad+e^{-\nu_E^0 t} \|h^{i}\|_{L^\infty_{x,v}}+Ck \Big\{\|wg\|_{L^\infty_{x,v}}+|wr|_{{L^\infty(\g_-)}} + |w\tilde{r}|_{L^\infty(\gamma_-)}\Big\}\nonumber\\
&\quad+C_{N,k}\  e^{-\f18|v|^2}\cdot \sup_{2\leq l\leq 2(k-1)+1}\left\{\left\|\frac{ h^{i-l}}{w}\right\|_{L^2_{x,v}}\right\}.\nonumber
\end{align}

Substituting \eqref{S3.45-1} into \eqref{S3.45}, one can get
\begin{align}\label{S3.45-2}
&|h^{i+1}(x,v)|\leq \int_{\max\{t_1,0\}}^t e^{-\nu_E^0(t-s)}  \int_{\mathbb{R}^3}|k_w(v,v')h^{i}(x-v_{0,3}(t-s),v')|\dd v'\dd s \nonumber\\
&\,\, + \Fi_{\{v_{0,3}<0\}}\Fi_{\{t_1> 0\}} \int_{\max\{t_2,0\}}^{t_1} e^{-\nu_E^0(t-s)}  \int_{\mathbb{R}^3}|k_w(\Mr v,v')h^{i}(2d-x+v_{0,3}(t-s),v')|\dd v'\dd s\nonumber\\
&\quad+e^{-\nu_{E}^0 t} \|h^{i+1}\|_{L^\infty_{x,v}} + A_{2,i}(t,v).
\end{align}

\medskip

\noindent{\it Step 3.} Combining \eqref{S3.44} and \eqref{S3.45-2}, for $t\in [0,T_0], (x,v)\in (\bar{\Omega}_d\times \R^3)\backslash (\g_0\cup \g_{-})$,  we obtain
\begin{align}\label{S3.45-3}
&|h^{i+1}(x,v)|\leq \int_{\max\{t_1,0\}}^t e^{-\nu_E^0(t-s)}  \int_{\mathbb{R}^3}|k_w(v,v')h^{i}(x-v_{0,3}(t-s),v')|\dd v'\dd s \nonumber\\
&\,+ \Fi_{\{v_{0,3}<0\}}\Fi_{\{t_1> 0\}} \int_{\max\{t_2,0\}}^{t_1} e^{-\nu_E^0(t-s)}  \int_{\mathbb{R}^3}|k_w(\Mr v,v')h^{i}(2d-x+v_{0,3}(t-s),v')|\dd v'\dd s\nonumber\\
&\quad+ A_{i}(t),
\end{align}
where we have denoted
\begin{align}
A_{i}:&=C \Big\{ k e^{-\nu_E^0 t}+ \left(\frac12\right)^{C_2T_0^{\frac43}}+\frac{k}N\Big\}\cdot \sup_{0\leq l\leq 2(k-1)+1}\{\|h^{i-l}\|_{L^\infty_{x,v}}{+|h^{i-l}|_{L^\infty(\g_+)}}\}\nonumber\\
&\,\,+e^{-\nu_{E}^0 t} \|h^{i+1}\|_{L^\infty_{x,v}} +Ck \Big\{\|\nu^{-1}wg\|_{L^\infty_{x,v}}+|wr|_{{L^\infty(\g_-)}} + |w\tilde{r}|_{L^\infty(\gamma_-)}\Big\}\nonumber\\
&\,\, +C_{N,k}\sup_{1\leq l\leq 2(k-1)+1}\left\{\left\|\frac{ h^{i-l}}{w}\right\|_{L^2_{x,v}}\right\}.
\end{align}

\noindent{\it Step 4.}
For later use, we denote $x'=x-v_{0,3}(t-s)\in \bar{\Omega}_d$ and $t_i'=t_i(s,x',v')$  for $s\in(\max\{t_1,0\},t)$. We also denote $\tilde{x}'=2d-x+v_{0,3}(t-s)\in\bar{\Omega}_d$ and $\tilde{t}_i'=t_i(s,\tilde{x},v')$. Noting \eqref{S3.45-3}, one has that
\begin{align}\label{S3.45-4}
&|h^{i}(x',v')|\leq \int_{\max\{t_1',0\}}^s e^{-\nu_E^0(s-\tau)}  \int_{\mathbb{R}^3}|k_w(v',v'')h^{i-1}(x'-v'_{0,3}(s-\tau),v'')|\dd v''\dd \tau \nonumber\\
&\,+ \Fi_{\{v'_{0,3}<0\}}\Fi_{\{t_1'> 0\}} \int_{\max\{t_2',0\}}^{t_1'} e^{-\nu_E^0(s-\tau)}  \int_{\mathbb{R}^3}|k_w(\Mr v',v'')h^{i-1}(2d-x'+v'_{0,3}(s-\tau),v'')|\dd v''\dd \tau\nonumber\\
&\,+ A_{i-1}(s),
\end{align}
and
\begin{align}\label{S3.45-5}
	&|h^{i}(\tilde{x}',v')|\leq \int_{\max\{\tilde{t}_1',0\}}^s e^{-\nu_E^0(s-\tau)}  \int_{\mathbb{R}^3}|k_w(v',v'')h^{i-1}(\tilde{x}'-v'_{0,3}(s-\tau),v'')|\dd v''\dd \tau \nonumber\\
	&\,+ \Fi_{\{v'_{0,3}<0\}}\Fi_{\{\tilde{t}_1'> 0\}} \int_{\max\{\tilde{t}_2',0\}}^{\tilde{t}_1'} e^{-\nu_E^0(s-\tau)}  \int_{\mathbb{R}^3}|k_w(\Mr v',v'')h^{i-1}(2d-\tilde{x}'+v'_{0,3}(s-\tau),v'')|\dd v''\dd \tau\nonumber\\
	&\,+ A_{i-1}(s).
\end{align}
Substituting \eqref{S3.45-4} and \eqref{S3.45-5} into \eqref{S3.45-3}, we can obtain
\begin{align}\label{S3.46}
|h^{i+1}(x,v)| &\leq A_i(t)+B_1+B_2+B_3+B_4+B_5,
\end{align}
where we have used the following notations
\begin{align}\label{S3.49-1}
B_1:=C\int_{0}^t e^{-\nu_E^0(t-s)}  \int_{\mathbb{R}^3_{v'}}\{|k_w(v,v')|+|k_w(\Mr v,v')|\} A_{i-1}(s)\dd v'\dd s,
\end{align}
\begin{align}
B_2:&=\int_{0}^tds \int_0^s e^{-\nu_E^0(t-\tau)} \dd\tau \int_{\mathbb{R}^3_{v'}}\int_{\mathbb{R}^3_{v''}}|k_w(v,v') k_w(v',v'')|\nonumber\\
&\qquad\times \Fi_{\{\max\{t_1,0\}<s<t\}}  \Fi_{\{\max\{t'_1,0\}<\tau<s\}}  |h^{i-1}(x'-v'_{0,3}(s-\tau),v'')|\dd v'' \dd v',\label{S3.49-2}\\
B_3:&=\int_{0}^tds \int_0^s e^{-\nu_E^0(t-\tau)} \dd\tau \int_{\mathbb{R}^3_{v'}}\int_{\mathbb{R}^3_{v''}}|k_w(v,v') k_w(\Mr  v',v'')|\Fi_{\{v_{0,3}'<0\}} \Fi_{\{t'_1>0\}} \nonumber\\
&\qquad\times \Fi_{\{\max\{t_1,0\}<s<t\}}  \Fi_{\{\max\{t'_2,0\}<\tau<t_1'\}}  |h^{i-1}(2d-x'+v'_{0,3}(s-\tau),v'')|\dd v'' \dd v',\label{S3.49-3}
\end{align}
\begin{align}
B_4:&=\Fi_{\{v_{0,3}<0\}} \Fi_{\{t_1>0\}}  \int_{0}^{t}ds\int_0^s \Fi_{\{\max\{t_2,0\}<s<t_1\}}  e^{-\nu_E^0(t-\tau)} d\tau\nonumber\\
& \times \int_{\R^3_{v'}} \int_{\R^3_{v''}}|k_w(\Mr v,v') k_w(v',v'')h^{i-1}(\tilde{x}'-v_{0,3}'(s-\tau),v'')|\cdot \Fi_{\{\max\{\tilde{t}'_1,0\}<\tau<s\}}\dd v''\dd v',\label{S3.49-4}
\end{align}
and
\begin{align}\label{S3.49-5}
B_5:&=\Fi_{\{v_{0,3}<0\}} \Fi_{\{t_1>0\}}  \int_{0}^{t}ds\int_0^s \Fi_{\{\max\{t_2,0\}<s<t_1\}}  e^{-\nu_E^0(t-\tau)} d\tau \nonumber\\
&\qquad \times \int_{\R^3_{v'}}\int_{\R^3_{v''}} |k_w(\Mr v,v') k_w(\Mr v',v'')h^{i-1}(2d-\tilde{x}'+v_{0,3}'(s-\tau),v'')|\nonumber\\
&\qquad\quad\times \Fi_{\{v'_{0,3}<0\}} \Fi_{\{\tilde{t}_1'>0\}} \Fi_{\{\max\{\tilde{t}'_2,0\}<\tau<\tilde{t}_1'\}}\dd v''\dd v'.
\end{align}

\medskip

For the term $B_1$, a direct calculation shows that
\begin{align}\label{S3.47}
B_1&\leq C \Big\{ k e^{-\f12 \nu_E^0 t}+ \left(\frac12\right)^{C_2T_0^{\frac43}}+\frac{k}N\Big\}\cdot \sup_{0\leq l\leq 2k}\{\|h^{i-l}\|_{L^\infty_{x,v}}+|h^{i-l}|_{L^\infty(\g_+)}  \}\nonumber\\
&\,\, +Ck \Big\{\|\nu^{-1}wg\|_{L^\infty_{x,v}}+|wr|_{{L^\infty(\g_-)}} + |w\tilde{r}|_{{L^\infty(\g_-)}}\Big\} +C_{N,k}\sup_{1\leq l\leq 2k}\left\{\left\|\frac{ h^{i-l}}{w}\right\|_{L^2_{x,v}}\right\}.
\end{align}

\smallskip

For the term $B_2$, we split the estimate into several cases.

\noindent{\it Case 1.} For $|v|\geq N$, we have
\begin{align}\label{S3.48}
B_2\leq  \frac{C}{N} \|h^{i-1}\|_{L^\infty_{x,v}}.
\end{align}

\noindent{\it Case 2.}  For $|v|\leq N, |v'|\geq 2N$ or $|v'|\leq 2N, |v''|\geq 3N$. In this case, we note from \eqref{2.15}  that
\begin{align*}%\label{S3.49}
\begin{cases}
\dis \int_{|v|\leq N, |v'|\geq 2N} \Big|k_w(v,v') e^{\frac{|v-v'|^2}{32}}\Big|\dd v'
\leq C (1+|v|)^{-1},\\[3mm]
\dis \int_{|v'|\leq 2N, |v''|\geq 3N} \Big|k_w(v',v'') e^{\frac{|v'-v''|^2}{32}}\Big|\dd v''
\leq C  (1+|v'|)^{-1}.
\end{cases}
\end{align*}
This yields  that
\begin{align}\label{S3.50}
&\int_{0}^t\dd s \int_0^s e^{-\nu_E^0(t-\tau)} \dd\tau \left\{\int_{|v|\leq N, |v'|\geq 2N} +\int_{|v'|\leq 2N, |v''|\geq 3N} \right\}(\cdots)\dd v'' \dd v'\nonumber\\
&\leq e^{-\frac{N^2}{32}}\|h^{i-1}\|_{L^\infty_{x,v}}\int_{|v|\leq N, |v'|\geq 2N} |k_w(v,v')e^{\frac{|v-v'|^2}{32}}| \cdot|k_w(v',v'')| \dd v'' \dd v'\nonumber\\
&\quad + e^{-\frac{N^2}{32}}\|h^{i-1}\|_{L^\infty_{x,v}} \int_{|v|\leq N, |v'|\geq 2N} |k_w(v,v')|\cdot |k_w(v',v'')e^{\frac{|v'-v''|^2}{32}}| \dd v'' \dd v'\nonumber\\
&\leq C e^{-\frac{N^2}{32}}\|h^{i-1}\|_{L^\infty_{x,v}}.
\end{align}

\noindent{\it Case 3.}  For $|v|\leq N$, $|v'|\leq 2N$, and $|v''|\leq 3N$, we first note that
\begin{align}\label{S3.51}
&\int_{0}^t\dd s \int_0^s e^{-\nu_E^0(t-\tau)} \dd\tau
\int_{|v'|\leq 2N, |v''|\leq 3N}(\cdots)\dd v'' \dd v'\nonumber\\
&\leq  \frac{C}{N} \|h^{i-1}\|_{L^\infty_{x,v}}+\int_{0}^t\dd s \int_0^{s-\frac1N} e^{-\nu_E^0(t-\tau)} \dd\tau  \int_{|v'|\leq 2N, |v''|\leq 3N} (\cdots)\dd v'' \dd v'\nonumber\\
&\leq \frac{C}{N}  \|h^{i-1}\|_{L^\infty_{x,v}} +C_{N,k}\int_{0}^t\dd s \int_0^{s-\frac1N} e^{-\nu_E^0(t-\tau)} \dd\tau
\left\{\iint |k_w(v,v') k_w(v',v'')|^2\dd v'' \dd v'\right\}^{\f12}\nonumber\\
&\quad \times\left\{\int_{|v'|\leq 2N, |v''|\leq 3N} \Fi_{\{\max\{t_1,0\}<s<t\}}  \Fi_{\{\max\{t'_1,0\}<\tau<s\}} \Big|\frac{h^{i-1}(y',v'')}{w(v'')}\Big|^2\dd v'' \dd v'\right\}^{\f12}\nonumber\\
&\leq \frac{C}N   \|h^{i-1}\|_{L^\infty_{x,v}}+C_{N,k}  \int_{0}^t\dd s \int_0^{s-\frac1N} e^{-\nu_E^0(t-\tau)} \dd\tau \nonumber\\
&\quad\times \left\{\int_{|v'|\leq 2N, |v''|\leq 3N} \Fi_{\{\max\{t_1,0\}<s<t\}}  \Fi_{\{\max\{t'_1,0\}<\tau<s\}} \Big|\frac{h^{i-1}(y',v'')}{w(v'')}\Big|^2\dd v'' \dd v'\right\}^{\f12},
\end{align}
where we have denoted $y'=x'-v'_{0,3}(s-\tau)\in \Omega_d$ for $s\in (\max\{t_1,0\}, t)$ and $\tau\in (\max\{t'_1,0\}, s)$. Similar to \eqref{S3.39}, we make  change of variable $v'_{0,3}\mapsto y'$, so that the second term on the right-hand side of \eqref{S3.51} is bounded as
\begin{align}\nonumber
&\int_{0}^t\dd s \int_0^{s-\frac1N} e^{-\nu_E^0(t-\tau)} \dd\tau  \left\{\int_{|v'|\leq 2N, |v''|\leq 3N} (\cdots)\dd v'' \dd v'\right\}^{\f12}\leq  CN^{\frac12}\left\|\frac{h^{i-1}}{w}\right\|_{L^2_{x,v}},
\end{align}
which together with \eqref{S3.51} yield that
\begin{equation}\label{S3.52}
\int_{0}^t\dd s \int_0^s e^{-\nu_E^0(t-\tau)} \dd\tau
\int_{|v'|\leq 2N, |v''|\leq 3N}(\cdots)\dd v'' \dd v'
\leq \frac{C}{N}  \|h^{i-1}\|_{L^\infty_{x,v}}+C_{N,k} \left\|\frac{\sqrt{\nu}h^{i-1}}{w}\right\|_{L^2_{x,v}}.
\end{equation}
Combining \eqref{S3.48}, \eqref{S3.50} and \eqref{S3.52}, we have
\begin{align}\label{S3.53}
B_2 \leq \frac{C}{N} \|h^{i-1}\|_{L^\infty_{x,v}}+C_{N,k} \left\|\frac{\sqrt{\nu}h^{i-1}}{w}\right\|_{L^2_{x,v}}.
\end{align}

By similar arguments as above, we can obtain
\begin{align}\label{S3.53-1}
B_3+B_4+B_5 \leq \frac{C}{N} \|h^{i-1}\|_{L^\infty_{x,v}}+C_{N,k} \left\|\frac{\sqrt{\nu}h^{i-1}}{w}\right\|_{L^2_{x,v}}.
\end{align}
Hence substituting \eqref{S3.53-1}, \eqref{S3.53}, \eqref{S3.47} into \eqref{S3.46}, we  have for  $t\in[0,T_0]$,
{\small\begin{align*}
|h^{i+1}(x,v)|&\leq C\left\{k e^{-\f12\nu_E^0 t} +\left(\frac12\right)^{{C_2T_0^{\frac43}}}
+\frac{k}N\right\} \cdot \sup_{0\leq l\leq 2k}\{\|h^{i-l}\|_{L^\infty_{x,v}}+|h^{i-l}|_{L^\infty(\g_+)}\}\nonumber\\
&\quad +e^{-\nu_E^0 t} \|h^{i+1}\|_{L^\infty_{x,v}}+C_{k}\Big\{\|\nu^{-1}wg\|_{L^\infty_{x,v}}+|wr|_{L^\infty(\g_-)}+ |w\tilde{r}|_{{L^\infty(\g_-)}}\Big\}\nonumber\\
&\quad+ C_{N,k}\  \sup_{0\leq l\leq 2k}\left\{\left\|\frac{h^{i-l}}{w}\right\|_{L^2_{x,v}}\right\}.
\end{align*}}
Now we take $k=C_1T_0^{\f43}$ and choose $t=T_0$. We first fix $T_0$ large enough, and then choose $N$ large enough, so that  one has $e^{-\nu_E^0t}\leq \frac12$ and
\begin{align}\nonumber
C\left\{k e^{-\f12\nu_E^0t}+\left(\frac12\right)^{{C_2T_0^{\frac43}}}
+\frac{k}N\right\}\leq \frac1{16}.
\end{align}
Therefore \eqref{3.13} follows. This completes  the proof of Lemma \ref{lemS3.3}. $\hfill\Box$

\subsection{Approximate solutions and uniform estimates}
Now we are in a position to  construct solutions of \eqref{S3.3} or equivalently \eqref{S3.8}. First of all, we consider the following approximate problem
\begin{align}\label{S3.54}
\begin{cases}
\dis \v f^n+ v_3 \partial_xf^n+p_E^0 \nu(v) f^n-p_E^0 Kf^n=g,\quad (x,v)\in \Omega_d\times \mathbb{R}^3,\\[1.5mm]
\dis f^n(0,v)|_{v_3>0}=(1-\frac{1}{n}) P_\g f^n(0,v)+r(v),\\[1.5mm]
\dis f^n(d,v)|_{v_3<0}=(1-\frac{1}{n}) f^n(d,\Mr v)
\end{cases}
\end{align}
where $\v\in(0,1]$ is arbitrary and $n>1$ is an integer. For later use, we choose $n_0>1$ large enough such that
$$
\frac18 (1-\frac2n+\f{3}{2n^2})^{-2k_0-1}\leq \f12, %\quad \mbox{and} \quad \f18(1-\frac1n)^{-\frac{k_0+1}{2}} \leq \frac12,
$$
for any $n\geq n_0$, where  $k_0>0$ is the one fixed  in Lemma \ref{lemS3.3}.

\begin{lemma}\label{lemS3.4}
Let $\v>0, d\geq 1$, $n\geq n_0$, and $\beta\geq 3$. Assume $\|\nu^{-1}wg\|_{L^\infty_{x,v}}+ |wr|_{L^\infty(\g_-)}<\infty$.
Then there exists a unique solution $f^n$ to \eqref{S3.54} satisfying
\begin{align*}%\label{S3.56}
\|wf^{n}\|_{L^\infty_{x,v}}+|wf^n|_{L^\infty(\gamma)}\leq C_{\v,n} \{\|\nu^{-1}wg\|_{L^\infty_{x,v}} + |wr|_{L^\infty(\g_-)} \},
\end{align*}
where the positive constant $C_{\v,n}>0$ depends only on $\v$ and $n$. Moreover, if
g is continuous in $\Omega_d\times\mathbb{R}^3$  and $r(v)$ is continuous in $\R^3_+$, then $f^n$ is continuous away from grazing set $\gamma_0$.
\end{lemma}

\noindent{\bf Proof.}  We consider the solvability of the  following boundary value problem
\begin{equation}\label{S3.56-1}
\begin{cases}
\dis \mathcal{L}_\lambda f:=\v f+v_3 \partial_x f+p_E^0 \nu(v) f - p_E^0 \lambda Kf=g,\\
\dis f(0,v)|_{v_3>0}=(1-\frac{1}{n}) P_\g f+r(v),\\
\dis f(d,v)|_{v_3<0}=(1-\frac{1}{n}) f(d,\Mr v)
\end{cases}
\end{equation}
for $\lambda\in[0,1]$. For brevity we denote $\mathcal{L}_\lambda^{-1}$ to be the solution operator associated to the problem, meaning that $f:=\mathcal{L}_\lambda^{-1} g$ is a solution to the BVP \eqref{S3.56-1}. Our idea is to prove the existence of $\mathcal{L}_0^{-1}$, and then extend to obtain the existence of  $\mathcal{L}_1^{-1}$ by a continuous argument on $\la$.  Since the proof is very long, we split it  into  four steps.   %\vspace{1mm}

\medskip

\noindent{\it Step 1.} In this step, we prove the existence of $\mathcal{L}_0^{-1}$. We consider the following approximate sequence
\begin{align}\label{S3.57}
\begin{cases}
\dis \mathcal{L}_0 f^{i+1}=\v f^{i+1}+v_3\partial_x f^{i+1} +p_E^0 \nu(v) f^{i+1}=g,\\
\dis f^{i+1}(0,v)|_{v_3>0}=(1-\frac{1}{n}) P_\g f^i+r(v),\\
\dis f^{i+1}(d,v)|_{v_3<0}=(1-\frac{1}{n}) f^i(d,\Mr v)
\end{cases}
\end{align}
for $i=0,1,2,\cdots$, where we have set $f^0\equiv0$. We will construct $L^\infty_{x,v}$ solutions to \eqref{S3.57} for $i=0,1,2,\cdots$, and
establish $L^\infty_{x,v}$-estimates.

\smallskip

Firstly, we will solve inductively the linear equation \eqref{S3.57} by %using the characteristic line defined \eqref{ode}.
the method of characteristics. Let $h^{i+1}(x,v)=w(v)f^{i+1}(x,v)$.  For %any  given
almost every $(x,v)\in\bar{\Omega}_d\times\mathbb{R}^3\backslash (\gamma_0\cup \gamma_-)$, %we have, for $t>t_{\mathbf{b}}(x,v)$, that
one can write
\begin{align}\label{S3.58}
h^{i+1}(x,v)
&=\Fi_{\{v_{0,3}>0\}} e^{-(\v+p_E^0\nu(v))t_{\mathbf{b}}}\cdot w(v)   \big[ (1-\frac{1}{n}) P_\g f^i(0,v)+r(v) \big]\nonumber\\
&\quad+\Fi_{\{v_{0,3}<0\}} e^{-(\v+p_E^0\nu(v))t_{\mathbf{b}}}\cdot w(v)   (1-\frac{1}{n}) f^{i}(d,\Mr v)   \nonumber\\
&\quad+\int_{t-t_{1}}^t  e^{-(\v+p_E^0\nu(v))(t-s)} (wg)(x-v_{0,3}(t-s),v) \dd s,
\end{align}
where $t_1=t-t_{\mathbf{b}}(x,v)$. %and for each $(x,v)\in\gamma_-$ we write
%\begin{align}
%h^{i+1}(x,v)=w(v)\Big[ (1-\frac{1}{n}) P_\g f^i+r(x,v)\Big].
%\end{align}
Noting the definition of $P_\g f$, we have
\begin{equation}\label{S3.59}
|w P_\g f(0)|_{L^\infty_v}\leq C|f(0)|_{L^\infty(\g_+)}.
\end{equation}

\smallskip

We consider \eqref{S3.58} with $i=0$. Since $f^0\equiv0$,
it is straightforward to see that
\begin{align*}%\label{S3.60}
\|h^{1}\|_{L^\infty_{x,v}}+|h^1|_{L^\infty(\gamma_+)}\leq   C \Big\{\|\nu^{-1}wg\|_{L^\infty_{x,v}}+|wr|_{L^\infty(\g_-)}\Big\}<\infty.
\end{align*}
Therefore we have obtained the solution to \eqref{S3.57} with $i=0$.  Assume that we have already solved \eqref{S3.57} for $i\leq l$ and obtained
\begin{equation}\label{S3.62}
\|h^{l+1}\|_{L^\infty_{x,v}}+|h^{l+1}|_{L^\infty(\gamma_+)}\leq C_{l+1}  \Big\{\|\nu^{-1}wg\|_{L^\infty_{x,v}}+|wr|_{L^\infty(\g_-)}\Big\}<\infty.
\end{equation}
We now consider \eqref{S3.57} for $i=l+1$. Noting \eqref{S3.62}, then we can solve \eqref{S3.57} by using \eqref{S3.58} with $i=l+1$. We still need to prove $h^{l+2}\in L^\infty_{x,v}$. Indeed, it follows from \eqref{S3.58} that
\begin{align*}%\label{S3.63}
\|h^{l+2}\|_{L^\infty_{x,v}}+{|h^{l+2}|_{L^\infty(\g_+)}}&\leq C |h^{l+1}|_{L^\infty{(\g_+)}} +C|wr|_{L^\infty(\g_-)} +C\|\nu^{-1}wg\|_{L^\infty_{x,v}}\nonumber\\
&\leq C_{l+2} \Big\{\|\nu^{-1}wg\|_{L^\infty_{x,v}}+|wr|_{L^\infty(\g_-)}\Big\}<\infty.
\end{align*}
Therefore, inductively, we have solved  \eqref{S3.57} for $i=0,1,2,\cdots$ and obtained
\begin{align}\label{S3.65}
\|h^{i}\|_{L^\infty_{x,v}}+|h^{i}|_{L^\infty(\gamma_+)}\leq C_{i} \Big\{\|\nu^{-1}wg\|_{L^\infty_{x,v}}+|wr|_{L^\infty(\g_-)}\Big\}<\infty.
\end{align}
The positive  constant $C_{i}$ may increase to infinity as $i\rightarrow \infty$.
Here, we emphasize that we first need to know the sequence $\{h^i\}_{i=0}^{\infty}$ is in $L^\infty_{x,v}$-space, otherwise one can not use Lemma \ref{lemS3.3} to get uniform $L^\infty_{x,v}$ estimates.

Let $(x,v)\in (\Omega_d\times\mathbb{R}^3)\backslash \gamma_0$, then it is easy to check that $t_{\mathbf{b}}(x,v)$ and $x_{\mathbf{b}}(x,v)$ are continuous for $v_{0,3}\neq0$. Therefore if $g$ and $r$ are continuous, by using \eqref{S3.58}, we can check that $f^i(x,v)$ is continuous away from grazing set $\g_0$.

\vspace{1mm}

Secondly, in order to take the limit $i\rightarrow\infty$, one has to get some uniform  estimates.  Multiplying \eqref{S3.57} by $f^{i+1}$ and integrating it over $\Omega_d\times\mathbb{R}^3$, one obtains that
\begin{align}\label{S3.66}
& \v\|f^{i+1}\|^2_{L^2_{x,v}}+\frac12\{|f^{i+1}(0)|^2_{L^2(\gamma_+)} +|f^{i+1}(d)|^2_{L^2(\gamma_+)}\}+\|f^{i+1}\|^2_{\nu}\nonumber\\
&\leq \frac12(1-\frac2n+\f{3}{2n^2})|f^{i}(0)|^2_{L^2(\gamma_+)} + \f12(1-\f1n)^2 |f^i(d)|^2_{L^2(\g_+)}\nonumber\\
&\quad+C_n |r|^2_{L^2(\g_-)} +\f{C}{\v}\|g\|^2_{L^2_{x,v}}+ \frac{\v}{4}\|f^{i+1}\|^2_{L^2_{x,v}}\nonumber\\
&\leq \frac12(1-\frac2n+\f{3}{2n^2})\{|f^{i}(0)|^2_{L^2(\gamma_+)} + |f^i(d)|^2_{L^2(\g_+)} \}\nonumber\\
&\quad+C_n |r|^2_{L^2(\g_-)} +\f{C}{\v}\|g\|^2_{L^2_{x,v}}+ \frac{\v}{4}\|f^{i+1}\|^2_{L^2_{x,v}}
\end{align}
where we have used the fact $|P_\g f(0)|_{L^2(\g_-)}=|P_\g f(0)|_{L^2(\g_+)}\leq |f(0)|^2_{L^2(\g_+)}$.
Then it follows from \eqref{S3.66} that
\begin{align}\label{S3.67}
&\f32\v\|f^{i+1}\|^2_{L^2_{x,v}}+\{|f^{i+1}(0)|^2_{L^2(\gamma_+)} +|f^{i+1}(d)|^2_{L^2(\gamma_+)}\}+2\|f^{i+1}\|^2_{\nu}\nonumber\\
&\leq (1-\frac2n+\f{3}{2n^2}) \{|f^{i}(0)|^2_{L^2(\gamma_+)} +|f^{i}(d)|^2_{L^2(\gamma_+)}\}+C_{\v,n} (\|g\|^2_{L^2_{x,v}}+|r|^2_{L^2(\g_-)}).
\end{align}
Now we take the difference $f^{i+1}-f^i$ in \eqref{S3.57}, then by similar energy estimate as above, we obtain
\begin{align}\label{S3.68}
&\f32\v\|f^{i+1}-f^i\|^2_{L^2_{x,v}}+\{|(f^{i+1}-f^i)(0)|^2_{L^2(\gamma_+)}+|(f^{i+1}-f^i)(d)|^2_{L^2(\gamma_+)}\}+2\|f^{i+1}-f^i\|^2_{\nu}\nonumber\\
&\leq (1-\frac2n+\f{3}{2n^2}) \{ |(f^{i}-f^{i-1})(0)|^2_{L^2(\gamma_+)} + |(f^{i}-f^{i-1})(d)|^2_{L^2(\gamma_+)} \}\nonumber\\
&\leq \cdots \leq (1-\frac2n+\f{3}{2n^2})^{i} \{|f^1(0)|^2_{L^2(\gamma_+)}+|f^1(d)|^2_{L^2(\gamma_+)}\}\nonumber\\
&\leq C_{\v,n} \cdot (1-\frac2n+\f{3}{2n^2})^{i} \, \Big( \|g\|^2_{L^2_{x,v}}+|r|^2_{L^2(\g_-)} \Big)\nonumber\\
&\leq  C_{\v,n,d} \cdot (1-\frac2n+\f{3}{2n^2})^{i}\, \Big\{\|\nu^{-1}wg\|_{L^\infty_{x,v}}+|wr|_{L^\infty(\g_-)}\Big\}<\infty.
\end{align}
Noting $1-\frac2n+\f{3}{2n^2}<1$, thus $\{f^{i}\}_{i=0}^\infty$ is a Cauchy sequence in $L^2_{x,v}$, i.e.,
\begin{align*}%\label{S3.69}
|f^{i}-f^j|^2_{L^2(\gamma_+)}+\|%f^{i+1}-f^i
{f^i-f^j}\|^2_{\nu}\rightarrow0,\quad\mbox{as} \  i,j\rightarrow\infty.
\end{align*}
And we also have, for $i=0,1,2,\cdots$, that
\begin{align}\label{S3.70}
|f^{i}|^2_{L^2(\gamma_+)}+\|f^{i}\|^2_{\nu}\leq C_{\v,n} \Big\{\|g\|^2_{L^2_{x,v}} + |r|^2_{L^2(\g_-)}\Big\}.
\end{align}

\vspace{1mm}

Next we consider the uniform  $L^\infty_{x,v}$ estimate. Here we point out that Lemma \ref{lemS3.3} still holds by replacing 1 by $1-\frac1n$ in the boundary condition of \eqref{2.9}, and the constants in Lemma \ref{lemS3.3} do not depend on $n\geq1$.  Thus we apply Lemma \ref{lemS3.3} to obtain that
\begin{align*}%\label{S3.71}
&\|h^{i+1}\|_{L^\infty_{x,v}}+|h^{i+1}|_{L^\infty{(\gamma_+)}}\nonumber\\
&\leq \frac18\sup_{0\leq l\leq 2 k_0} \{\|h^{i-l}\|_{L^\infty_{x,v}}+|h^{i-l}|_{L^\infty(\gamma_+)}\}
+C\Big\{\|\nu^{-1}wg\|_{L^\infty_{x,v}}+|wr|_{L^\infty(\g_-)}\Big\}\nonumber\\
&\quad+C\sup_{0\leq l\leq 2 k_0}  \|f^{i-l}\|_{L^2_{x,v}} \nonumber\\
&\leq\frac18\sup_{0\leq l\leq 2k_0} \{\|h^{i-l}\|_{L^\infty_{x,v}}+|h^{i-l}|_{L^\infty(\gamma_+)}\}
+C_{\v,n,d}\Big\{\|\nu^{-1}wg\|_{L^\infty_{x,v}}+|wr|_{L^\infty(\g_-)}\Big\},
\end{align*}
where we have used \eqref{S3.70} in the second inequality.  Now we apply Lemma \ref{lemA.1} to obtain that for $i\geq 2k_0+1$,
\begin{align}\label{S3.72}
\|h^{i}\|_{L^\infty_{x,v}}+|h^{i}|_{L^\infty{(\gamma_+)}}
&\leq\left(\frac18\right)^{[\frac{i}{2k_0+1}]}  \max_{1\leq l\leq 4k_0} \Big\{\|h^{l}\|_{L^\infty_{x,v}}+|h^{l}|_{L^\infty{(\gamma_+)}}\Big\}\nonumber\\
&\qquad+\frac{8+2k_0}{7} C_{\v,n,d} \Big\{\|\nu^{-1}wg\|_{L^\infty_{x,v}}+|wr|_{L^\infty(\g_-)}\Big\}\nonumber\\
&\leq C_{k_0,\v,n,d}  \Big\{\|\nu^{-1}wg\|_{L^\infty_{x,v}}+|wr|_{L^\infty(\g_-)}\Big\},
%\quad\mbox{for} \  i\geq k+1,
\end{align}
where we have used \eqref{S3.65} in the second inequality.
% of \eqref{S3.72}.
Hence it follows from \eqref{S3.72} and \eqref{S3.65} that
\begin{align}\label{S3.73}
\|h^{i}\|_{L^\infty_{x,v}}+|h^{i}|_{L^\infty{(\gamma_+)}}\leq C_{k_0,\v,n,d}  \Big\{\|\nu^{-1}wg\|_{L^\infty_{x,v}}+|wr|_{L^\infty(\g_-)}\Big\}, \quad\mbox{for} \  i\geq 1.
\end{align}
Taking the difference $h^{i+1}-h^i$ and then applying Lemma \ref{lemS3.3} to $h^{i+1}-h^i$, we have that for $i\geq 2k_0$,
\begin{align}\label{S3.74}
&\|h^{i+2}-h^{i+1}\|_{L^\infty_{x,v}}+|h^{i+2}-h^{i+1}|_{L^\infty(\g_+)}\nonumber\\
&\leq \frac18 \max_{0\leq l\leq 2k_0} \Big\{\|h^{i+1-l}-h^{i-l}\|_{L^{\infty}_{x,v}}+|h^{i+1-l}-h^{i-l}|_{L^\infty(\g_+)}\Big\}\nonumber\\
&\qquad+C \sup_{0\leq l\leq 2k_0} \Big\{\|f^{i+1-l}-f^{i-l}\|_{L^2_{x,v}}\Big\}\nonumber\\
&\leq\frac18 \max_{0\leq l\leq 2k_0} \Big\{\|h^{i+1-l}-h^{i-l}\|_{L^{\infty}_{x,v}}+|h^{i+1-l}-h^{i-l}|_{L^\infty(\g_+)}\Big\}\nonumber\\
&\quad+C_{\v,n,d} \cdot (1-\frac2n+\f{3}{2n^2})^{i-2k_0}\, \Big\{\|\nu^{-1}wg\|_{L^\infty_{x,v}}+|wr|_{L^\infty(\g_-)}\Big\}\nonumber\\
&\leq \frac18 \max_{0\leq l\leq 2k_0} \Big\{\|h^{i+1-l}-h^{i-l}\|_{L^{\infty}_{x,v}}+|h^{i+1-l}-h^{i-l}|_{L^\infty(\g_+)}\Big\}\nonumber\\
&\quad+C_{k_0,\v,n,d} \cdot  \Big\{\|\nu^{-1}wg\|_{L^\infty_{x,v}}+|wr|_{L^\infty(\g_-)}\Big\} \eta_n^{i+1},
\end{align}
where we have used \eqref{S3.68} and  denoted  $\eta_n:=1-\frac2n+\f{3}{2n^2}<1$. Here we choose $n$ large enough so that $\frac18 \eta_n^{-2k_0-1}\leq \frac12$, then it follows from \eqref{S3.74} and Lemma \ref{lemA.1}  that
\begin{align}\label{S3.75}
&\|h^{i+2}-h^{i+1}\|_{L^\infty_{x,v}}+|h^{i+2}-h^{i+1}|_{L^\infty(\gamma_+)}\nonumber\\
&\leq \left(\frac18\right)^{\left[\frac{i}{2k_0+1}\right]} \max_{1\leq l\leq 4k_0+1} \Big\{\|h^l\|_{L^\infty_{x,v}}+|h^{l}|_{L^\infty(\gamma_+)}\Big\}\nonumber\\
&\quad+C_{k_0,\v,n,d}  \Big\{\|\nu^{-1}wg\|_{L^\infty_{x,v}}+|wr|_{L^\infty(\g_-)}\Big\} \cdot  \eta_n^{i}\nonumber\\
&\leq C_{k_0,\v,n,d}  \Big\{\|\nu^{-1}wg\|_{L^\infty_{x,v}}+|wr|_{L^\infty(\g_-)}\Big\} \cdot\Bigg\{\left(\frac18\right)^{\left[\frac{i}{2k_0+1}\right]}  +\eta_n^{i}\Bigg\},
\end{align}
for $i\geq 2k_0+1$. Then \eqref{S3.75} implies immediately that $\{h^i\}_{i=0}^\infty$ is a Cauchy sequence in $L^\infty_{x,v}$, i.e., there exists a limit function $%h^n
{h}\in L^\infty_{x,v}$ so that $\|h^{i}-{h}\|_{L^\infty_{x,v}}+|h^{i}-{h}|_{L^\infty(\gamma_+)}\rightarrow 0$ as $i\rightarrow \infty$. Thus we obtained a function
{$f:=\f{h}{w}$} solves
\begin{equation*}%\label{S3.77}
\begin{cases}
\mathcal{L}_0 f=\v f+ v_3\partial_xf+\nu(v) f=g,\\[1.5mm]
\dis f(0,v)|_{v_3>0}=(1-\frac{1}{n}) P_\g f(0,v)+ r(v),\\[1.5mm]
\dis f(d,v)|_{v_3<0}=(1-\f1n) f(d,\Mr v),
\end{cases}
%\mbox{for} \  n\gg1.
\end{equation*}
with $n\geq n_0$ large enough.
Moreover, from \eqref{S3.73}, %and trace theorem,
there exists a constant $C_{\v,n,k_0,d}>0$ such that
\begin{align*}%\label{S3.78}
\|h\|_{L^\infty_{x,v}}+|h|_{L^\infty(\g)}\leq C_{\v,n,k_0,d} \Big\{\|\nu^{-1}wg\|_{L^\infty_{x,v}}+|wr|_{L^\infty(\g_-)}\Big\}.
\end{align*}

\medskip
\noindent{\it Step 2.} {\it A priori uniform estimates.} For any given $\lambda\in[0,1]$, let $f^n$ be the solution of \eqref{S3.56-1}, i.e.,
\begin{align}\label{S3.79}
\begin{cases}
\dis \mathcal{L}_\lambda f^n=\v f^{n}+v_3\partial_x f^{n}+p_E^0 \nu(v) f^{n}- p_E^0 \lambda  Kf^n=g,\\[1.5mm]
\dis f^{n}(0,v)|_{v_3>0}=(1-\frac{1}{n}) P_\g f^n(0,v)+ r(v),\\[1.5mm]
\dis f^n(d,v)|_{v_3<0}=(1-\f1n) f^n(d,\Mr v),
\end{cases}
\end{align}
Moreover we also assume that $\|wf^{n}\|_{L^\infty_{x,v}}+|wf^n|_{L^\infty(\gamma)}<\infty$. We note that $\langle \FL f^n, f^{n}\rangle\geq0 $, which implies that
\begin{equation}\label{S3.81}
	\langle Kf^n, f^{n}\rangle \leq \|f^n\|_{\nu}.
\end{equation}
On the other hand, a direct calculation shows that
\begin{align}\label{S3.82}
|f^{n}(0)|^2_{L^2(\gamma_-)}&=\left|(1-\frac{1}{n}) P_\g f^n (0,v)+ r(v) \right|_{L^2(\g_-)}^2\nonumber\\
&\leq (1-\f2n+\f{3}{2n^2}) |f^n|^2_{L^2(\g_{+})} +C_n |r|_{L^2(\g_-)}^2,
\end{align}
and
\begin{align}\label{S3.82-1}
|f^{n}(d)|^2_{L^2(\gamma_-)}=(1-\f1n)^2 |f^{n}(d)|^2_{L^2(\gamma_+)}.
\end{align}
Multiplying \eqref{S3.79} by $f^{n}$,   one has that
\begin{align}
&\v\|f^{n}\|^2_{L^2_{x,v}}+\frac12|f^{n}|^2_{L^2(\gamma_+)}-\frac12|f^{n}|^2_{L^2(\gamma_-)} + p_E^0 [\|f^n\|_{\nu}^2-\lambda \langle Kf^n, f^{n}\rangle]\nonumber\\
&\leq \frac{\v}{4} \|f^n\|_{L^2_{x,v}}+\frac{C}{\v}\|g\|^2_{L^2_{x,v}},\nonumber
\end{align}
which, together with \eqref{S3.81}, \eqref{S3.82} and \eqref{S3.82-1}, yields that
%\begin{align}\label{S3.80}
%\|f^{n}\|^2_{L^2_{x,v}}+|f^{n}|^2_{L^2(\gamma_+)} \leq C_{\v,n} \{\|g\|^2_{L^2_{x,v}} +   |r|^2_{L^2(\g_-)}\}.
%\end{align}
%Substituting \eqref{S3.81} and \eqref{S3.82} into \eqref{S3.80}, one has that
\begin{align}\label{S3.83}
\|\mathcal{L}^{-1}_\lambda g\|^2_{L^2_{x,v}}+|\mathcal{L}^{-1}_\lambda g|_{L^2(\g_{+})}^2&=\|f^{n}\|^2_{L^2_{x,v}} +|f^n|_{L^2(\g_{+})}^2
 \leq  C_{\v,n} \{\|g\|^2_{L^2_{x,v}} + |r|^2_{L^2(\g_{-})}\}.
\end{align}
Let $h^n:=w f^n$.  Then, by using \eqref{3.14} and \eqref{S3.83}, we  obtain
\begin{align}\label{S3.86}
\|w\mathcal{L}^{-1}_\lambda g\|_{L^\infty_{x,v}}+|w\mathcal{L}^{-1}_\lambda g|_{L^\infty(\gamma)}&=\|h^n\|_{L^\infty_{x,v}}+|h^n|_{L^\infty(\gamma)}\nonumber\\
&\leq C_{\v,n,k_0,d} \Big\{\|\nu^{-1}wg\|_{L^\infty_{x,v}}+|wr|_{L^\infty(\g_-)}\Big\}.
\end{align}

\vspace{1mm}

On the other hand, let $wg_1 \in L^\infty_{x,v}$ and $wg_2 \in L^\infty_{x,v}$.  Let  $f^n_1=\mathcal{L}^{-1}_\lambda g_1$ and $f^n_2=\mathcal{L}^{-1}_\lambda g_2$ be the solutions to \eqref{S3.79} with $g$ replaced by $g_1$ and $g_2$, respectively. Then we have that
\begin{align*}%\label{S3.87}
\begin{cases}
\dis \v (f^{n}_2-f^n_1)+v_3\partial_x (f^{n}_2-f^n_1)+p_E^0 \nu(v) (f^{n}_2-f^n_1)-p_E^0 \lambda K(f^{n}_2-f^n_1)=g_2-g_1,\\[2mm]
\dis (f^{n}_2-f^n_1)(0,v)|_{v_3>0}=(1-\frac1n) P_\g (f^{n}_2-f^n_1)(0,v),\\[2mm]
\dis (f^{n}_2-f^n_1)(d,v)|_{v_3<0}=(1-\f1n)(f^{n}_2-f^n_1)(d,\Mr v)
\end{cases}
\end{align*}
By similar arguments as in \eqref{S3.79}-\eqref{S3.86}, we can obtain
\begin{equation}\label{S3.88}
\|\mathcal{L}^{-1}_\lambda g_2-\mathcal{L}^{-1}_\lambda g_1\|^2_{L^2_{x,v}}
+ |\mathcal{L}^{-1}_\lambda g_2-\mathcal{L}^{-1}_\lambda g_1|^2_{L^2(\g_+)}
\leq  C_{\v,n,k_0} \|g_2-g_1\|^2_{L^2_{x,v}},
\end{equation}
and
\begin{equation}\label{S3.89}
\|w(\mathcal{L}^{-1}_\lambda g_2-\mathcal{L}^{-1}_\lambda g_1)\|_{L^\infty_{x,v}}
+|w(\mathcal{L}^{-1}_\lambda g_2-\mathcal{L}^{-1}_\lambda g_1)|_{L^\infty(\gamma)}\leq C_{\v,n,k_0,d}  \|\nu^{-1}w(g_2-g_1)\|_{L^\infty_{x,v}}.
\end{equation}
The uniqueness of solution to \eqref{S3.79} also follows from \eqref{S3.88}.
We point out that the constants $C_{\v,n,k_0,d}$ in \eqref{S3.83}-\eqref{S3.89} do not depend on $\lambda\in[0,1]$. This property is crucial for us to extend $\mathcal{L}_0^{-1}$ to $\mathcal{L}_1^{-1}$ by a bootstrap argument.

\medskip
\noindent{\it Step 3.}  In this step,  we shall prove the existence of solution $f^n$ to \eqref{S3.56-1} for sufficiently small $0<\lambda\ll1$, i.e., to prove the existence of operator $\mathcal{L}_{\lambda}^{-1}$. Firstly, we define the Banach space
\begin{align}\nonumber
\mathbf{X}:&=\Big\{f=f(x,v) :  \   wf\in L^\infty(\Omega_d\times\mathbb{R}^3), \ wf\in L^\infty(\gamma), \\  &\quad\quad \mbox{and} \  f(0,v)|_{v_3>0}=(1-\frac1n) P_\g f(0,v)+r(v),\,\,\, f(d,v)|_{v_3<0}=(1-\f1n) f(d,\Mr v) \Big\}.\nonumber
\end{align}
Now we define
\begin{align}\nonumber
T_\lambda f=\mathcal{L}_0^{-1} \Big(p_E^0\lambda K f+g\Big).
\end{align}
For any $f_1, f_2\in \mathbf{X}$, by using \eqref{S3.89}, we have  that
\begin{align}
&\|w(T_\lambda f_1-T_\lambda f_2)\|_{L^\infty_{x,v}}+|w(T_\lambda f_1-T_\lambda f_2)|_{L^\infty(\gamma)}\nonumber\\
&= \left\|w\{\mathcal{L}_0^{-1}(p_E^0\lambda Kf_1+g)-\mathcal{L}_0^{-1}(p_E^0\lambda Kf_2+g)\} \right\|_{L^\infty_{x,v}}\nonumber\\
&\qquad\quad+ \left|w\{\mathcal{L}_0^{-1}(p_E^0\lambda Kf_1+g)-\mathcal{L}_0^{-1}(p_E^0\lambda Kf_2+g)\} \right|_{L^\infty(\gamma)}\nonumber\\
&\leq C_{\v,n\,k_0,d}\|\nu^{-1}w\{(p_E^0\lambda Kf_1+g)-(p_E^0\lambda Kf_2+g)\}\|_{L^\infty_{x,v}}\nonumber\\
&\leq \lambda p_E^0 C_{\v,n,k_0,d} \|\nu^{-1} w (Kf_1-Kf_2)\|_{L^\infty_{x,v}} \nonumber\\
&\leq \lambda p_E^0 C_{\v,n,k_0,d} \|w(f_1-f_2)\|_{L^\infty_{x,v}}.\nonumber
\end{align}
We take $\lambda_\ast>0$ sufficiently small such that $\lambda_\ast p_E^0 C_{\v,k_0,d,n}\leq 1/2$, then $T_\lambda : \mathbf{X}\rightarrow \mathbf{X}$ is a contraction mapping for $\lambda\in[0,\lambda_\ast]$. Thus $T_\lambda$ has a fixed point, i.e.,  $\exists\, f^\lambda\in \mathbf{X}$ such that
\begin{equation}\nonumber
f^\lambda=T_\lambda f^\lambda=\mathcal{L}_0^{-1} \Big( p_E^0 \lambda K f^\lambda+g\Big),
\end{equation}
which yields immediately that
\begin{align}\nonumber
\mathcal{L}_\lambda f^\lambda=\v f^\lambda + v_3 \partial_x f^\lambda+p_E^0\nu(v) f^\lambda- p_E^0\lambda Kf^{\lambda}=g.
\end{align}
Hence, for any $\lambda\in[0,\lambda_\ast]$,  we have solved \eqref{S3.56-1} with $f^\lambda=\mathcal{L}_\lambda^{-1}g\in\mathbf{X}$.  Therefore we have obtained the existence of $\mathcal{L}_\lambda^{-1}$ for $\lambda\in[0,\lambda_\ast]$. Moreover the operator $\mathcal{L}_\lambda^{-1}$  has the properties \eqref{S3.83}-\eqref{S3.89}.

\medskip

\noindent{\it Step 4.} Finally we define
\begin{align}\nonumber
T_{\lambda_\ast+\lambda}f=\mathcal{L}_{\lambda_\ast}^{-1}\Big(\lambda K f+g\Big).
\end{align}
Noting the estimates for $\mathcal{L}_{\lambda_\ast}^{-1}$ are independent of $\lambda_\ast$. By similar arguments, we can prove $T_{\lambda_\ast+\lambda} : \mathbf{X}\rightarrow \mathbf{X}$ is a contraction mapping  for $\lambda\in[0,\lambda_\ast]$.  Then we obtain the existence of operator $\mathcal{L}_{\lambda_\ast+\lambda}^{-1}$, and  \eqref{S3.83}-\eqref{S3.89}.   Step by step, we can finally obtain the existence of operator $\mathcal{L}_1^{-1}$, and $\mathcal{L}_1^{-1}$ satisfies the estimates in \eqref{S3.83}-\eqref{S3.89}. The continuity is easy to obtain since the convergence of sequence under consideration is always in $L^\infty_{x,v}$.  Therefore we complete the proof of Lemma \ref{lemS3.4}.  $\hfill\Box$

\

Before taking the limit $n\to +\infty$, we first introduce a useful lemma which will be used later.
\begin{lemma}[\cite{EGKM2013,EGM2018}]\label{lem2.7}
Define the near grazing set $\g_{\pm}^{\epsilon}$ as
\begin{equation}\label{2.71}
\g_{\pm}^{\epsilon}:=\Big\{ (x,v)\in \g_{\pm}\,\, : \,\, |v|\leq \epsilon \,\, \mbox{or}\,\, |v|\geq \f{1}{\epsilon}  \Big\}.
\end{equation}
Let $\epsilon>0$ be a small  positive constant, then it holds that
\begin{align}\label{2.72}
|f(0) \Fi_{\g_{\pm}\backslash \g_{\pm}^{\epsilon}}|_{L^1(\g)} \leq \f{d}{\epsilon^2}\Big\{ \|f\|_{L^{1}_{x,v}} + \|v_3\pa_x f\|_{L^1_{x,v}}\Big\}.
\end{align}
More precisely, we can have
\begin{align}\label{2.73}
|f(0) \Fi_{\g_{\pm}\backslash \g_{\pm}^{\epsilon}}|_{L^1(\g)} \leq \f{4}{\epsilon^2}\Big\{ \|f\|_{L^{1}(\tilde{\Omega}\times \R^3)} + \|v_3\pa_x f\|_{L^1(\tilde{\Omega}\times \R^3)}\Big\},
\end{align}
where  $\tilde{\Omega}=(0,\f12)$.
\end{lemma}

\

\begin{lemma}\label{lemS3.6}
Let $\v>0, d\geq 1$ and $\beta\geq3$, and  assume $\|\nu^{-1}wg\|_{L^\infty_{x,v}} + |wr|_{L^\infty(\g_-)}<\infty$. Then there exists a unique solution $f^\v$ to solve the approximate linearized steady Boltzmann equation \eqref{S3.3}. Moreover, it satisfies
\begin{align}\label{S3.93}
\|wf^\v\|_{L^\infty_{x,v}} +|wf^\v|_{L^\infty(\gamma)} \leq C_{\v,d} \cdot  \Big\{\|\nu^{-1}wg\|_{L^\infty_{x,v}} + |wr|_{L^\infty(\g_-)} \Big\},
\end{align}
where the positive constant $C_{\v,d}>0$ depends only on $\v$ and $d$. Moreover, if
$g(x,v)$ is continuous in $\Omega_d\times \mathbb{R}^3$ and $r(v)$ is continuous in $\R^3_+$, then $f^\v$ is continuous away from the  grazing set $\gamma_0$.
\end{lemma}

\noindent{\bf Proof.} Let $f^n$ be the solution of \eqref{S3.54} constructed in Lemma \ref{lemS3.4} for $n\geq n_0$ with $n_0$ large enough.   Multiplying \eqref{S3.54} by $f^n$, one obtains that
\begin{align}\label{S3.94}
\v\|f^n\|^2_{L^2_{x,v}}+|f^n(0)|^2_{L^2(\g_+)}+2c_0 \|{(\mathbf{I-P})}f^n\|^2_{\nu}
\leq C_{\v} \|g\|^2_{L^2_{x,v}}+\left|(1-\f1n) P_\g f^n(0)+r\right|^2_{L^2(\g_-)}.
\end{align}
A direct calculation shows that
\begin{align}
\left|(1-\f1n) P_\g f^n(0)+r\right|^2_{L^2(\g_-)} &\leq |P_\g f^n(0)|^2_{L^2(\g_-)} +2 |P_\g f^n(0)|_{L^2(\g_-)} \cdot |r|_{L^2(\g_-)} + |r|^2_{L^2(\g_-)}\nonumber\\
&\leq  |P_\g f^n(0)|^2_{L^2(\g_+)} + \delta  |P_\g f^n(0)|^2_{L^2(\g_+)}  +C_{\delta} |r|^2_{L^2(\g_-)},\nonumber
\end{align}
which, together with \eqref{S3.94}, yields that
\begin{align}\label{S3.95}
&\v\|f^n\|^2_{L^2_{x,v}}+|(I-P_\g)f^n(0)|^2_{L^2(\g_+)}+2c_0 \|{(\mathbf{I-P})}f^n\|^2_{\nu}\nonumber\\
&\leq \delta  |P_\g f^n(0)|^2_{L^2(\g_+)} + C_{\v,\d} \big\{\|g\|^2_{L^2_{x,v}} + |r|^2_{L^2(\g_-)} \big\},
\end{align}
where $\delta>0$ is a small constant to be chosen later.

\smallskip

We still need to bound the first term on RHS of \eqref{S3.95}. Firstly, a direct calculation shows that
\begin{equation}\label{2.76}
\f12 |P_\g f^n(0)|^2_{L^2(\g_+)} \leq |P_\g f^n(0) I_{\g_+}|^2_{L^2(\g_+\backslash \g_+^{\epsilon})},
\end{equation}
provided that $\epsilon\ll 1$. We notice that
\begin{equation*}
f^n(0,v)=(I-P_\g) f^n(0,v) + P_\g f^n(0,v), \quad \forall \, (0,v)\in \g_{+},
\end{equation*}
which yields that
\begin{align}\label{2.77}
|P_\g f^n(0,v) I_{\g_+\backslash \g_+^{\epsilon}} |_{L^2(\g)}^2\leq 2 |f^n(0,v) I_{\g_+\backslash \g_+^{\epsilon}} |_{L^2(\g)}^2 + 2 |(I-P_\g) f^n(0,v) I_{\g_+\backslash \g_+^{\epsilon}} |_{L^2(\g)}^2.
\end{align}
On the other hand, it follows from \eqref{S3.54} that
\begin{equation*}
\f12 v_3 \pa_x(|f^n|^2)=-\v |f^n|^2- p_E^0 f^n \FL f^n +g f^n,
\end{equation*}
which immediately implies that
\begin{align}\label{2.78}
\|v_3 \pa_x(|f^n|^2) \|_{L^1_{x,v}}\leq C\Big\{ \|f^n \|_{L^2_{x,v}}^2 + \|{(\mathbf{I-P})} f^n \|_{\nu}^2 + \|g \|_{L^2_{x,v}}^2 \Big\}.
\end{align}
Thus, using Lemma \ref{lem2.7} and \eqref{S3.95} to obtain
\begin{align}\label{2.79-1}
& |f^n(0) I_{\g_+\backslash \g_+^{\epsilon}} |_{L^2(\g)}^2= |(f^n(0))^2 I_{\g_+\backslash \g_+^{\epsilon}} |_{L^1(\g)} \nonumber\\
& \leq   \f{4}{\epsilon^2}\Big\{ \|(f^n)^2\|_{L^{1}_{x,v}} + \|v_3\cdot\pa_x (f^n)^2\|_{L^1_{x,v}}\Big\} \nonumber\\
 &\leq \f{C}{\epsilon^2}\Big\{ \|f^n\|_{L^{2}_{x,v}}^2 + \|{(\mathbf{I-P})} f^n\|_{\nu}^2 +  \|g\|_{L^2_{x,v}}^2\Big\}\\
 &\leq C_{\v,\epsilon}\cdot \delta  |P_\g f^n(0)|^2_{L^2(\g_+)} + C_{\v,\epsilon,\d} \big\{\|g\|^2_{L^2_{x,v}} + |r|^2_{L^2(\g_-)} \big\},\nonumber
\end{align}
which, together with \eqref{S3.95}-\eqref{2.77}, yields that
\begin{align*}
 |P_\g f^n(0)|^2_{L^2(\g_+)}\leq C_{\v,\epsilon}\cdot \delta  |P_\g f^n(0)|^2_{L^2(\g_+)} + C_{\v,\epsilon,\d} \big\{\|g\|^2_{L^2_{x,v}} + |r|^2_{L^2(\g_-)} \big\}.
\end{align*}
Now taking $0<\delta\ll 1$ so that $C_{\v,\epsilon}\cdot \delta\leq \f12$, one obtains that
\begin{align}\label{2.79}
|P_\g f^n(0)|^2_{L^2(\g_+)}\leq  C_{\v,\epsilon} \big\{\|g\|^2_{L^2_{x,v}} + |r|^2_{L^2(\g_-)} \big\},
\end{align}
which, together with \eqref{S3.95}, shows that
\begin{equation}\label{2.80}
\|f^n\|^2_{L^2_{x,v}}+|f^n(0)|^2_{L^2(\g_+)}+2c_0 \|{(\mathbf{I-P})}f^n\|^2_{\nu} \leq C_{\v,\epsilon} \big\{\|g\|^2_{L^2_{x,v}} + |r|^2_{L^2(\g_-)} \big\}.
\end{equation}
We apply \eqref{3.14} and use \eqref{2.80} to obtain
\begin{align}\label{S3.101-1}
\|wf^n\|_{L^\infty_{x,v}}{+|wf^n|_{L^\infty(\g)}}
&\leq C\Big\{\|\nu^{-1}wg\|_{L^\infty_{x,v}}+ |wr|_{L^\infty(\g_-)} + \|f^n\|_{L^2_{x,v}}\Big\}\nonumber\\
&\leq C_{\v,\epsilon,d} \Big\{\|\nu^{-1}wg\|_{L^\infty_{x,v}} + |wr|_{L^\infty(\g_-)} \Big\}.
\end{align}

\smallskip

Taking the difference $f^{n_1}-f^{n_2}$ with $n_1,n_2\geq n_0$, we know that
\begin{align}\label{S3.102}
\begin{cases}
\dis\v (f^{n_1}-f^{n_2})+ v_3\partial_x (f^{n_1}-f^{n_2})+p_E^0 \FL (f^{n_1}-f^{n_2})=0,\\[2mm]
\dis (f^{n_1}-f^{n_2})(0,v)|_{v_{3}>0}=(1-\frac{1}{n_1}) P_\g  (f^{n_1}-f^{n_2})(0,v)+(\frac{1}{n_2}-\frac{1}{n_1}) P_\g f^{n_2}(0,v),\\[2mm]
\dis (f^{n_1}-f^{n_2})(d,v)|_{v_{3}<0}=(1-\f{1}{n_1}) (f^{n_1}-f^{n_2})(d,\Mr v) +(\f{1}{n_2}-\f{1}{n_1}) f^{n_2}(d,\Mr v)
\end{cases}
\end{align}
Multiplying \eqref{S3.102} by $ f^{n_1}-f^{n_2}$, and integrating it over $\Omega_d\times\mathbb{R}^3$,  by similar arguments as in \eqref{S3.94}-\eqref{2.80}, we can obtain
\begin{align}\label{S3.103}
&\|(f^{n_1}-f^{n_2})\|^2_{L^2_{x,v}} + | (f^{n_1}-f^{n_2})(0) |^2_{L^2(\g_+)}  + 2c_0 \|{(\mathbf{I-P})}(f^{n_1}-f^{n_2})\|^2_{\nu}\nonumber\\
&\leq C_{\v,\epsilon} \left| (\frac{1}{n_2}-\frac{1}{n_1}) P_\g f^{n_2}(0) \right|_{L^2(\g_-)}^2 +C_{\v,\epsilon} [\frac{1}{n_2}+\frac{1}{n_1}] \big\{|f^{n_1}|_{L^2(\g_+)}^2+|f^{n_2}|_{L^2(\g_+)}^2\big\}\nonumber\\
&\leq C_{\v,\epsilon} \cdot \Big\{ \f1{n_1}  + \f1{n_2} \Big\}\cdot \{|f^{n_1}|_{L^2(\g_+)}^2+|f^{n_2}|_{L^2(\g_+)}^2\}\nonumber\\
%&\leq C_{\v,\epsilon} \cdot \Big\{ \f1{n_1} + \f1{n_2}\Big\}\cdot  \Big\{\|g\|^2_{L^2} + |r|^2_{L^2(\g_-)} \Big\}\nonumber\\
&\leq  C_{d,\v,\epsilon} \cdot \Big\{\f1{n_1} + \f1{n_2}\Big\}\cdot  \Big\{\|\nu^{-1}wg\|_{L^\infty_{x,v}} + |wr|_{L^\infty(\g_-)} \Big\}^2\to 0,
\end{align}
as $n_1$, $n_2 \rightarrow \infty$,
where we have used the uniform estimate \eqref{S3.101-1} in the last inequality.  Applying \eqref{3.14} to $  f^{n_1}-f^{n_2}$  and using \eqref{S3.103}, then one has
\begin{align*}
&\|w(f^{n_1}-f^{n_2})\|_{L^\infty_{x,v}}+|w(f^{n_1}-f^{n_2})|_{L^\infty(\g)}\nonumber\\
&\leq C (\frac{1}{n_2}+\frac{1}{n_1}) | wf^{n_2} |_{L^\infty(\gamma_+)}+C\|f^{n_1}-f^{n_2}\|_{L^2_{x,v}}
\nonumber\\
&\leq  C_{\v,\epsilon,d} \cdot (\frac1{n_1}+\frac1{n_2}) \cdot  \Big\{\|\nu^{-1}wg\|_{L^\infty_{x,v}} + |wr|_{L^\infty(\g_-)} \Big\} \to 0,
\end{align*}
as $n_1,\ n_2 \rightarrow \infty$, which yields that $wf^n$ is a Cauchy sequence in $L^\infty_{x,v}$. We denote $f^\v=\lim_{n\rightarrow\infty} f^n$, then it is direct to check that $f^\v$ is a solution to \eqref{S3.3}, and \eqref{S3.93} holds. The continuity of $f^\v$ is easy to obtain since the convergence of  sequences is always in $L^\infty_{x,v}$ and $f^n$ is continuous away from the grazing set.  Therefore we have completed the proof of Lemma \ref{lemS3.6}. $\hfill\Box$

\

From now on, we assume
\begin{align}\label{3.57-1}
\int_0^d\int_{\R^3} \sqrt{\mu} g(x,v)dx dv=0\quad \mbox{and}\quad \int_{v_3>0}r(v) \sqrt{\mu} |v_3|dv=0.
\end{align}
Let $f^\v$ be the solution constructed in Lemma \ref{lemS3.6}, and we denote
\begin{equation*}
	\mathbf{P}f^\v(x,v)=\{a^\v(x)+b^\v\cdot v+\f12 c^\v(x) (|v|^2-3)\}\sqrt{\mu}.
\end{equation*}
Multiplying \eqref{S3.3} by $\sqrt{\mu}$ and integrating over $[0,d]\times \R^3$ to obtain
\begin{equation}\label{3.57}
\v \int_0^d\int_{\R^3} \sqrt{\mu} f^\v(x,v) dvdx\equiv \v\int_0^d a^\v(x)dx=0,
\end{equation}
where we have used \eqref{3.57-1}.

\smallskip

\begin{lemma}\label{lemS3.5}
Let $d\geq 1$. Assume \eqref{3.57-1} and  let $f^\v$ be the solution of \eqref{S3.3} constructed in Lemma \ref{lemS3.6},  then it holds that
\begin{align}\label{S3.92}
	\|\mathbf{P}f^\v\|^2_{L^2_{x,v}}\leq C d^6 \Big\{ \|{(\mathbf{I-P})}f^\v\|^2_{\nu}+\|g\|_{L^2_{x,v}}^2 + |(I-P_\g)f^\v|^2_{L^2(\g_+)} + |r|^2_{L^2(\g_-)}\Big\}.
\end{align}
\end{lemma}

\noindent{\bf Proof.} The weak formulation of \eqref{S3.3} is
\begin{align}\label{3.62}
&\v\int_0^d\int_{\R^3} f^\v(x,v) \psi(x,v) dvdx-\int_0^d\int_{\R^3}  v_3 f^\v(x,v)\partial_x \psi(x,v)  dvdx\nonumber\\
&=-\int_{\R^3} v_3f^\v(d,v) \psi(d,v) dv+\int_{\R^3} v_3f^\v(0,v) \psi(0,v) dv\nonumber\\
&\qquad-\int_0^d\int_{\R^3} \psi(x,v) \FL f^\v(x,v)  dvdx+\int_0^d\int_{\R^3} g(x,v) \psi(x,v) dvdx.
\end{align}
Motivated by  \cite{EGKM2013,EGM2018}, we choose some special test function $\psi$ to calculate the macroscopic part of $f^\v$.

\vspace{1mm}

\noindent{\it Step 1. Estimate on $c^\v$.} Define
\begin{equation*}
\zeta^\v_c(x)=\int_x^d c^\v(z)dz.
\end{equation*}
It is easy to check that
\begin{equation}\label{3.64}
\zeta^\v_c(d) =0,\quad |\zeta_c^\v(0)|\leq d^{\f12} \|c^\v\|_{L^2_x} \quad \mbox{and}\quad \|\zeta_c^\v\|_{L^2_x}\leq d \|c^\v\|_{L^2_x}.
\end{equation}
We define the test function $\psi$ in  \eqref{3.62} to be
\begin{equation*}%\label{3.65}
\psi=\psi_c^\v(x,v)=v_3 (|v|^2 -5) \sqrt{\mu} \, \zeta_c^\v(x).
\end{equation*}
Then the second term on LHS of \eqref{3.62} is estimated as
\begin{align}\label{3.66}
&-\int_0^d\int_{\R^3}  v_3 f^\v(x,v)\,\partial_x \psi_c^\v(x,v)  dvdx\nonumber\\
&=\int_0^d\int_{\R^3}   \big[a^\v+b^\v\cdot v+\f12 c^\v(x)(|v|^2-3) \big] v_3^2 (|v|^2 -5) \mu(v) c^\v(x) dvdx\nonumber\\
&\qquad+\int_0^d\int_{\R^3}    (\mathbf{I-P})f^\v(x,v) \, v_3^2 (|v|^2 -5) \sqrt{\mu} \, c^\v(x) dvdx\nonumber\\
&\geq 5 \|c^\v\|_{L^2_x}^2-C\|(\mathbf{I-P})f^\v\|_{\nu}\|c^\v\|_{L^2_x}
\geq 4 \|c^\v\|_{L^2_x}^2-C\|(\mathbf{I-P})f^\v\|_{\nu}^2,
\end{align}
where we have used
\begin{equation}\label{3.67}
\int_{\R^3} v_3^2 (|v|^2-3) (|v|^2 -5) \mu(v) dv=10, \quad
\int_{\R^3} v_3^2 (|v|^2 -5) \mu(v) dv=0.
\end{equation}

By using \eqref{3.67}, the first term on LHS of \eqref{3.62} is bounded as
\begin{align*}%\label{3.68}
\v\left|\int_0^d\int_{\R^3} f^\v(x,v) \, \psi_c^\v(x,v) dvdx\right|
&\leq C\v\|(\mathbf{I-P})f^\v\|_{\nu} \|\zeta_c^\v\|_{L^2_x} %\nonumber\\
%+ \v \left|\int_0^d\int_{\R^3}  b_3^\v \zeta_c(x) v_3^2  (|v-\fu|^2 -5) \mu(v)  dvdx\right|\nonumber\\
\leq C\v d \|(\mathbf{I-P})f^\v\|_{\nu} \|c^\v\|_{L^2_x}.
\end{align*}

For the boundary terms, we note
\begin{equation}\label{2.101}
f^\v(0,v)=P_\g f^\v(0,v) + I_{\g_+} \cdot (I-P_\g)f^\v(0,v) + I_{\g_-} \cdot r(v)\quad\mbox{for}\,\, (0,v)\in \gamma,
\end{equation}
where
\begin{align}
P_\g f^\v(0,v)=\sqrt{2\pi \mu(v)} \, z_{\g}(f^\v)\quad\mbox{with}\quad
 z_{\g}(f^\v):=\int_{u_3<0} f^\v(0,u) \sqrt{\mu(u)} |u_3| du.
\end{align}
Noting \eqref{3.64}, one has that
\begin{align}\label{3.69}
&\left|-\int_{\R^3} v_3f^\v(d,v) \psi_c^\v(d,v) dv+\int_{\R^3} v_3f^\v(0,v) \psi_c^\v(0,v) dv\right|\nonumber\\
&=\left|\int_{\R^3} v_3^2 (|v|^2-5) \sqrt{\mu(v)} \zeta^\v_c(0)\, \Big[\sqrt{2\pi \mu(v)} \, z_\g(f^\v) + I_{\g_+} \cdot (I-P_\g)f^\v(0,v) + I_{\g_-} \cdot r(v)\Big]dv\right|\nonumber\\
&=| \zeta^\v_c(0)|\cdot\left|\int_{\R^3} v_3^2 (|v|^2-5) \sqrt{\mu(v)}\, \Big[ I_{\g_+} \cdot (I-P_\g)f^\v(0,v) + I_{\g_-} \cdot r(v)\Big]dv\right|\nonumber\\
&\leq Cd^{\f12} \|c^\v\|_{L^2_x} \, \big[ |(I-P_\g)f^\v(0)|_{L^2(\g_+)} + |r|_{L^2(\g_-)}\big].
\end{align}
Hence the RHS of \eqref{3.62} is bounded as
\begin{align}\label{3.70}
\mbox{RHS of } \eqref{3.62} \leq C d \Big( \|(\mathbf{I-P})f^\v\|_{\nu} +\|g\|_{L^2_{x,v}} + |(I-P_\g)f^\v(0)|_{L^2(\g_+)} + |r|_{L^2(\g_-)}\Big)\|c^\v\|_{L^2_x}.
\end{align}
Combining \eqref{3.66}, \eqref{2.101} and \eqref{3.70}, one obtains
\begin{align}\label{3.71}
\|c^\v\|_{L^2_x}^2\leq C d^2 \Big( \|(\mathbf{I-P})f^\v\|_{\nu}^2 +\|g\|_{L^2_{x,v}}^2 + |(I-P_\g)f^\v(0)|_{L^2(\g_+)}^2 + |r|_{L^2(\g_-)}^2 \Big).
\end{align}

\vspace{1mm}

\noindent{\it Step 2. Estimate on $b^\v$.} We define
\begin{align*}
\zeta_{b,i}^\v(x)=\int_x^d b^\v_{i}(z)dz,\quad i=1,2,3.
\end{align*}
It holds that
\begin{equation}\label{3.72}
\begin{split}
&\zeta^\v_{b,i}(d)=0, \quad |\zeta^\v_{b,i}(0)|\leq d^{\f12} \|b_i^\v\|_{L^2_x} \quad \mbox{and}\quad \|\zeta^\v_{b,i}\|_{L^2_x}\leq d \|b^\v_i\|_{L^2_x}, \quad i=1,2,3.
\end{split}
\end{equation}
%For $b^\v_3$, we define
%\begin{align}
%\zeta_{b,3}(x)=-\int_0^x b_3^\v(z)dz+\frac1d \int_0^d(d-z) b_3^\v(z) dz.
%\end{align}
Now we take the test function $\psi$ in \eqref{3.62} to be
\begin{align*}
\psi=\psi_{b,3}^\v(x,v)=(v_3^2-1) \sqrt{\mu(v)} \, \zeta_{b,3}^\v(x).
\end{align*}
Then the second term on LHS of \eqref{3.62} is controlled as
\begin{align}\label{3.75}
&-\int_0^d\int_{\R^3}  v_3 f^\v(x,v) \, \partial_x \psi_{b,3}^\v(x,v)  dvdx\nonumber\\
&=\int_0^d\int_{\R^3}   [a^\v+b^\v\cdot v+\f12 c^\v(x)(|v|^2-3)] \, v_3 (v_3^2-1)  \mu(v) \, b_3^\v(x)  dvdx\nonumber\\
&\qquad+\int_0^d\int_{\R^3} (\mathbf{I-P})f^\v(x,v)\cdot v_3(v_3^2-1) \sqrt{\mu(v)} \, b_3^\v(x)  dvdx\nonumber\\
&=2\|b_3^\v\|_{L^2_x}^2 + \int_0^d\int_{\R^3} (\mathbf{I-P})f^\v(x,v)\cdot v_3(v_3^2-1) \sqrt{\mu(v)} \, b_3^\v(x)  dvdx\nonumber\\
&\geq \frac74 \|b_3^\v\|_{L^2_x}^2-C\|(\mathbf{I-P})f^\v\|_{\nu}^2,
\end{align}
where we have used
\begin{equation*}
\begin{split}
\int_{\R^3} v_3^2 (v_3^2-1) \mu(v)dv=2.
\end{split}
\end{equation*}
The first term on LHS of \eqref{3.62} is bounded as
\begin{align*}%\label{3.77}
 \v\left|\int_0^d\int_{\R^3} f^\v(x,v) \psi_{b,3}^\v(x,v) dvdx\right|
&\leq C\v d \|c^\v\|_{L^2_x}\|b_3^\v\|_{L^2_x}+Cd\v\|(\mathbf{I-P})f^\v\|_{\nu} \|b_3^\v\|_{L^2_x} \nonumber\\
&\leq \frac14 \|b_3^\v\|_{L^2_x}^2+C\v^2 d^2 \Big(\|(\mathbf{I-P})f^\v\|_{\nu}^2+ \|c^\v\|_{L^2_x}^2\Big).
\end{align*}
For the boundary terms on RHS of \eqref{3.62}, it follows from \eqref{3.72} that
\begin{align}\label{3.78}
&\left|-\int_{\R^3} v_3f^\v(d,v) \, \psi_{b,3}^\v(d,v) dv+\int_{\R^3} v_3f^\v(0,v) \, \psi_{b,3}^\v(0,v) dv\right|\nonumber\\
&=\left|\int_{\R^3} v_3 (v_3^2-1) \sqrt{\mu(v)} \, \zeta^\v_{b,3}(0)\, \Big[\sqrt{2\pi \mu(v)} \, z^\v_\g(0) + I_{\g_+} \cdot (I-P_\g)f^\v(0,v) + I_{\g_-} \cdot r(v)\Big]dv\right|\nonumber\\
&=\left|\int_{\R^3} v_3 (v_3^2-1) \sqrt{\mu(v)} \, \zeta^\v_{b,3}(0)\, \Big[I_{\g_+} \cdot (I-P_\g)f^\v(0,v) + I_{\g_-} \cdot r(v)\Big]dv\right|\nonumber\\
&\leq Cd^{\f12} \|b_3^\v\|_{L^2_x} \, \big[ |(I-P_\g)f^\v(0)|_{L^2(\g_+)} + |r|_{L^2(\g_-)}\big].
\end{align}
Thus   we have
\begin{equation}\label{3.79}
\mbox{RHS of } \eqref{3.62} \leq Cd \|b_3^\v\|_{L^2_x}\cdot  \Big\{\|(\mathbf{I-P})f^\v\|_{\nu}+ \|g\|_{L^2_{x,v}} +  |(I-P_\g)f^\v(0)|_{L^2(\g_+)} + |r|_{L^2(\g_-)}\Big\}.
\end{equation}
Then combining \eqref{3.75}-\eqref{3.79} and using \eqref{3.71}, one obtains that
\begin{align}\label{3.80}
\|b_3^\v\|_{L^2_x}^2\leq Cd^4\Big\{\|(\mathbf{I-P})f^\v\|_{\nu}^2+ \|g\|_{L^2_{x,v}}^2 + |(I-P_\g)f^\v(0)|^2_{L^2(\g_+)} + |r|^2_{L^2(\g_-)}\Big\}.
\end{align}

\medskip

For the estimate of  $b_i^\v(x), i=1,2$, we define
\begin{align}
\psi=\psi_{b,i}(x,v)=|v|^2 v_i v_3 \sqrt{\mu(v)} \, \zeta_{b,i}^\v(x), \quad i=1,2.
\end{align}
By a tedious calculation, one has
\begin{align}\label{3.102-1}
&-\int_0^d\int_{\R^3}  v_3 f^\v(x,v) \, \partial_x \psi_{b,i}^\v(x,v)  dvdx\nonumber\\
&=\int_0^d\int_{\R^3}   [a^\v+b^\v\cdot v+\f12 c^\v(x)(|v|^2-3)] \, |v|^2 v_i v_3^2 \, \mu(v) \, b_i^\v(x)  dvdx\nonumber\\
&\qquad+\int_0^d\int_{\R^3} (\mathbf{I-P})f^\v(x,v)\cdot \, |v|^2 v_i v_3^2 \, \mu(v) \, b_i^\v(x)  dvdx\nonumber\\
&=7\|b^\v_i\|^2_{L^2_x} \int_{\R^3} |v|^2 v_i^2 v_3^2 \, \mu(v) dv +\int_0^d\int_{\R^3} (\mathbf{I-P})f^\v(x,v)\cdot \, |v|^2 v_i v_3^2 \, \mu(v) \, b_i^\v(x)  dvdx  \nonumber\\
&\geq 6\|b_i^\v\|_{L^2_x}^2-C\|(\mathbf{I-P})f^\v\|_{\nu}^2,
\end{align}
where we have used
\begin{equation*}
\int_{\R^3} |v|^2 v_i^2 v_3^2 \, \mu(v) dv=7 ,\quad i=1,2.
\end{equation*}
The first term on LHS of \eqref{3.62} is bounded as
\begin{align}\label{3.77}
	\v\left|\int_0^d\int_{\R^3} f^\v(x,v) \, \psi_{b,i}^\v(x,v) dvdx\right|
	&\leq Cd\v\|(\mathbf{I-P})f^\v\|_{\nu} \|b_i^\v\|_{L^2_x} \nonumber\\
	&\leq \frac14 \|b_i^\v\|_{L^2_x}^2 + C\v^2 d^2 \|(\mathbf{I-P})f^\v\|_{\nu}^2.
\end{align}
For the boundary term, we have
\begin{align}
&\left|-\int_{\R^3} v_3f^\v(d,v) \, \psi_{b,i}^\v(d,v) dv+\int_{\R^3} v_3f^\v(0,v) \, \psi_{b,i}^\v(0,v) dv\right|\nonumber\\
&=\left|\int_{\R^3} |v|^2 v_iv_3 \sqrt{\mu(v)} \sqrt{\mu(v)} \, \zeta^\v_{b,i}(0)\, \Big[\sqrt{2\pi \mu(v)} \, z^\v_\g(0) + I_{\g_+} \cdot (I-P_\g)f^\v(0,v) + I_{\g_-} \cdot r(v)\Big]dv\right|\nonumber\\
&=\left|\int_{\R^3} |v|^2 v_iv_3 \sqrt{\mu(v)} \sqrt{\mu(v)} \, \zeta^\v_{b,i}(0)\, \Big[I_{\g_+} \cdot (I-P_\g)f^\v(0,v) + I_{\g_-} \cdot r(v)\Big]dv\right|\nonumber\\
&\leq Cd^{\f12} \|b_i^\v\|_{L^2_x} \, \big[ |(I-P_\g)f^\v(0)|_{L^2(\g_+)} + |r|_{L^2(\g_-)}\big].\nonumber
\end{align}
Thus the terms on RHS  of \eqref{3.62} is bounded as
\begin{align}\label{3.102-2}
\mbox{RHS of } \eqref{3.62} \leq Cd \|b_i^\v\|_{L^2_x}\cdot  \Big\{\|(\mathbf{I-P})f^\v\|_{\nu}+ \|g\|_{L^2_{x,v}} +  |(I-P_\g)f^\v(0)|_{L^2(\g_+)} + |r|_{L^2(\g_-)}\Big\}.
\end{align}
Substituting \eqref{3.102-1}, \eqref{3.77} and \eqref{3.102-2} into \eqref{3.62}, it holds that
\begin{equation}\label{3.102-3}
\|b_i^\v\|_{L^2_x}^2\leq Cd^2\Big\{\|(\mathbf{I-P})f^\v\|_{\nu}^2+ \|g\|_{L^2_{x,v}}^2 + |(I-P_\g)f^\v(0)|^2_{L^2(\g_+)} + |r|^2_{L^2(\g_-)}\Big\}, \,\,\, i=1,2.
\end{equation}

\

\noindent{\it Step 3. Estimate on $a^\v$.} Define
\begin{equation*}
\zeta_a^\v(x)=-\int_0^x a^\v(z)dz.
\end{equation*}
From \eqref{3.57}, it is easy to check that
\begin{equation}\label{3.81}
\zeta_a^\v(x)|_{x=0,d}=0,\quad \mbox{and}\quad \|\zeta^\v_a\|_{L^2_x}\leq d \|a^\v\|_{L^2_x}.
\end{equation}
We define the test function $\psi$ in  \eqref{3.62} to be
\begin{equation*}%\label{3.82}
\psi=\psi^\v_a(x,v)=(|v|^2 -10) \sqrt{\mu}\, v_3\, \zeta^\v_a(x).
\end{equation*}
Then the second term on LHS of \eqref{3.62} is estimated as
\begin{align}\label{3.84}
&-\int_0^d\int_{\R^3}  v_3 f^\v(x,v)\partial_x \psi_a(x,v)  dvdx\nonumber\\
&=\int_0^d\int_{\R^3}   [a^\v+b^\v\cdot v+\f12 c^\v(x)(|v|^2-3)]\cdot (|v|^2 -10) v_3^2\mu(v) a^\v(x) dxdv\nonumber\\
&\qquad+\int_0^d\int_{\R^3}  v_3 (\mathbf{I-P})f^\v(x,v)(|v|^2 -10) v_3^2\mu(v) a^\v(x) dvdx\nonumber\\
&\geq 5\|a^\v\|_{L^2_x}^2-C\|(\mathbf{I-P})f^\v\|_{\nu}\|a^\v\|_{L^2_x} \geq 4\|a^\v\|_{L^2_x}^2-C\|(\mathbf{I-P})f^\v\|_{\nu}^2,
\end{align}
where we have used
\begin{align}\nonumber
\int_{\R^3} v_3^2 \cdot (|v|^2 -10) \mu(v) dv=5\quad\mbox{and}\quad\int_{\R^3} (|v|^2-3) v_3^2 \cdot (|v|^2 -10) \mu(v) dv=0.
\end{align}
A direct calculation shows that
\begin{align*}%\label{3.85}
&\v\left|\int_0^d\int_{\R^3} f^\v(x,v) \psi_a(x,v) dvdx\right| \nonumber\\
&\leq \left|\v\int_0^d\int_{\R^3} b_3^\v(x) \zeta_{a}^\v(x) v_3^2(|v|^2 -10) \mu(v)dvdx \right|+C\v d\|(\mathbf{I-P})f^\v\|_{\nu} \|a^\v\|_{L^2_x} \nonumber\\
&\leq \|a^\v\|_{L^2_x}^2+Cd^2\v^2 (\|(\mathbf{I-P})f^\v\|_{\nu}^2+\|b^\v_3\|_{L^2_x}^2).
\end{align*}
It follows from  \eqref{3.81} that
\begin{equation*}%\label{3.86}
-\int_{\R^3} v_3f^\v(d,v) \psi_a^\v(d,v) dv+\int_{\R^3} v_3f^\v(0,v) \psi_a^\v(0,v) dv=0.
\end{equation*}
Hence, for the RHS of \eqref{3.62}, it holds that
\begin{align}\label{3.87}
\mbox{RHS of } \eqref{3.62} \leq Cd \|a^\v\|_{L^2_x} \Big\{\|(\mathbf{I-P})f^\v\|_{\nu}+ \|g\|_{L^2_{x,v}}\Big\}.
\end{align}
Combining \eqref{3.84}-\eqref{3.87} and using \eqref{3.80}, we obtain
\begin{align}\label{3.88}
\|a^\v\|_{L^2_x}^2&\leq Cd^2 \{\|(\mathbf{I-P})f^\v\|_{\nu}^2+ \|g\|_{L^2_{x,v}}^2
+ \|b^\v_3\|^2_{L^2_x}\}\nonumber\\
&\leq Cd^6\Big\{\|(\mathbf{I-P})f^\v\|_{\nu}^2+ \|g\|_{L^2_{x,v}}^2 + |(I-P_\g)f^\v(0)|^2_{L^2(\g_+)} + |r|^2_{L^2(\g_-)}\Big\}.
\end{align}
Therefore, \eqref{S3.92} follows directly from \eqref{3.88}, \eqref{3.102-3}, \eqref{3.80} and \eqref{3.71}. The proof of Lemma \ref{lemS3.5} is complete. $\hfill\Box$

\vspace{2mm}

\begin{lemma}\label{lemS3.7}
Let $d\geq 1$, $\beta\geq 3$.  Assume \eqref{3.57-1} and  $\|\nu^{-1}wg\|_{L^\infty_{x,v}}+|wr|_{L^\infty(\g_-)}<\infty$.  Under the condition $\dis\int_0^d\int_{\R^3} f(x,v) \sqrt{\mu} dvdx=0$, there exists a  unique solution $f=f(x,v)$  to the linearized steady Boltzmann equation
\begin{align}\label{S3.106}
\begin{cases}
v_3 \partial_xf+p_E^0 \FL f=g,\quad (x,v)\in \Omega_d\times \R^3, \\
f(0,v)|_{v_3>0}=P_\g f(0,v) + r(v),\\
f(d,v)|_{v_3<0}=f(x,\Mr v),
\end{cases}
\end{align}
with
\begin{equation}\label{S3.107}
\|wf\|_{L^\infty_{x,v}} +|wf|_{L^\infty(\gamma)} \leq C_d \{\|\nu^{-1}wg\|_{L^\infty_{x,v}} + |wr|_{L^\infty(\g_-)}\}.
\end{equation}
Moreover, if $g$ is continuous in $\Omega_d\times\mathbb{R}$  and $r(v)$ is continuous in $\R^3_+$,
then $f$ is continuous away from the grazing set $\gamma_0$.
\end{lemma}

\noindent{\bf Proof.} Let $f^\v$ be the solution of \eqref{S3.3} constructed in Lemma \ref{lemS3.6} for $\v\in(0,1]$.  Multiplying  $\eqref{S3.3}_1$ by $f^\v$ and integrating it over $\Omega_d\times\mathbb{R}^3$, we have
\begin{align}\label{S3.109}
&\v\|f^\v\|^2_{L^2_{x,v}}+\f12 |(I-P_\g)f^\v(0)|^2_{L^2(\g_+)} +c_0p_E^0 \|(\mathbf{I-P})f^\v\|^2_{\nu}\nonumber\\
&\leq \delta |P_\g f^\v(0)|^2_{L^2(\g_+)}+C_\delta |r|^2_{L^2(\g_-)}+\|g\|_{L^2_{x,v}} \|f^\v\|_{L^2_{x,v}},
\end{align}
which, together with  Lemma \ref{lemS3.5}, yields that
\begin{align}\label{S3.110}
&\f{1}{d^6} \|f^\v\|^2_{L^2_{x,v}} + \|({\bf I-P})f^\v\|^2_{\nu} +  |(I-P_\g)f^\v(0)|^2_{L^2(\g_+)} \nonumber\\
&\leq C \delta |P_\g f^\v(0)|^2_{L^2(\g_+)} + C d^6 \|g\|^2_{L^2_{x,v}} + C_\delta |r|^2_{L^2(\g_+)}.
\end{align}
By similar arguments as \eqref{2.76}-\eqref{2.79}, one can obtain
\begin{equation}\label{S3.110-1}
|P_\g f^\v(0)|^2_{L^2(\g_+)} \leq \frac{C}{\epsilon^2} \Big\{ \|f^\v\|_{L^2_{x,v}}^2 + \|({\bf I-P})f^\v\|^2_{\nu} + |(I-P_\g)f^\v(0)|^2_{L^2(\g_+)}+ \|g\|^2_{L^2_{x,v}}\Big\}.
\end{equation}
Substituting \eqref{S3.110-1} into \eqref{S3.110}, then taking $\delta$ suitably small, one has
\begin{equation}\label{3.94}
\|f^\v\|^2_{L^2_{x,v}} + \|({\bf I-P})f^\v\|^2_{\nu} +  |(I-P_\g)f^\v(0)|^2_{L^2(\g_+)}
\leq C_{\epsilon,d} \Big\{ \|g\|^2_{L^2_{x,v}} + |r|^2_{L^2(\g_-)}\Big\}.
\end{equation}
Now applying \eqref{3.14} to $f^\v$ and using \eqref{3.94}, then we obtain
\begin{equation}\label{S3.119}
\|wf^\v\|_{L^\infty_{x,v}}+|wf^\vep|_{L^\infty(\g)}\leq C_{\epsilon,d} \Big\{\|\nu^{-1}wg\|_{L^\infty_{x,v}} + |wr|_{L^\infty(\g_-)} \Big\}.
\end{equation}

\smallskip

Next we consider the convergence of $f^\v$ as $\v\rightarrow0+$. For any $\v_1,\v_2>0$, we consider the difference $f^{\v_2}-f^{\v_1}$ satisfying
\begin{equation}\label{S3.120}
\begin{cases}
v_3\partial_x (f^{\v_2}-f^{\v_1})+ \FL(f^{\v_2}-f^{\v_1})=-\v_2 f^{\v_2}+\v_1 f^{\v_1},\\[2mm]
(f^{\v_2}-f^{\v_1})(0,v)|_{v_3>0}= P_\g (f^{\v_2}-f^{\v_1})(0,v),\\[2mm]
(f^{\v_2}-f^{\v_1})(d,v)|_{v_3<0}=(f^{\v_2}-f^{\v_1})(d,\Mr v)
\end{cases}
\end{equation}
Multiplying \eqref{S3.120} by $f^{\v_2}-f^{\v_1}$,  by similar arguments as in \eqref{S3.109}-\eqref{3.94}, and using \eqref{S3.119}, one gets
\begin{align}\label{S3.121}
&\|f^{\v_2}-f^{\v_1}\|^2_{L^2_{x,v}} + \|({\bf I-P})(f^{\v_2}-f^{\v_1})\|^2_{\nu} +  |(I-P_\g) (f^{\v_2}-f^{\v_1})(0)|^2_{L^2(\g_+)}\nonumber\\
&\leq C_{\epsilon,d} \|\v_2 f^{\v_2}-\v_1 f^{\v_1}\|^2_{L^2_{x,v}} \leq C_{\epsilon,d} (\v_1+\v_2)^2 \cdot \Big\{\|\nu^{-1}wg\|^2_{L^\infty_{x,v}} + |wr|^2_{L^\infty(\g_-)} \Big\} \to 0,
\end{align}
as $\v_1$, $\v_2\rightarrow 0+$. Finally, applying \eqref{3.14} to $f^{\v_2}-f^{\v_1}$ and using \eqref{S3.121}, then we obtain
\begin{align}\label{S3.122}
&\|w(f^{\v_2}-f^{\v_1})\|_{L^\infty_{x,v}}+| w(f^{\v_2}-f^{\v_1})|_{L^\infty(\g)}\nonumber\\
&\leq C\Big\{ \|\nu^{-1}w (\v_2 f^{\v_2}-\v_1 f^{\v_1})\|_{L^\infty_{x,v}} +\|f^{\v_2}-f^{\v_1}\|_{L^2_{x,v}} \Big\}\nonumber\\
&\leq C_{\epsilon,d} \cdot (\v_1+\v_2)\Big\{\|\nu^{-1}wg\|_{L^\infty_{x,v}} + |wr|_{L^\infty(\g_-)} \Big\}  \rightarrow 0,
\end{align}
as $\v_1$, $\v_2\rightarrow 0+$,
 With \eqref{S3.122}, we know that there exists a function $f$ so that $\|w(f^{\v}-f)\|_{L^\infty_{x,v}}\rightarrow0$ as $\v\rightarrow 0+$. And it is  direct to see that $f$ solves \eqref{S3.106}. Also,  \eqref{S3.107} follows immediately from \eqref{S3.119}. The continuity of $f$
 follows directly from the $L^\infty_{x,v}$-convergence and the continuity of $f^\v$. Therefore the proof of Lemma \ref{lemS3.7} is complete. $\hfill\Box$

\

To obtain the solution for half-space problem, we need some uniform estimates independent of $d$, then we can take the limit $d\rightarrow\infty$. Let $f(x,v)$ be the solution of \eqref{S3.106}, we denote
\begin{align*}
\mathbf{P} f(x,v)=\big[a(x)+b(x)\cdot v+c(x) (\f{|v|^2}{2}-\f32)\big] \sqrt{\mu}.
\end{align*}
It follows from \eqref{3.57} that
\begin{equation}\label{S3.123}
\int_0^d\int_{\R^3} \sqrt{\mu(v)} f(x,v) dvdx=\int_0^d a(x) dx=0.
\end{equation}
Noting  $\dis z_{\g_+}(f):=\int_{v_3<0} |v_3| \sqrt{\mu(v)} f(0,v) dv,$
we define
\begin{equation}\label{S3.137-1}
\bar{f}(x,v)=f(x,v)-\sqrt{2\pi \mu(v)}\, z_{\g_+}(f),
\end{equation}
then it is easy to check that $\bar{f}$ satisfies
\begin{align}\label{S3.123-5}
\begin{cases}
v_3 \partial_x\bar{f}+p_E^0 \FL \bar{f}=g,\quad (x,v)\in \Omega_d\times \R^3, \\
\bar{f}(0,v)|_{v_3>0}=P_\g \bar{f}(0,v) + r(v),\\
\bar{f}(d,v)|_{v_3<0}=\bar{f}(d,\Mr v).
\end{cases}
\end{align}
It is direct to check that
\begin{equation}\label{S3.123-6}
z_{\g_{+}}(\bar{f})=\int_{v_3<0} |v_3| \sqrt{\mu(v)} \bar{f}(0,v) dv=0.
\end{equation}
Moreover it follows from \eqref{S3.107} and \eqref{S3.137-1} that
\begin{align}\label{S3.140-2}
\|w\bar{f}\|_{L^\infty_{x,v}} + |w\bar{f}|_{L^\infty(\g)} \leq C_d \{\|\nu^{-1}wg\|_{L^\infty_{x,v}} + |wr|_{L^\infty(\g_-)}\}.
\end{align}

\

For later use, we denote
\begin{align*}
\mathbf{P} \bar{f}(x,v)=\big[\bar{a}(x)+\bar{b}(x)\cdot v+\bar{c}(x) (\f{|v|^2}{2}-\f32)\big] \sqrt{\mu}.
\end{align*}
From now on, we assume $g\in \mathcal{N}^{\perp}$, then multiplying \eqref{S3.123-5} by $\sqrt{\mu}$, one has that
\begin{equation}\label{S3.123-1}
\f{d}{dx} \int_{\R^3} v_3\sqrt{\mu} \bar{f}(x,v) dv= \int_{\R^3} g(x,v) \sqrt{\mu}dv\equiv0.
\end{equation}
It follows from the specular reflection boundary condition that
\begin{equation*}
\int_{\R^3} v_3 \sqrt{\mu} \bar{f}(d,v)dv=0,
\end{equation*}
which, together with \eqref{S3.123-1}, yields
\begin{equation}\label{S3.123-2}
\bar{b}_3(x)=0,\quad \mbox{for}\quad x\in [0,d].
\end{equation}
Similarly, multiplying \eqref{S3.123-5} by $v_i\sqrt{\mu}, i=1,2,$ and $(|v|^2-5)\sqrt{\mu}$, we can obtain, for $x\in[0,d]$, that
\begin{align}\label{S3.123-4}
\begin{split}
0=&\int_{\R^3} v_3v_i \sqrt{\mu} \bar{f}(x,v) dv=\int_{\R^3} v_3v_i \sqrt{\mu} ({\bf I-P})\bar{f}(x,v) dv,\\
0=&\int_{\R^3} v_3 (|v|^2-5) \sqrt{\mu} \bar{f}(x,v) dv=\int_{\R^3} v_3 (|v|^2-5) \sqrt{\mu} ({\bf I-P})\bar{f}(x,v) dv,
\end{split}
\end{align}
Utilizing \eqref{S3.123-2} and \eqref{S3.123-4}, a direct calculation shows that
\begin{align}\label{S3.124}
\int_{\R^3} v_3 |{\bf P} \bar{f}(x,v)|^2dv=0\quad\mbox{and}\quad
\int_{\R^3} v_3 {\bf P} \bar{f}(x,v) \, ({\bf I-P})\bar{f}(x,v)dv=0,\quad \forall x\in[0,d],
\end{align}
which implies that
\begin{align}\label{S3.125}
\int_{\R^3} v_3 |\bar{f}(x,v)|^2 dv= \int_{\R^3} v_3 |({\bf I-P})\bar{f}(x,v)|^2 dv,\quad \forall x\in[0,d].
\end{align}

\

\begin{lemma}\label{lem2.11}
For the solution of \eqref{S3.123-5}, it holds that
\begin{align}\label{S3.126}
&|(I-P_\g)\bar{f}(0)|^2_{L^2(\g_+)} + \int_0^d e^{2\sigma_1 x} \|({\bf I-P})\bar{f}(x)\|_{\nu}^2dx\nonumber\\
&\leq C\Big\{|r|^2_{L^2(\g_-)} + \int_0^d e^{2\sigma_1 x} \|g(x)\|^2_{L^2_{x,v}} dx \Big\}\quad\mbox{for}\quad \sigma_1\in[0,\s_0].
\end{align}
\end{lemma}

\noindent{\bf Proof.} Multiplying \eqref{S3.123-5} by $e^{2\sigma_1 x} \bar{f}(x,v)$, one has
\begin{align*}
\f{d}{dx}\Big\{e^{2\sigma_1 x} \int_{\R^3} v_3 |\bar{f}(x,v)|^2 dx\Big\}+ (c_0p_E^0-C\sigma_1) e^{2\sigma_1 x} \|({\bf I-P})\bar{f}(x)\|_{\nu}^2
\leq C e^{2\sigma_1 x} \|g(x)\|^2_{L^2},
\end{align*}
where we have used \eqref{S3.125}.
Integrating above inequality over $[0,d]$, using \eqref{S3.125} and taking $\sigma_1>0$ suitably small, one obtains
\begin{align}\label{S3.127}
-\int_{\R^3}  v_3 |\bar{f}(0,v)|^2 dv+\f12 c_0p_E^0 \int_0^d e^{2\sigma_1 x} \|({\bf I-P})\bar{f}(x)\|_{\nu}^2dx\leq C\int_0^d e^{2\sigma_1 x} \|g(x)\|^2_{L^2} dx,
\end{align}
where we have used the fact $\int_{\R^3} v_3|\bar{f}(d,v)|^2dv=0$ due to the specular boundary condition.

Using the diffuse reflection boundary condition, a direct calculation shows
\begin{align*}%\label{S3.131}
-\int_{\R^3}  v_3 |\bar{f}(0,v)|^2 dv&=|(I-P_\g)\bar{f}(0)|^2_{L^2(\g_+)}-|r|^2_{L^2(\g_-)}-2\int_{v_3>0}z_{\g_+}(\bar{f}) \sqrt{2\pi \mu(v)} v_3 r(v)dv\nonumber\\
&=|(I-P_\g)\bar{f}(0)|^2_{L^2(\g_+)}-|r|^2_{L^2(\g_-)},
\end{align*}
which, together with \eqref{S3.127}, yields \eqref{S3.126}. Therefore the proof is complete. $\hfill\Box$

\

To obtain the uniform in $d$ estimate for macroscopic term, motivated by \cite{GPS-1988,HJW}, we define
\begin{equation}\label{S3.132}
\tilde{f}(x,v):=\bar{f}(x,v)+\Phi(v)
\end{equation}
with
\begin{equation*}
\Phi(v):=\big[\phi_0+\phi_1 v_1+\phi_2v_2+\phi_3(\frac{|v|^2}{2}-\f32)\big] \sqrt{\mu(v)},
\end{equation*}
where $\phi_0, \phi_1, \phi_2,\phi_3$ are four constants determined later. Clearly, $\tilde{f}$ satisfies
\begin{align}\label{S3.133}
\begin{cases}
v_3 \partial_x\tilde{f}+p_E^0 \FL \tilde{f}=g,\quad (x,v)\in \Omega_d\times \R^3, \\
\tilde{f}(0,v)|_{v_3>0}=P_\g \tilde{f}(0,v) + (I-P_\g)\Phi + r(v),\\
\tilde{f}(d,v)|_{v_3<0}=\tilde{f}(x,\Mr v),
\end{cases}
\end{align}
Also we have
\begin{align}\label{S3.134}
\int_{v_3>0} \{(I-P_\g)\Phi+r\}v_3\sqrt{\mu(v)} dv=0,
\end{align}
and
\begin{align}\label{S3.134-1}
\begin{split}
&{\bf P}\tilde{f}(x,v)=\big[\tilde{a}(x)+\tilde{b}_1(x)v_1+\tilde{b}_2(x)v_2+\tilde{c}(x)(\f{|v|^2}{2}-\f32)\big]\, \sqrt{\mu(v)},\\
&({\bf I-P})\tilde{f}(x,v)=({\bf I-P})\bar{f}(x,v)=({\bf I-P})f(x,v),
\end{split}
\end{align}
with
\begin{align}\label{S3.134-2}
\begin{cases}
\tilde{a}(x)=\bar{a}(x)+\phi_0=a(x)+\phi_0-\sqrt{2\pi} z_{\g_+}(f),\\
\tilde{b}_i(x)=\bar{b}_i(x)+\phi_i=b_i(x)+\phi_i, \, i=1,2, \\
\tilde{c}(x)=\bar{c}(x)+\phi_3=c(x)+\phi_3.
\end{cases}
\end{align}

For later use, we denote
\begin{align}\nonumber
\begin{split}
\mathcal{A}_{ij}(v)&=\{v_iv_j-\f{\d_{ij}}{3} |v|^2\}\sqrt{\mu(v)},\\
\mathcal{B}_i(v)&=\f12 v_i (|v|^2-5) \, \sqrt{\mu(v)}.
\end{split}
\end{align}
It is obvious to know $\mathcal{A}_{ij},\, \mathcal{B}_i \in \mathcal{N}^{\perp}$. We define
\begin{align}\label{S3.135}
\begin{split}
\kappa_1&:=\langle \mathcal{A}_{ij}, \, {\bf L}^{-1} \mathcal{A}_{ij} \rangle>0,\,\, i\neq j,\\
\kappa_2&:=\langle \mathcal{B}_{i}, \, {\bf L}^{-1} \mathcal{B}_{i} \rangle>0,
\end{split}
\qquad  i,j=1,2,3.
\end{align}

\begin{lemma}\label{lem2.12}
There exist constants $\phi_0,\phi_1,\phi_2,\phi_3$  such that
\begin{align}\label{S3.136}
\begin{split}
&\int_{\R^3} v_3\tilde{f}(d,v)\cdot v_3 \sqrt{\mu} dv=0,\\
&\int_{\R^3} v_3\tilde{f}(d,v)\cdot \FL^{-1} (\mathcal{A}_{3i}) dv=0,\  i=1,2,\\
&\int_{\R^3} v_3\tilde{f}(d,v)\cdot \FL^{-1} (\mathcal{B}_{3}) dv=0.
\end{split}
\end{align}
\end{lemma}

\noindent{\bf Proof.} Motivated by \cite{GPS-1988}, a direct calculation shows that
\begin{align}
	\int_{\R^3} v_3\tilde{f}(x,v)\cdot v_3 \sqrt{\mu} dv
	&=\tilde{a}(x) +\tilde{c}(x) +\int_{\R^3} \mathcal{A}_{33}(v)\cdot (\mathbf{I-P})\tilde{f}(x,v) dv\nonumber\\
	&=\phi_0+\phi_3 +\bar{a}(x)+\bar{c}(x)+\int_{\R^3} \mathcal{A}_{33}(v)\cdot (\mathbf{I-P})\bar{f}(x,v) dv,\label{S3.137}\\
	\int_{\R^3} v_3\tilde{f}(x,v)\cdot \FL^{-1} (\mathcal{A}_{31}) dv
	&= \kappa_1 \tilde{b}_1(x)+\int_{\R^3} v_3 (\mathbf{I-P})\tilde{f}(x,v)\cdot \FL^{-1} (\mathcal{A}_{31})  dv\nonumber\\
	&=\kappa_1 \phi_1+\kappa_1 \bar{b}_1(x) +\int_{\R^3} v_3 (\mathbf{I-P})\bar{f}(x,v)\cdot \FL^{-1} (\mathcal{A}_{31})  dv,\label{S3.138} \\
	\int_{\R^3} v_3\tilde{f}(x,v)\cdot \FL^{-1} (\mathcal{A}_{32}) dv
	&=\kappa_1 \tilde{b}_2(x)+\int_{\R^3} v_3 (\mathbf{I-P})\tilde{f}(x,v)\cdot \FL^{-1} (\mathcal{A}_{32})  dv\nonumber\\
	&=\kappa_1 \phi_2+\kappa_1 \bar{b}_2(x) +\int_{\R^3} v_3 (\mathbf{I-P})\bar{f}(x,v)\cdot \FL^{-1} (\mathcal{A}_{32})  dv,\label{S3.139}\\
	\int_{\R^3} v_3\tilde{f}(x,v)\cdot \FL^{-1} (\mathcal{B}_{3}) dv
	&=\kappa_2 \tilde{c}(x)+\int_{\R^3} v_3 (\mathbf{I-P})\tilde{f}(x,v)\cdot \FL^{-1} (\mathcal{B}_{3})  dv\nonumber\\
	&=\kappa_2 \phi_3+\kappa_2 \bar{c}(x)+\int_{\R^3} v_3 (\mathbf{I-P})\bar{f}(x,v)\cdot \FL^{-1} (\mathcal{B}_{3})  dv,\label{S3.140}
\end{align}
where we have used $\eqref{S3.134-1}_2$. Using  \eqref{S3.137}-\eqref{S3.140}, then \eqref{S3.136} is equivalent  as
\begin{equation}\label{S3.140-1}
	\left(
	\begin{array}{cccc}
		1 & 0 & 0 & 1  \\
		0 & \kappa_1 & 0 & 0  \\
		0 & 0 & \kappa_1 & 0  \\
		0 & 0 & 0 & \kappa_2  \\
	\end{array}
	\right)
	\left(
	\begin{array}{c}
		\phi_0 \\
		\phi_1\\
		\phi_2\\
		\phi_3\\
	\end{array}
	\right)
	=-\left(
	\begin{array}{c}
	\dis	\bar{a}(d)+\bar{c}(d)+\int_{\R^3}(\mathbf{I-P})\bar{f}(d,v) \cdot \mathcal{A}_{33}(v)dv \\[2mm]
	\dis\kappa_1 \bar{b}_1(d) +\int_{\R^3} v_3 (\mathbf{I-P})\bar{f}(d,v)\cdot \FL^{-1} (\mathcal{A}_{31})  dv\\[2mm]
	\dis\kappa_1 \bar{b}_2(d) +\int_{\R^3} v_3 (\mathbf{I-P})\bar{f}(d,v)\cdot \FL^{-1} (\mathcal{A}_{32})  dv\\[2mm]
	\dis\kappa_2 \bar{c}(d)+\int_{\R^3} v_3 (\mathbf{I-P})\bar{f}(d,v)\cdot \FL^{-1} (\mathcal{B}_{3})  dv
	\end{array}
	\right).
\end{equation}
Noting the matrix is non-singular, hence  $(\phi_0,\phi_1,\phi_2,\phi_3)$  is found. Therefore the proof of Lemma \ref{lem2.12} is complete. $\hfill\Box$

\smallskip

\begin{lemma}\label{lem2.13}
	Let $\phi_0,\phi_1,\phi_2,\phi_3$ be the constants determined  in Lemma \ref{lem2.12}, then it holds that
	\begin{equation}\label{S3.141}
		\|e^{\s x}\tilde{f}\|_{L^2_{x,v}}\leq C\left\{|r|_{L^2(\g_-)}+\frac{1}{\s_1-\s}\|e^{\s_1 x} g\|_{L^2_{x,v}}\right\},
	\end{equation}
	with $0<\s<\s_1\leq \s_0$, and the constant $C>0$ is independent of $d$. It is important that the right hand side of \eqref{S3.141} is independent of $\Phi$.
\end{lemma}

\noindent{\bf Proof.}  It follows from \eqref{S3.123-4} that
\begin{align}
\begin{split}\label{S3.142}
\int_{\R^3} \FL (\mathbf{I-P})\bar{f} \cdot \FL^{-1} (\mathcal{A}_{3i}) dv
&=\int_{\R^3}  (\mathbf{I-P})\bar{f} \cdot \mathcal{A}_{3i} dv=0,\,\, i=1,2,\\
\int_{\R^3} \FL (\mathbf{I-P})\bar{f}\cdot \FL^{-1} (\mathcal{B}_{3}) dv
&=\int_{\R^3}  (\mathbf{I-P})\bar{f}\cdot \mathcal{B}_{3} dv=0.
\end{split}
\end{align}
Multiplying \eqref{S3.133} by  $\FL^{-1} (\mathcal{A}_{31}), \FL^{-1} (\mathcal{A}_{32}) $ and $\FL^{-1} (\mathcal{B}_{3}) $, respectively, and using \eqref{S3.142},  we can obtain
\begin{align}
	\frac{d}{dx}\left(
	\begin{array}{c}
	\dis	\int_{\R^3} v_3\tilde{f}(x,v)\cdot \FL^{-1} (\mathcal{A}_{31}) dv\\[1mm]
	\dis	\int_{\R^3} v_3\tilde{f}(x,v)\cdot \FL^{-1} (\mathcal{A}_{32}) dv\\[1mm]
	\dis	\int_{\R^3} v_3\tilde{f}(x,v)\cdot \FL^{-1} (\mathcal{A}_{3}) dv\\
	\end{array}
	\right)
	&=
	\left(
	\begin{array}{c}
	\dis	\int_{\R^3}  [g-p_E^0 \FL (\mathbf{I-P})\bar{f}]\cdot \FL^{-1} (\mathcal{A}_{31}) dv\\[1mm]
	\dis	\int_{\R^3} [g-p_E^0 \FL (\mathbf{I-P})\bar{f}]\cdot \FL^{-1} (\mathcal{A}_{32}) dv\\[1mm]
	\dis	\int_{\R^3} [g-p_E^0 \FL (\mathbf{I-P})\bar{f}]\cdot \FL^{-1} (\mathcal{B}_{3}) dv\\
	\end{array}
	\right) \nonumber\\
&	=\left(
	\begin{array}{c}
	\dis	\int_{\R^3}  g \FL^{-1} (\mathcal{A}_{31}) dv\\[1mm]
	\dis	\int_{\R^3} g  \FL^{-1} (\mathcal{A}_{32}) dv\\[1mm]
	\dis	\int_{\R^3} g \FL^{-1} (\mathcal{B}_{3}) dv\\
	\end{array}
	\right).\label{S3.143}
\end{align}
Integrating above system over $[x,d]$ and using $\eqref{S3.136}$, we get
\begin{align}\label{S3.143-1}
	\left(
	\begin{array}{c}
	\dis	\int_{\R^3} v_3\tilde{f}(x,v)\cdot \FL^{-1} (\mathcal{A}_{31}) dv\\[1mm]
	\dis	\int_{\R^3} v_3\tilde{f}(x,v)\cdot \FL^{-1} (\mathcal{A}_{32}) dv\\[1mm]
	\dis	\int_{\R^3} v_3\tilde{f}(x,v)\cdot \FL^{-1} (\mathcal{B}_{3}) dv\\
	\end{array}
	\right)
	=-\int_x^d \left(
	\begin{array}{c}
	\dis	\int_{\R^3}  g \FL^{-1} (\mathcal{A}_{31}) dv\\[1mm]
	\dis	\int_{\R^3} g  \FL^{-1} (\mathcal{A}_{32}) dv\\[1mm]
	\dis	\int_{\R^3} g \FL^{-1} (\mathcal{B}_{3}) dv\\
	\end{array}
	\right)(z)dz,
\end{align}
which, together with \eqref{S3.138}-\eqref{S3.140}, yields that
\begin{align}\label{S3.144}
	|(\kappa_1 \tilde{b}_{1},\kappa_1  \tilde{b}_{2},\kappa_2 \tilde{c})(x)|\leq C\|(\mathbf{I-P})\bar{f}(x)\|_{\nu}+C\int_x^d \|g(z)\|_{L^2_v} dz.
\end{align}
Multiplying \eqref{S3.144} by $e^{\s x}$ with $0<\s<\s_1\leq\s_0$ and using \eqref{S3.126}, one has
\begin{align}\label{S3.145}
	\int_0^d e^{2\s x}|(\tilde{b}_{1},\tilde{b}_{2},\tilde{c})(x)|^2d x
	&\leq  C\int_0^d e^{2\s x} \|(\mathbf{I-P})\bar{f}(x)\|_{\nu}^2 dx\nonumber\\
	&\quad +C\int_0^d e^{2\s x} \left|\int_x^d \|g(z)\|_{L^2_v}dz\right|^2 dx\nonumber\\
	&\leq C\left\{ |r|^2_{L^2(\g_-)} + \frac{1}{\s_1-\s}\|e^{\s_1 x} g\|_{L^{2}_{x,v}}^2\right\}.
\end{align}

Finally we consider the case for $\tilde{a}$. In fact, multiplying \eqref{S3.133} by $v_3\sqrt{\mu}$, we get
\begin{align}\label{S3.146}
\frac{d}{dx} \int_{\R^3} \tilde{f}(x,v)\cdot v_3^2 \sqrt{\mu} dv=\int_{\R^3}g\cdot v_3\sqrt{\mu}  dv=0.
\end{align}
Integrating above equation over $[x,d]$,  and using  \eqref{S3.136}-\eqref{S3.137} , one obtains
\begin{align}\label{S3.148}
\tilde{a}(x)=-\tilde{c}(x)+\int_{\R^3}   (\mathbf{I-P})\bar{f}(x,v)\cdot v_3^2 \sqrt{\mu} dv.
\end{align}
Multiplying \eqref{S3.148} by $e^{\s x}$ with $0<\s<\s_1\leq \s_0$ and using \eqref{S3.126}, \eqref{S3.145}, it holds that
\begin{equation}\label{S3.149}
	\int_0^d e^{2\s x}|\tilde{a}(x)|^2d x\leq  C\left\{ |r|^2_{L^2(\g_-)} + \frac{1}{\s_1-\s}\|e^{\s_1 x} g\|_{L^{2}_{x,v}}^2\right\}.
\end{equation}
Combining \eqref{S3.149}, \eqref{S3.145} and  \eqref{S3.126}, we prove \eqref{S3.141}. Therefore the proof of Lemma \ref{lem2.13} is complete. $\hfill\Box$

\

\begin{lemma}\label{lem2.16}
Let $\tilde{f}$ the solution of \eqref{S3.133}, then it holds that
\begin{align}\label{S3.166}
&\|e^{\s x}w\tilde{f}\|_{L^\infty_{x,v}}+|e^{\s x}w\tilde{f}|_{L^\infty(\g)}\nonumber\\
&\leq C\Big\{\frac{1}{\s_0-\s}\|e^{\s_0 x}\nu^{-1}wg\|_{L^\infty_{x,v}} + |(\phi_0,\phi_1,\phi_2,\phi_3)|+ |wr|_{L^\infty(\g_-)}  \Big\}.
\end{align}
\end{lemma}

\noindent{\bf Proof.} Let $\tilde{h}:=e^{\s x} w \tilde{f}$. Multiplying \eqref{S3.133} by $e^{\s x}w $ to have
\begin{equation*}%\label{3.129-1}
\begin{cases}
v_3 \partial_x\tilde{h}+p_E^0 \nu_{\s}(v) \tilde{h}=p_E^0 K_w \tilde{h}+e^{\s x}wg,\\
\tilde{h}(0,v)|_{v_3>0}=w(v)P_\g \tilde{f}(0,v) + w(v)\{(I-P_\g)\Phi + r(v)\},\\
\tilde{h}(d,v)|_{v_3<0}=\tilde{h}(d,\Mr v),
\end{cases}
\end{equation*}
where $\nu_{\s}(v):=\nu(v) -\s v_3$. We further take $\s_0>0$ small such that $\nu_{\s}(v) \geq \frac12 \nu(v)>0$. By the same arguments as  in Lemma \ref{lemS3.3} and using \eqref{S3.141}, we can obtain
\begin{align*}
&\|\tilde{h}\|_{L^\infty_{x,v}} +|\tilde{h}|_{L^\infty(\gamma)} \nonumber\\
& \leq C\Big\{C\|e^{\s x} \tilde{f}\|_{L^2_{x,v}} + \|e^{\s x}\nu^{-1} wg\|_{L^\infty_{x,v}} + |(\phi_0,\phi_1,\phi_2,\phi_3)|+ |wr|_{L^\infty(\g_-)} \Big\}\nonumber\\
&\leq C\Big\{\frac{C}{\s_1-\s}\|e^{\s_1 x} g\|_{L^2_{x,v}}+ \|e^{\s x}\nu^{-1} wg\|_{L^\infty_{x,v}}+  |(\phi_0,\phi_1,\phi_2,\phi_3)|+ |wr|_{L^\infty(\g_-)} \Big\}\nonumber\\
&\leq C\Big\{\frac{1}{\s_0-\s}\|e^{\s_0 x}\nu^{-1}wg\|_{L^\infty_{x,v}} + |(\phi_0,\phi_1,\phi_2,\phi_3)|+ |wr|_{L^\infty(\g_-)}  \Big\},
\end{align*}
where we have chosen $\s_1=\s+\frac{\s_0-\s}{2}$ such that  $0<\s<\s_1<\s_0$. Hence the proof of Lemma \ref{lem2.16} is completed. $\hfill\Box$

\begin{remark}
From \eqref{S3.140-1}, we know that the constants $(\phi_0,\phi_1,\phi_2, \phi_3)$ depend on $d$, hence the right hand side of \eqref{S3.166} still depends on $d$. In the following, we shall write down them as $(\phi_0,\phi_1,\phi_2, \phi_3)(d)$ to emphasize the dependence of $d$. To obtain the uniform estimate for $\tilde{f}$, we   need to establish the uniform in $d$ estimate  for $(\phi_0,\phi_1,\phi_2, \phi_3)(d)$. This is the key part of the present paper. % To take the limit $d\to +\infty$, we need first to study the limit of $(\phi_0,\phi_1,\phi_2, \phi_3)(d)$ as $d\to \infty$.
\end{remark}

\medskip

Hereafter we also denote the solutions of \eqref{S3.106}, \eqref{S3.123-5} and \eqref{S3.133} as $f_d, \bar{f}_d$ and $\tilde{f}_d$, respectively,  to emphasize the dependence of $d$. For later use, we denote
\begin{equation}
{\bf P}\bar{f}_d(x,v)=\big\{\bar{a}_d(x)+\bar{b}_{d,1}(x)v_1+ \bar{b}_{d,2}v_2+ \bar{c}_{d}(x) (\f{|v|^2}{2}-\f32)\big\} \sqrt{\mu(v)}.
\end{equation}

\smallskip

Firstly we give a  formula of  $(\bar{a}_d,\bar{b}_{d,1},\bar{b}_{d,2}, \bar{c}_d)(x)$ by using the  boundary data at $x=0$.
\begin{lemma}\label{lem2.14}
It holds that
\begin{align}
\kappa_1 \, \bar{b}_{d,i}(x)&=-\Big\langle (\mathbf{I-P})\bar{f}_d(x,v),\,  v_3\FL^{-1} (\mathcal{A}_{3i}) \Big\rangle+\int_0^x\int_{\R^3}  g(z) \FL^{-1} (\mathcal{A}_{3i}) dvdz\nonumber\\
&\quad+ \int_{\g_+} (I-P_\g)\bar{f}_d(0,v)\cdot v_3\FL^{-1} (\mathcal{A}_{3i})dv + \int_{\g_-}v_3\FL^{-1} (\mathcal{A}_{3i}) \cdot r(v)dv,\,\, i=1,2,\label{S3.151}\\
\kappa_2\, \bar{c}_d(x)&=-\Big\langle (\mathbf{I-P})\bar{f}_d(x,v),\,  v_3\FL^{-1} (\mathcal{B}_{3}) \Big\rangle+\int_0^x\int_{\R^3} g(z) \FL^{-1} (\mathcal{B}_{3}) dvdz\nonumber\\
&\qquad+\int_{\g_+} (I-P_\g)\bar{f}_d(0,v)\cdot v_3\FL^{-1} (\mathcal{B}_3)dv + \int_{\g_-}v_3\FL^{-1} (\mathcal{B}_3) \cdot r(v)dv,\label{S3.151-1}\\
\bar{a}_d(x)&=-\Big\langle (\mathbf{I-P})\bar{f}_d(x,v),\,  \mathcal{A}_{33}-\f{1}{\kappa_2}v_3\FL^{-1} (\mathcal{B}_{3})\Big\rangle -\f{1}{\kappa_2}\int_0^x\int_{\R^3} g(z) \FL^{-1} (\mathcal{B}_{3}) dvdz\nonumber\\
&\qquad +\int_{\g_+} (I-P_\g)\bar{f}_d(0,v)\cdot v_3 \big\{v_3\sqrt{\mu(v)}-\f{1}{\kappa_2}\FL^{-1} (\mathcal{B}_{3})\big\}dv\nonumber\\
&\qquad+ \int_{\g_-}   r(v)\cdot v_3 \big\{v_3\sqrt{\mu(v)}-\f{1}{\kappa_2}\FL^{-1} (\mathcal{B}_{3})\big\} dv.\label{S3.151-2}
\end{align}
\end{lemma}

\noindent{\bf Proof.} Multiplying $\eqref{S3.123-5}_1$ by $\FL^{-1} (\mathcal{A}_{31}), \FL^{-1} (\mathcal{A}_{32})$ and $\FL^{-1} (\mathcal{B}_{3})$, respectively, and integrating over $[0,x]$ and using \eqref{S3.142}, one gets
\begin{align}\nonumber
\begin{split}
\int_{\R^3} v_3\bar{f}_d(x,v)\cdot \FL^{-1} (\mathcal{A}_{31}) dv
&=\int_{\R^3} v_3\bar{f}_d(0,v)\cdot \FL^{-1} (\mathcal{A}_{31}) dv+\int_0^x\int_{\R^3}  g(z) \FL^{-1} (\mathcal{A}_{31}) dvdz,\\
\int_{\R^3} v_3\bar{f}_d(x,v)\cdot \FL^{-1} (\mathcal{A}_{32}) dv
&=\int_{\R^3} v_3\bar{f}_d(0,v)\cdot \FL^{-1} (\mathcal{A}_{32}) dv + \int_0^x\int_{\R^3} g(z)  \FL^{-1} (\mathcal{A}_{32}) dvdz,\\
\int_{\R^3} v_3\bar{f}_d(x,v)\cdot \FL^{-1} (\mathcal{B}_{3}) dv
&=\int_{\R^3} v_3\bar{f}_d(0,v)\cdot \FL^{-1} (\mathcal{B}_{3}) dv+\int_0^x\int_{\R^3} g(z) \FL^{-1} (\mathcal{B}_{3}) dvdz,
\end{split}
\end{align}
which, together with a tedious calculation, gives
\begin{align}\label{S3.152-1}
\kappa_1\, \bar{b}_{d,1}(x)&=-\Big\langle (\mathbf{I-P})\bar{f}_d(x,v),\,  v_3\FL^{-1} (\mathcal{A}_{31}) \Big\rangle+\int_{\R^3} v_3\bar{f}_d(0,v)\cdot \FL^{-1} (\mathcal{A}_{31}) dv\nonumber\\
&\quad+\int_0^x\int_{\R^3}  g(z) \FL^{-1} (\mathcal{A}_{31}) dvdz,
\end{align}
\begin{align}\label{S3.152-2}
\kappa_1\, \bar{b}_{d,2}(x)&=-\Big\langle (\mathbf{I-P})\bar{f}_d(x,v),\,  v_3\FL^{-1} (\mathcal{A}_{32}) \Big\rangle+\int_{\R^3} v_3\bar{f}_d(0,v)\cdot \FL^{-1} (\mathcal{A}_{32}) dv\nonumber \\
&\quad+ \int_0^x\int_{\R^3} g(z)  \FL^{-1} (\mathcal{A}_{32}) dvdz,
\end{align}
and
\begin{align}\label{S3.152}
\kappa_2 \bar{c}_{d}(x)&=-\Big\langle (\mathbf{I-P})\bar{f}_d(x,v),\, v_3 \FL^{-1} (\mathcal{B}_{3}) \Big\rangle+\int_{\R^3} v_3\bar{f}_d(0,v)\cdot \FL^{-1} (\mathcal{B}_{3}) dv\nonumber\\
&\quad+\int_0^x\int_{\R^3} g(z) \FL^{-1} (\mathcal{B}_{3}) dvdz.
\end{align}
For the expression of $\bar{a}(x)$, multiplying $\eqref{S3.123-5}_1$ by $v_3\sqrt{\mu}$, and integrating  over $[0,x]$, one obtains
\begin{equation}\label{S3.153}
\bar{a}_d(x)=-\bar{c}_d(x)-\Big\langle (\mathbf{I-P})\bar{f}_d(x,v),\,  \mathcal{A}_{33}\Big\rangle+\int_{\R^3}\bar{f}_d(0,v)v_3^2\sqrt{\mu(v)} dv.
\end{equation}

\smallskip

It follows from \eqref{S3.123-6}  that
\begin{equation}\label{S3.154}
\bar{f}_d(0,v)= I_{\g_+}\cdot (I-P_\g) \bar{f}_d(0,v) + I_{\g_-}\cdot r(v).
\end{equation}
Noting \eqref{S3.123-6} and  \eqref{S3.154}, a direct calculation shows
\begin{align}
\begin{split}\label{S3.154-1}
\int_{\R^3} v_3\bar{f}_d(0,v)\cdot \FL^{-1} (\mathcal{A}_{3i}) dv
&=\int_{\g_+} (I-P_\g)\bar{f}_d(0,v)\cdot v_3\FL^{-1} (\mathcal{A}_{3i})dv\\
&\qquad + \int_{\g_-}v_3\FL^{-1} (\mathcal{A}_{3i}) \cdot r(v)dv,\, i=1,2,\\
\int_{\R^3} v_3\bar{f}_d(0,v)\cdot \FL^{-1} (\mathcal{B}_3) dv
&=\int_{\g_+} (I-P_\g)\bar{f}_d(0,v)\cdot v_3\FL^{-1} (\mathcal{B}_3)dv\\
&\qquad + \int_{\g_-}v_3\FL^{-1} (\mathcal{B}_3) \cdot r(v)dv,
\end{split}
\end{align}
and
\begin{align}\label{S3.154-2}
\int_{\R^3}\bar{f}_d(0,v)v_3^2\sqrt{\mu(v)} dv
=\int_{\g_+} (I-P_\g)\bar{f}_d(0,v)\cdot v_3^2\sqrt{\mu(v)}dv + \int_{\g_-}   r(v)\cdot v_3^2\sqrt{\mu(v)}dv,
\end{align}
which, together with \eqref{S3.152-1}-\eqref{S3.153}, yields \eqref{S3.151}-\eqref{S3.151-2}. Therefore the proof of Lemma \ref{lem2.14} is complete. $\hfill\Box$

\medskip

\begin{lemma}\label{lem2.17}
For $d\geq 1$, it holds that
\begin{align}\label{S3.180}
|(\phi_0,\phi_1,\phi_2,\phi_3)(d)|\leq C\Big\{\|e^{\s_0 x}\nu^{-1}wg\|_{L^\infty_{x,v}} + |wr|_{L^\infty(\g_-)}  \Big\},
\end{align}
and
\begin{align}\label{S3.181}
\|e^{\s x}w\tilde{f}\|_{L^\infty_{x,v}}+|e^{\s x}w\tilde{f}|_{L^\infty(\g)}
\leq C\Big\{\frac{1}{\s_0-\s}\|e^{\s_0 x}\nu^{-1}wg\|_{L^\infty_{x,v}} +  |wr|_{L^\infty(\g_-)}  \Big\}.
\end{align}
where the constants $C>0$ are independent of $d$. Hence we have obtained the uniform in $d$ estimates for both $\tilde{f}$ and $(\phi_0,\phi_1,\phi_2,\phi_3)(d)$.
\end{lemma}

\noindent{\bf Proof.} Noting that \eqref{S3.181} follows directly from \eqref{S3.166} and \eqref{S3.180}, hence we need only to prove \eqref{S3.180}.

\smallskip

Noting $(\mathbf{I-P})\bar{f}_d(x,v)\equiv (\mathbf{I-P})\tilde{f}_d(x,v)$, it follows from \eqref{S3.151}-\eqref{S3.151-2} and \eqref{S3.126} that
\begin{align}\label{S3.182}
&|(\bar{a}_d,\bar{b}_{d,1},\bar{b}_{d,2},\bar{c}_d)(d)|\nonumber\\
&\leq C \Big\{ |(\mathbf{I-P})\tilde{f}_d(d)|_{L^\infty_v} + \|g\|_{L^2_{x,v}} + |(I-P_\g)\bar{f}_d(0)|_{L^2(\g_+)} + |r|_{L^2(\g_-)}\Big\}\nonumber\\
&\leq C \Big\{ |(\mathbf{I-P})\tilde{f}_d(d)|_{L^\infty_v} + \|g\|_{L^2_{x,v}} + |r|_{L^2(\g_-)}\Big\}.
\end{align}
On the other hand, it follows from \eqref{S3.140-1} and \eqref{S3.182}   that
\begin{align}\label{S3.183}
&|(\phi_0,\phi_1,\phi_2,\phi_3)(d)|
\leq C |(\bar{a}_d,\bar{b}_{d,1},\bar{b}_{d,2},\bar{c}_d)(d)|+ C |(\mathbf{I-P})\bar{f}_d(d)|_{L^\infty_v}\nonumber\\
&\leq C|(\mathbf{I-P})\tilde{f}_d(d)|_{L^\infty_v} + C \Big\{\|g\|_{L^2_{x,v}}  + |r|_{L^2(\g_-)}\Big\}\nonumber\\
&\leq Ce^{-\sigma d} |e^{\sigma d} w\tilde{f}_d(d)|_{L^\infty_{v}} + C \Big\{\|g\|_{L^2_{x,v}}  + |r|_{L^2(\g_-)}\Big\}\nonumber\\
&\leq Ce^{-\sigma d} |(\phi_0,\phi_1,\phi_2,\phi_3)(d)| + C\Big\{\frac{1}{\s_0-\s}\|e^{\s_0 x}\nu^{-1}wg\|_{L^\infty_{x,v}} +  |wr|_{L^\infty(\g_-)}  \Big\},
\end{align}
where we have used \eqref{S3.166}  in the last inequality.

Taking $\s=\f12\s_0$, and $ d_0$ suitably large so that $Ce^{-\sigma d_0}\leq \f12$, then \eqref{S3.183} is reduced to be
\begin{equation}\label{S3.184}
|(\phi_0,\phi_1,\phi_2,\phi_3)(d)| \leq C\Big\{\|e^{\s_0 x}\nu^{-1}wg\|_{L^\infty_{x,v}} +  |wr|_{L^\infty(\g_-)}  \Big\},\quad \mbox{for}\quad d\geq d_0\geq 1.
\end{equation}
For $d\in[1,d_0]$, using \eqref{S3.140-1} and  \eqref{S3.140-2},  one has
\begin{equation}\nonumber
|(\phi_0,\phi_1,\phi_2,\phi_3)(d)| \leq C_{d_0} \, \{\|\nu^{-1}wg\|_{L^\infty_{x,v}} + |wr|_{L^\infty(\g_-)}\},
\end{equation}
which, together with \eqref{S3.184}, yields \eqref{S3.180}. Hence the proof of Lemma \ref{lem2.17} is complete. $\hfill\Box$

\medskip

In the following, we study the asymptotic behavior of $(\phi_0,\phi_1,\phi_2,\phi_3)(d)$ as $d\to\infty$.
Let $1\leq d_1\leq d_2<\infty$, then it holds that
\begin{align}\label{S3.155}
	\begin{cases}
		v_3 \partial_x(\bar{f}_{d_2}-\bar{f}_{d_1})+p_E^0 \FL (\bar{f}_{d_2}-\bar{f}_{d_1})=0,\quad (x,v)\in (0,d_1)\times \R^3, \\
		%\bar{f}_{d_1}(d_1,v)|_{v_3<0}=\bar{f}_{d_1}(d_1,v_{\sp},-v_3),\\
		(\bar{f}_{d_2}-\bar{f}_{d_1})(0,v)|_{v_3>0}=P_\g (\bar{f}_{d_2}-\bar{f}_{d_1})(0,v).
		%		(\bar{f}_{d_2}-\bar{f}_{d_1})(d,v)|_{v_3<0}=  (\bar{f}_{d_2}-\bar{f}_{d_1})(d,\Mr v).
	\end{cases}
\end{align}

\begin{lemma}\label{lem2.15}
	It holds that
	\begin{align}\label{S3.158}
		|({I-P_\g})(\bar{f}_{d_2}-\bar{f}_{d_1})(0)|_{L^2(\g_+)}
		\leq  Ce^{-\sigma_1 d_1 }\Big\{|r|^2_{L^2(\g_-)} + \int_0^{d_1} e^{2\sigma_1 x} \|g(x)\|^2_{L^2_{x,v}} dx \Big\}^{\f12}.
	\end{align}
\end{lemma}

\noindent{\bf Proof.} Multiplying \eqref{S3.155} by $\bar{f}_{d_2}-\bar{f}_{d_1}$, one can obtain
\begin{align}\label{S3.159}
	&\f12 \int_{\R^3} v_3 |(\bar{f}_{d_2}-\bar{f}_{d_1})(x,v)|^2 dv+c_0p_E^0\int_0^{x} \|({\bf I-P})(\bar{f}_{d_2}-\bar{f}_{d_1})(z)\|_{\nu}^2 dz\nonumber\\
	&\leq \f12 \int_{\R^3} v_3 |(\bar{f}_{d_2}-\bar{f}_{d_1})(0,v)|^2 dv, \,\, \forall x\in [0,d_1].
\end{align}
Using $\eqref{S3.155}_2$, we get that
\begin{align}\label{S3.160}
	\int_{\R^3} v_3 |(\bar{f}_{d_2}-\bar{f}_{d_1})(0,v)|^2 dv
	=-|({I-P_\g})(\bar{f}_{d_2}-\bar{f}_{d_1})(0)|^2_{L^2(\g_+)}.
\end{align}
Using  \eqref{S3.124} and \eqref{S3.125},  one has
\begin{align}\label{S3.161}
	&\int_{\R^3} v_3 |(\bar{f}_{d_2}-\bar{f}_{d_1})(x,v)|^2 dv\nonumber\\
	&=\int_{\R^3} v_3 |\bar{f}_{d_1}(x,v)|^2 dv -2\int_{\R^3} v_3 \bar{f}_{d_2}(x,v)\bar{f}_{d_1}(x,v) dv + \int_{\R^3} v_3|\bar{f}_{d_2}(x,v)|^2 dv\nonumber\\
	&=\int_{\R^3} v_3|({\bf I-P})\bar{f}_{d_1}(x,v)|^2 dv-2\int_{\R^3} v_3 ({\bf I-P})\bar{f}_{d_2}(d_1,v)\cdot ({\bf I-P})\bar{f}_{d_1}(d_1,v) dv \nonumber\\
	&\quad+\int_{\R^3} v_3|({\bf I-P})\bar{f}_{d_2}(x,v)|^2 dv.
\end{align}

Substituting \eqref{S3.160} and \eqref{S3.161} into \eqref{S3.159}, one gets
\begin{align}\label{S3.162}
	|({I-P_\g})(\bar{f}_{d_2}-\bar{f}_{d_1})(0)|^2_{L^2(\g_+)}\leq C\Big\{ \|({\bf I-P})\bar{f}_{d_1}(x)\|^2_{\nu} + \|({\bf I-P})\bar{f}_{d_2}(x)\|^2_{\nu} \Big\},\,\, \forall x\in [0,d_1].
\end{align}
Integrating \eqref{S3.162} over $x\in [d_1-1, d_1]$, and using \eqref{S3.126}, then we obtain
\begin{align}
	&|({I-P_\g})(\bar{f}_{d_2}-\bar{f}_{d_1})(0)|^2_{L^2(\g_+)}\nonumber\\
	&\leq C\Big\{ \int_{d_1-1}^{d_1} \|({\bf I-P})\bar{f}_{d_1}(x)\|^2_{\nu} dx + \int_{d_1-1}^{d_1} \|({\bf I-P})\bar{f}_{d_2}(x)\|^2_{\nu} dx\Big\}\nonumber\\
	&\leq Ce^{-2\sigma_1 d_1 }\Big\{ \int_{d_1-1}^{d_1} e^{2\sigma_1 x }\|({\bf I-P})\bar{f}_{d_1}(x)\|^2_{\nu} dx + \int_{d_1-1}^{d_1} e^{2\sigma_1 x }\|({\bf I-P})\bar{f}_{d_2}(x)\|^2_{\nu} dx\Big\}\nonumber\\
	&\leq Ce^{-2\sigma_1 d_1 }\Big\{|r|^2_{L^2(\g_-)} + \int_0^{d_1} e^{2\sigma_1 x} \|g(x)\|^2_{L^2_{x,v}} dx \Big\},\nonumber
\end{align}
which  immediately yields \eqref{S3.158}. Therefore the proof of Lemma \ref{lem2.15} is complete. $\hfill\Box$

\medskip

\begin{lemma}\label{lem2.18}
There exist constants $(\phi_0^\infty,\phi_1^\infty,\phi_2^\infty, \phi_3^\infty)$ such that
\begin{align}\label{S3.186}
\lim_{d\to +\infty}(\phi_0,\phi_1,\phi_2, \phi_3)(d)=(\phi_0^\infty,\phi_1^\infty,\phi_2^\infty, \phi_3^\infty),
\end{align}
with
\begin{equation}\label{S3.187}
|(\phi_0^\infty,\phi_1^\infty,\phi_2^\infty, \phi_3^\infty)|\leq C\Big\{\|e^{\s_0 x} \nu^{-1}wg\|_{L^\infty_{x,v}} +  |wr|_{L^\infty(\g_-)}  \Big\}.
\end{equation}
\end{lemma}

\noindent{\bf Proof.} Let $1\leq d_1\leq d_2<\infty$, it follows from \eqref{S3.140-1}, \eqref{S3.151}-\eqref{S3.151-2}, \eqref{S3.158}, $\eqref{S3.134-1}_2$ and \eqref{S3.181} that
\begin{align}\label{S3.188}
&\Big|\Big(\phi_0(d_2)-\phi_0(d_1), \phi_1(d_2)-\phi_1(d_1), \phi_2(d_2)-\phi_2(d_1),\phi_3(d_2)-\phi_3(d_1) \Big)\Big|\nonumber\\
&\leq \Big|\Big(\bar{a}_{d_2}(d_2) - \bar{a}_{d_1}(d_1),\, \bar{b}_{d_2,1}(d_2)-\bar{b}_{d_1,1}(d_1), \, \bar{b}_{d_2,2}(d_2)- \bar{b}_{d_1,2}(d_1),\, \bar{c}_{d_2}(d_2)-\bar{c}_{d_1}(d_1) \Big)\Big|\nonumber\\
&\quad +C |(\mathbf{I-P})\tilde{f}_{d_1}(d_1)|_{L^\infty_v} + C |(\mathbf{I-P})\tilde{f}_{d_2}(d_2)|_{L^\infty_v} \nonumber\\
&\leq 	|({I-P_\g})(\bar{f}_{d_2}-\bar{f}_{d_1})(0)|_{L^2(\g_+)} + \int_{d_1}^{d_2} \|g(z)\|_{L^2_v} dz \nonumber\\
&\quad+ C e^{-\sigma d_1 }\Big\{\frac{1}{\s_0-\s}\|e^{\s_0 x} \nu^{-1}wg\|_{L^\infty_{x,v}} +  |wr|_{L^\infty(\g_-)}  \Big\}\nonumber\\
&\leq C e^{-\sigma d_1 }\Big\{\frac{1}{\s_0-\s}\|e^{\s_0 x} \nu^{-1}wg\|_{L^\infty_{x,v}} +  |wr|_{L^\infty(\g_-)}  \Big\}\to 0,\quad \mbox{as}\quad d_1\to +\infty,
\end{align}
which immediately gives \eqref{S3.186}. And \eqref{S3.187} follows directly from \eqref{S3.186} and \eqref{S3.180}. Therefore the proof of Lemma \ref{lem2.18} is complete. $\hfill\Box$

\medskip

\subsection{Proof of Theorem \ref{thm3.1}}
Let $1\leq d_1\leq d_2<\infty$,  we denote
\begin{equation}
\mathfrak{f}(x,v):=(\tilde{f}_{d_2}-\tilde{f}_{d_1})(x,v),\quad\mbox{and}\quad \mathfrak{h}(x,v)=w(v)\mathfrak{f}(x,v),\quad \forall (x,v)\in[0,d_1]\times \R^3,
\end{equation}
then it follows from \eqref{S3.133} that
\begin{equation}\label{S3.189}
	\begin{cases}
		v_3 \partial_x\mathfrak{f}+p_E^0 \FL \mathfrak{f}=0,\quad (x,v)\in (0,d_1)\times \R^3, \\
		\mathfrak{f}(0,v)|_{v_3>0}=P_\g \mathfrak{f}(0,v) + (I-P_\g)(\Phi(d_2)-\Phi(d_1)).
	\end{cases}
\end{equation}
We divide the proof into three steps.

\medskip

\noindent{\it Step 1. Convergence in $L^2$-norm.} Multiplying \eqref{S3.189} by $\mathfrak{f}$ to obtain
\begin{align}\label{S3.100}
& |(I-P_\g)\mathfrak{f}(0)|^2_{L^2(\g_+)}+p_E^0c_0\int_0^{d_1}\int_{\R^3} (1+|v|) |({\bf I-P})\mathfrak{f}(x,v)|^2 dvdx\nonumber\\
&\leq C\int_{\R^3} |v_3|\cdot |\mathfrak{f}(d_1,v)|^2dv + C |(\Phi(d_2)-\Phi(d_1))|^2_{L^2(\g)}\nonumber\\
&\leq C\Big|\Big(\phi_0(d_2)-\phi_0(d_1), \phi_1(d_2)-\phi_1(d_1), \phi_2(d_2)-\phi_2(d_1),\phi_3(d_2)-\phi_3(d_1) \Big)\Big|^2\nonumber\\
&\qquad + |w\mathfrak{f}(d_1)|^2_{L^\infty_v}\nonumber\\
&\leq C e^{-2\sigma d_1 }\Big\{\frac{1}{\s_0-\s}\|e^{\s_0 x} \nu^{-1} wg\|_{L^\infty_{x,v}} +  |wr|_{L^\infty(\g_-)}\Big\}^2,
\end{align}
where we have used \eqref{S3.181} and  \eqref{S3.188} in the last inequality.

\medskip

We still need to control the macroscopic part of $\mathfrak{f}$. We denote
\begin{equation}
{\bf P}\mathfrak{f}(x,v)=\Big[\mathfrak{a}(x) +   \mathfrak{b}_1(x) v_1 +   \mathfrak{b}_2(x) v_2 + \mathfrak{c}(x) (\frac{1}{2}|v|^2-\f32) \Big]\sqrt{\mu(v)}.
\end{equation}
Noting from \eqref{S3.143}, one has
\begin{align}
\left(
	\begin{array}{c}
	\dis	\int_{\R^3} v_3 \mathfrak{f}(x,v)\cdot \FL^{-1} (\mathcal{A}_{31}) dv\\[1mm]
	\dis	\int_{\R^3} v_3 \mathfrak{f}(x,v)\cdot \FL^{-1} (\mathcal{A}_{32}) dv\\[1mm]
	\dis	\int_{\R^3} v_3 \mathfrak{f}(x,v)\cdot \FL^{-1} (\mathcal{B}_{3}) dv\\
	\end{array}
	\right)
	&=\left(
	\begin{array}{c}
	\dis	\int_{\R^3} v_3 \mathfrak{f}(d_1,v)\cdot \FL^{-1} (\mathcal{A}_{31}) dv\\[1mm]
	\dis	\int_{\R^3} v_3 \mathfrak{f}(d_1,v)\cdot \FL^{-1} (\mathcal{A}_{32}) dv\\[1mm]
	\dis	\int_{\R^3} v_3 \mathfrak{f}(d_1,v)\cdot \FL^{-1} (\mathcal{B}_{3}) dv\\
	\end{array}
	\right),\nonumber
\end{align}
which, together with \eqref{S3.138}-\eqref{S3.140}, yields
\begin{align}
|(\mathfrak{b}_1(x), \mathfrak{b}_2(x), \mathfrak{c}(x))|
&\leq C |\mathfrak{f}(d_1)|_{L^\infty_v} + \Big(\int_{\R^3}|({\bf I-P})\mathfrak{f}(x,v)|^2dv\Big)^{\f12}.\nonumber
\end{align}
Hence it follows from  \eqref{S3.181} and \eqref{S3.100} that
\begin{align}\label{S3.101}
\int_0^{d_1}|(\mathfrak{b}_1 , \mathfrak{b}_2 , \mathfrak{c})(x)|^2dx  &\leq Cd_1|\mathfrak{f}(d_1)|^2_{L^\infty_v}+\int_0^{d_1}\int_{\R^3}|({\bf I-P})\mathfrak{f}(x,v)|^2dvdx\nonumber\\
&\leq C d_1 \,e^{-2\sigma d_1 }\Big\{\frac{1}{\s_0-\s}\|e^{\s_0 x}\nu^{-1}wg\|_{L^\infty_{x,v}} +  |wr|_{L^\infty(\g_-)}\Big\}^2.
\end{align}
Noting \eqref{S3.148} and using \eqref{S3.101}, \eqref{S3.100},  one can obtain
\begin{align}
\int_0^{d_1}|\mathfrak{a}(x)|^2dx &\leq \int_0^{d_1}|\mathfrak{c}(x)|^2dx +\int_0^{d_1}\int_{\R^3}|({\bf I-P})\mathfrak{f}(x,v)|^2dvdx\nonumber\\
&\leq C d_1 \,e^{-2\sigma d_1 }\Big\{\frac{1}{\s_0-\s}\|e^{\s_0 x}\nu^{-1}wg\|_{L^\infty_{x,v}} +  |wr|_{L^\infty(\g_-)}\Big\}^2,\nonumber
\end{align}
which, together with \eqref{S3.100} and \eqref{S3.101}, gives
\begin{align}\label{S3.105}
\int_0^{d_1}\int_{\R^3}|\mathfrak{f}(x,v)|^2dvdx\leq C d_1 \,e^{-2\sigma d_1 }\Big\{\frac{1}{\s_0-\s}\|e^{\s_0 x} \nu^{-1}wg\|_{L^\infty_{x,v}} +  |wr|_{L^\infty(\g_-)}\Big\}^2.
\end{align}

\medskip

\noindent{\it Step 2. Convergence in $L^\infty_{x,v}$-norm.} In the following, we use $t_{k}=t_{k}(t,x,v), X_{cl}(s;t,x,v), x_k=x_k(v,x)$ to denote the back-time cycles defined for domain $[0,d_1]\times \R^3$.  Let $(x,v)\in [0,d_1]\times\R^3\backslash (\g_0\cup\g_{-})$, it follows from \eqref{S3.189} that
\begin{align}\label{S3.190}
\mathfrak{h}(x,v)
%&=\Fi_{\{v_{0,3}<0\}} \bigg\{ e^{-\nu_{E}(v) (t-t_1)}\cdot \mathfrak{h}(d_1,v) \nonumber\\
%&\quad + \int_{t_1}^t e^{-\nu_{E}(v) (t-s)}\big(K_w \mathfrak{h}\big)(x-v_{0,3}(t-s),v) ds \bigg\}\nonumber\\
%&\quad + \Fi_{\{v_{0,3}>0\}} \bigg\{e^{-\nu_{E}(v) (t-t_1)}\cdot  (I-P_\g)(\Phi(d_2)-\Phi(d_1))(0,v)\nonumber\\
%&\quad+\int_{t_1}^t e^{-\nu_{E}(v) (t-s)}\big(K_w \mathfrak{h}\big)(x-v_{0,3}(t-s),v) ds \nonumber\\
%&\quad + \frac{e^{-\nu_{E}(v) (t-t_1)}}{\tilde{w}(v)} \int_{\mathcal{V}_1}\int_{t_2}^{t_1} e^{-\nu_{E}(v) (t_1-s)}\big(K_w \mathfrak{h}\big)(-v_{1,3}(t_1-s),v_1) \tilde{w}(v_1) ds d\sigma_1 \nonumber\\
%&\quad+ \frac{e^{-\nu_{E}(v) (t-t_1)}}{\tilde{w}(v)} \int_{\mathcal{V}_1} e^{-\nu_{E}(v_1) (t_1-t_2)} \mathfrak{h}(d_1,v_1) \tilde{w}(v_1) d\sigma_1 \bigg\}\nonumber\\
&=\Fi_{\{v_{0,3}>0\}} \frac{e^{-\nu_{E}(v) (t-t_1)}}{\tilde{w}(v)} \int_{\mathcal{V}_1}\int_{t_2}^{t_1} e^{-\nu_{E}(v_1) (t_1-s)}\big(K_w \mathfrak{h}\big)(-v_{1,3}(t_1-s),v_1) \tilde{w}(v_1) ds d\sigma_1\nonumber\\
&\quad+ \int_{t_1}^t e^{-\nu_{E}(v) (t-s)}\big(K_w \mathfrak{h}\big)(x-v_{0,3}(t-s),v) ds +\mathfrak{B},
\end{align}
where we have used the notations $\nu_E(v):=p_E^0\, \nu(v)$ and
\begin{align}\label{S3.191}
\mathfrak{B}:&= \Fi_{\{v_{0,3}>0\}} \bigg\{e^{-\nu_{E}(v) (t-t_1)} w(v)\cdot  (I-P_\g)(\Phi(d_2)-\Phi(d_1))(v)\nonumber\\
&\quad+ \frac{e^{-\nu_{E}(v) (t-t_1)}}{\tilde{w}(v)} \int_{\mathcal{V}_1} e^{-\nu_{E}(v_1) (t_1-t_2)} \mathfrak{h}(d_1,v_1) \tilde{w}(v_1) d\sigma_1 \bigg\}\nonumber\\
&\quad+ \Fi_{\{v_{0,3}\leq0\}}\cdot  e^{-\nu_{E}(v) (t-t_1)}\cdot \mathfrak{h}(d_1,v).
\end{align}

Using \eqref{S3.188}, \eqref{S3.181}, it holds that
\begin{align}\label{S3.192}
|\mathfrak{B}|&\leq C|(\Phi(d_2)-\Phi(d_1))| + |\mathfrak{h}(d_1)|_{L^\infty} \nonumber\\
&\leq C e^{-\sigma d_1 }\Big\{\frac{1}{\s_0-\s}\|e^{\s_0 x}\nu^{-1}wg\|_{L^\infty_{x,v}} +  |wr|_{L^\infty(\g_-)}  \Big\}.
\end{align}
By similar arguments as in \eqref{S3.34}-\eqref{S3.41}, one can obtain
\begin{align}\label{S3.194}
&\Fi_{\{v_{0,3}>0\}} \frac{e^{-\nu_{E}(v) (t-t_1)}}{\tilde{w}(v)} \left|\int_{\mathcal{V}_1}\int_{t_2}^{t_1} e^{-\nu_{E}(v) (t_1-s)}\big(K_w \mathfrak{h}\big)(-v_{1,3}(t_1-s),v_1) \tilde{w}(v_1) ds d\sigma_1\right|\nonumber\\
&\leq \frac{C}{N} \|\mathfrak{h}\|_{L^\infty([0,d_1]\times\R^3)} +C_N \|\mathfrak{f}\|_{L^2([0,d_1]\times\R^3)}.
\end{align}
Combining \eqref{S3.190},\eqref{S3.192} and \eqref{S3.194}, one gets, for $(x,v)\in [0,d_1]\times\R^3$
\begin{align}\label{S3.195}
|\mathfrak{h}(x,v)|&\leq \int_{t_1}^t e^{-\nu_{E}^0(t-s)} ds\int_{\R^3} |k_w(v,v') \mathfrak{h}(x',v')| dv' + \frac{C}{N} \|\mathfrak{h}\|_{L^\infty([0,d_1]\times\R^3)} \nonumber\\
&\quad +C_N \|\mathfrak{f}\|_{L^2([0,d_1]\times\R^3)}+C e^{-\sigma d_1 }\Big\{\frac{1}{\s_0-\s}\|e^{\s_0 x}\nu^{-1}wg\|_{L^\infty_{x,v}} +  |wr|_{L^\infty(\g_-)}  \Big\}
\end{align}
where we have denoted $x':=x-v_{0,3}(t-s)\in [0,d_1]$.

\medskip

For the first term on right hand side of \eqref{S3.195}, we use \eqref{S3.195} again to obtain
\begin{align}\label{S3.196}
&\int_{t_1}^t e^{-\nu_{E}^0(t-s)} ds\int_{\R^3} |k_w(v,v') \mathfrak{h}(x',v)| dv'\nonumber\\
&\leq \int_{t_1}^t e^{-\nu_{E}^0(t-s)} ds \int^s_{t_1'} e^{-\nu_{E}^0(s-\tau)} d\tau \int_{\R^3_{v'}}\int_{\R^3_{v''}} |k_w(v,v') k_w(v',v'') \mathfrak{h}(x'',v'')|dv'' dv'\nonumber\\
&\quad +\frac{C}{N} \|\mathfrak{h}\|_{L^\infty([0,d_1]\times\R^3)} +C_N \|\mathfrak{f}\|_{L^2([0,d_1]\times\R^3)}\nonumber\\
&\quad +C e^{-\sigma d_1 }\Big\{\frac{1}{\s_0-\s}\|e^{\s_0 x}\nu^{-1}wg\|_{L^\infty_{x,v}} +  |wr|_{L^\infty(\g_-)}  \Big\}
\end{align}
where $t_1'=t_1(s,x',v')$ and $x''=x'-v_{0,3}'(s-\tau)$. Applying similar arguments as \eqref{S3.48}-\eqref{S3.53} to the first term on right hand side of \eqref{S3.196}, then we can have
 \begin{align}
 	\int_{t_1}^t e^{-\nu_{E}^0(t-s)} ds\int_{\R^3} |k_w(v,v') \mathfrak{h}(x',v)| dv'
 	&\leq  \frac{C}{N} \|\mathfrak{h}\|_{L^\infty([0,d_1]\times\R^3)} +C_N \|\mathfrak{f}\|_{L^2([0,d_1]\times\R^3)}\nonumber\\
 	&+C e^{-\sigma d_1 }\Big\{\frac{1}{\s_0-\s}\|e^{\s_0 x}\nu^{-1}wg\|_{L^\infty_{x,v}} +  |wr|_{L^\infty(\g_-)}  \Big\},\nonumber
 \end{align}
which, together with \eqref{S3.195}, yields that, for $(x,v)\in ([0,d_1]\times \R^3)\backslash (\g_0\cup \g_-)$
\begin{align}
|\mathfrak{h}(x,v)|  &\leq  \frac{C}{N} \|\mathfrak{h}\|_{L^\infty([0,d_1]\times\R^3)} +C_N \|\mathfrak{f}\|_{L^2([0,d_1]\times\R^3)}\nonumber\\
&+C e^{-\sigma d_1 }\Big\{\frac{1}{\s_0-\s}\|e^{\s_0 x}\nu^{-1}wg\|_{L^\infty_{x,v}} +  |wr|_{L^\infty(\g_-)}  \Big\}.\nonumber
\end{align}
Hence it holds
\begin{align}
\|\mathfrak{h}\|_{L^\infty([0,d_1]\times\R^3)} + |\mathfrak{h}(0)|_{L^\infty(\g)}  &\leq  \frac{C}{N} \|\mathfrak{h}\|_{L^\infty([0,d_1]\times\R^3)} +C_N \|\mathfrak{f}\|_{L^2([0,d_1]\times\R^3)}\nonumber\\
&+C e^{-\sigma d_1 }\Big\{\frac{1}{\s_0-\s}\|e^{\s_0 x}\nu^{-1}wg\|_{L^\infty_{x,v}} +  |wr|_{L^\infty(\g_-)}  \Big\}.\nonumber
\end{align}
Taking $N>0$ suitably large so that $\frac{C}{N}\leq \f12$ and using \eqref{S3.181}, \eqref{S3.105}, then one obtains
\begin{align}\label{S3.197}
&\|\mathfrak{h}\|_{L^\infty([0, d_1]\times\R^3)} + |\mathfrak{h}(0)|_{L^\infty(\g)}  \nonumber\\
&\leq C  \|\mathfrak{f}\|_{L^2([0,d_1]\times\R^3)} +C e^{-\sigma d_1 }\Big\{\frac{1}{\s_0-\s}\|e^{\s_0 x}\nu^{-1}wg\|_{L^\infty_{x,v}} +  |wr|_{L^\infty(\g_-)}  \Big\}\nonumber\\
&\leq C d_1 e^{-\f12\sigma d_1 }\Big\{\frac{1}{\s_0-\s}\|e^{\s_0 x}\nu^{-1}wg\|_{L^\infty_{x,v}} +  |wr|_{L^\infty(\g_-)}  \Big\}\,\,\mbox{for}\,\,d_1\to \infty.
\end{align}
%Combining \eqref{S3.197} and \eqref{S3.181}, one gets
%\begin{align}\label{S3.206}
%&\|\mathfrak{h}\|_{L^\infty([0, d_1]\times\R^3)} + |\mathfrak{h}(0)|_{L^\infty(\g)}  \nonumber\\
%&\leq  C d_1 e^{-\f12\sigma d_1 }\Big\{\frac{1}{\s_0-\s}\|e^{\s_0 x}\nu^{-1}wg\|_{L^\infty_{x,v}} +  |wr|_{L^\infty(\g_-)}  \Big\}\to 0, \,\,\mbox{for}\,\,d_1\to \infty.
%\end{align}
With the help of \eqref{S3.197},
there exists a function $f(x,v)$ with $(x,v)\in\R_+\times \R^3$ so that \begin{equation*}
\|w(\tilde{f}_{d}-f)\|_{L^\infty([0, d]\times\R^3)} + |w(\tilde{f}_{d}-f)(0)|_{L^\infty(\g_+)}\rightarrow0,\,\,\,\mbox{as}\,\, d \rightarrow \infty,
\end{equation*}
It follows  from \eqref{S3.181} and the strong convergence in $L^\infty_{x,v}$ that
\begin{equation}\label{S3.207}
\|e^{\s x}wf\|_{L^\infty_{x,v}}+|e^{\s x}wf(0)|_{L^\infty(\g)}
\leq C\Big\{\frac{1}{\s_0-\s}\|e^{\s_0 x}\nu^{-1}wg\|_{L^\infty_{x,v}} +  |wr|_{L^\infty(\g_-)}  \Big\}.
\end{equation}
 The continuity of $f$  follows directly from the $L^\infty_{x,v}$-convergence and the continuity of $\tilde{f}_{d}$. It is  direct to see that $f(x,v)$ solves
\begin{align}\label{S3.208}
\begin{cases}
v_3 \partial_xf+p_E^0 \FL f=g,\quad (x,v)\in (0,\infty)\times \R^3, \\[2mm]
f(0,v)|_{v_3>0}=P_\g f(0,v) + (I-P_\g)(\Phi^\infty)(v) + r(v),\\[2mm]
\lim_{x\rightarrow\infty} f(x,v)=0,
\end{cases}
\end{align}
with
\begin{align}
\Phi^\infty(v)=\Big\{\phi_0^\infty+  \phi_1^\infty\, v_1+   \phi_2^\infty \,v_2 + \phi_3^\infty \,  (\f{|v|^2}{2}-\f32)\Big\}\sqrt{\mu(v)},\label{S3.209}\\
|(\phi_0^\infty,\phi_1^\infty,\phi_2^\infty, \phi_3^\infty)|\leq C\Big\{\|e^{\s_0 x}\nu^{-1}wg\|_{L^\infty_{x,v}} +  |wr|_{L^\infty(\g_-)}  \Big\}.\label{S3.210}
\end{align}
where the constants $(\phi_0^\infty,\phi_1^\infty,\phi_2^\infty, \phi_3^\infty)$ are the ones determined in Lemma \ref{lem2.18}.

\smallskip

We note that
\begin{equation}\nonumber
(I-P_\g)(\Phi^\infty)\equiv (I-P_\g)(\alpha \sqrt{\mu}+ \Phi^\infty),\quad \forall\,\, \alpha\in \R,
\end{equation}
that means the system \eqref{S3.208} is invariant under such transformation. So we always normalized $\Phi^{\infty}$ such that
\begin{equation}\label{S3.212-3}
\phi_0^\infty=0.
\end{equation}
%In the following we shall prove that the solution $(\hat{f}, \Phi^\infty)$ is unique under constraints \eqref{S3.207}, \eqref{S3.209}-\eqref{S3.212-3}.

\

\noindent{\it Step 3. Uniqueness.} For any given $g\in \mathcal{N}^{\perp}$ and $r$ with
\begin{align}
\|e^{\s_0 x}\nu^{-1}wg\|_{L^\infty_{x,v}} +  |wr|_{L^\infty(\g_-)}<\infty,\quad \mbox{and} \quad \int_{v_3>0} v_3 \sqrt{\mu(v)} r(v) dv=0.
\end{align}
Let $(\hat{f}_1, \Phi_1^\infty)$ and $(\hat{f}_2, \Phi_2^\infty)$ be two solutions of \eqref{S3.208}, respectively, with \eqref{S3.207}, \eqref{S3.209}-\eqref{S3.212-3}.
Noting \eqref{S3.209} and \eqref{S3.212-3}, we denote
\begin{align}\label{S3.212-1}
\begin{split}
\Phi_1^\infty&:=\Big\{ \phi_{1,1}^\infty\, v_1+   \phi_{1,2}^\infty \,v_2 + \phi_{1,3}^\infty \,  (\f{|v|^2}{2}-\f32)\Big\}\sqrt{\mu(v)},\\
\Phi_2^\infty&:=\Big\{\phi_{2,1}^\infty\, v_1+   \phi_{2,2}^\infty \,v_2 + \phi_{2,3}^\infty \,  (\f{|v|^2}{2}-\f32)\Big\}\sqrt{\mu(v)}.
\end{split}
\end{align}
Moreover, for the solutions $\hat{f}_i, i=1,2$, it is direct to check, for $x\in[0,+\infty)$, that
\begin{align}\label{S3.210-1}
	\begin{split}
		&\int_{\R^3} v_3 \sqrt{\mu} \hat{f}_i(x,v)dv=0,\\
		&\int_{\R^3} \mathcal{A}_{3j}(v)\hat{f}_i(x,v) dv=\int_{\R^3} \mathcal{A}_{3j}(v) ({\bf I-P})\hat{f}_i(x,v) dv=0,\\
		&\int_{\R^3} \mathcal{B}_3(v) \hat{f}_i(x,v) dv=\int_{\R^3} \mathcal{B}_3(v) ({\bf I-P})\hat{f}_i(x,v) dv=0,
	\end{split}
\,\, \forall \,\, i,j=1,2.
\end{align}

\medskip

We define
\begin{align}\label{S3.211-1}
\hat{\mathfrak{f}}(x,v):=[\hat{f}_1(x,v)-\Phi_1^\infty]-[\hat{f}_2(x,v)-\Phi_2^\infty],
\end{align}
then it follows from \eqref{S3.208} that
\begin{align}\label{S3.211}
	\begin{cases}
		v_3 \partial_x\hat{\mathfrak{f}}+p_E^0 \FL\hat{\mathfrak{f}}=0,\quad (x,v)\in (0,\infty)\times \R^3, \\[2mm]
		\hat{\mathfrak{f}}(0,v)|_{v_3>0}=P_\g \hat{\mathfrak{f}}(0,v),\\[2mm]
		\lim_{x\rightarrow\infty} e^{\s x}\|w(v)\big[\hat{\mathfrak{f}}-(-\Phi_1^\infty+\Phi_2^\infty)\big]\|_{L^\infty_v}=0.
	\end{cases}
\end{align}
Multiplying $\eqref{S3.211}_1$ by $\hat{\mathfrak{f}}$, and integrating  over $(x,v)\in\R_+\times\R^3$, one can obtain
\begin{align}\label{S3.215}
\f12 \int_{\R^3} v_3 |(\Phi_1^\infty-\Phi_2^\infty)(v)|^2 dv
+c_0p_E^0\int_0^\infty \|({\bf I-P}) \hat{\mathfrak{f}}(x)\|_{\nu}^2
\leq \f12 \int_{\R^3} v_3 |\hat{\mathfrak{f}}(0,v)|^2 dv,
\end{align}
where we have used  $\eqref{S3.211}_3$.

Noting  \eqref{S3.212-1}, one can check that
\begin{equation}\label{S3.216}
\int_{\R^3} v_3 |(\Phi_1^\infty-\Phi_2^\infty)(v)|^2 dv=0.
\end{equation}
It follows from  $\eqref{S3.211}_2$ that
\begin{equation}\label{S3.217}
\int_{\R^3} v_3 |\hat{\mathfrak{f}}(0,v)|^2 dv= - |(I-P_\g)\hat{\mathfrak{f}}(0)|_{L^2(\g_+)}^2.
\end{equation}
Substituting \eqref{S3.216} and \eqref{S3.217} into \eqref{S3.215}, one has
\begin{equation}\label{S3.218}
|(I-P_\g)\hat{\mathfrak{f}}(0)|_{L^2(\g_+)}^2+c_0p_E^0\int_0^\infty \|({\bf I-P}) \hat{\mathfrak{f}}(x)\|_{\nu}^2\leq 0,
\end{equation}
which immediately yields that
\begin{align}\label{S3.219}
0\equiv ({\bf I-P}) \hat{\mathfrak{f}}(x,v)
%= ({\bf I-P}) \hat{f}_1(x,v) - ({\bf I-P}) \hat{f}_1(x,v)
,\,\,\, \forall \,\, (x,v)\in [0,\infty)\times \R^3,
\end{align}
and
\begin{align}\label{S3.220}
\hat{\mathfrak{f}}(0,v)=P_\g \hat{\mathfrak{f}}(0,v), \,\,\, \forall\,\, v\in \R^3.
\end{align}

\medskip

We still need to prove ${\bf P}\hat{\mathfrak{f}}(x,v)\equiv0$. We denote
\begin{equation}\nonumber
	{\bf P}\hat{\mathfrak{f}}(x,v)=\Big[\hat{\mathfrak{a}}(x) +   \hat{\mathfrak{b}}_1(x) v_1 +   \hat{\mathfrak{b}}_2(x) v_2 + \hat{\mathfrak{c}}(x) (\frac{1}{2}|v|^2-\f32) \Big]\sqrt{\mu(v)}.
\end{equation}
Using \eqref{S3.220}, one can check that
\begin{align}\label{S3.221}
	\begin{split}
		&\int_{\R^3} \hat{\mathfrak{f}}(0,v) \cdot v_3^2 \sqrt{\mu(v)}   dv=\sqrt{2\pi}z_{\g_+}(\hat{f}),\\
		&\int_{\R^3} \hat{\mathfrak{f}}(0,v)\cdot v_3 \FL^{-1} (\mathcal{A}_{3i}) dv=0,\,\, i=1,2,\\
		&\int_{\R^3}  \hat{\mathfrak{f}}(0,v)\cdot v_3\FL^{-1} (\mathcal{B}_{3}) dv=0,
	\end{split}
\end{align}
where $\dis z_{\g_+}(\hat{f})=\int_{v_3<0} |v_3|\sqrt{\mu(v)} \, \hat{\mathfrak{f}}(0,v) dv.$

Noting \eqref{S3.219}, multiplying $\eqref{S3.211}$ by $v_3\sqrt{\mu}, \FL^{-1} (\mathcal{A}_{31}), \FL^{-1} (\mathcal{A}_{32})$ and $\FL^{-1} (\mathcal{B}_{3})$, respectively, and integrating over $[0,x]$ and using \eqref{S3.221}, one gets
\begin{align}\label{S3.223}
\begin{split}
&\dis\int_{\R^3} \hat{\mathfrak{f}}(x,v) \cdot v_3^2 \sqrt{\mu(v)}   dv=\int_{\R^3} \hat{\mathfrak{f}}(0,v) \cdot v_3^2 \sqrt{\mu(v)}   dv=\sqrt{2\pi}z_{\g_+}(\hat{f}),\\
&\dis\int_{\R^3} \hat{\mathfrak{f}}(x,v)\cdot v_3 \FL^{-1} (\mathcal{A}_{3i}) dv=\int_{\R^3} \hat{\mathfrak{f}}(0,v)\cdot v_3 \FL^{-1} (\mathcal{A}_{3i}) dv=0,\,\,\, i=1,2,\\
&\dis\int_{\R^3}  \hat{\mathfrak{f}}(x,v)\cdot v_3\FL^{-1} (\mathcal{B}_{3}) dv=\int_{\R^3}  \hat{\mathfrak{f}}(0,v)\cdot v_3\FL^{-1} (\mathcal{B}_{3}) dv=0.
\end{split}
\end{align}
Noting \eqref{S3.219}, by similar arguments as in \eqref{S3.137}-\eqref{S3.140}, then one has
\begin{align}\nonumber
\begin{split}
\dis\sqrt{2\pi}z_{\g_+}(\hat{f})&=\int_{\R^3} \hat{\mathfrak{f}}(x,v) \cdot v_3^2 \sqrt{\mu(v)}   dv=\hat{\mathfrak{a}}(x)+\hat{\mathfrak{c}}(x),\\
\dis0&=\int_{\R^3} \hat{\mathfrak{f}}(x,v)\cdot v_3 \FL^{-1} (\mathcal{A}_{3i}) dv=\kappa_1 \,\hat{\mathfrak{b}}_i(x),\,\,\, i=1,2,\\
\dis0&=\int_{\R^3}  \hat{\mathfrak{f}}(x,v)\cdot v_3\FL^{-1} (\mathcal{B}_{3}) dv=\kappa_2 \, \hat{\mathfrak{c}}(x),
\end{split}
\end{align}
that is
\begin{equation}\label{S3.224}
\hat{\mathfrak{a}}(x)\equiv \sqrt{2\pi} z_{\g_+}(\hat{f}),\qquad
\hat{\mathfrak{b}}_i(x)\equiv 0,\,\, i=1,2,\qquad \hat{\mathfrak{c}}(x)\equiv 0,\,\,\forall\,\,\, x\in[0,+\infty).
\end{equation}
It follows from \eqref{S3.224} and \eqref{S3.219} that
\begin{equation*}
\hat{\mathfrak{f}}(x,v)\equiv\sqrt{2\pi}z_{\g_+}(\hat{f}) \sqrt{\mu(v)},
\end{equation*}
which, together with $\eqref{S3.211}_3$, yields
\begin{equation*}
\hat{z}_{\g_+}(0)=0\quad \mbox{and}\quad \Phi_1^\infty=\Phi_2^\infty.
\end{equation*}
Therefore we have proved
\begin{equation*}
\hat{f}_1(x,v)\equiv \hat{f}_2(x,v)\quad \mbox{and}\quad \Phi_1^\infty=\Phi_2^\infty.
\end{equation*}
That is the solution $f,\, \Phi^\infty$ of \eqref{S3.208} is unique under the constraints  \eqref{S3.207}, \eqref{S3.209}-\eqref{S3.212-3}.

\smallskip

Finally we define
\begin{align*}
\mathcal{G}(g,r) \equiv	f^\infty(v):=\Big\{ b_{1}^\infty\, v_1+   b_{2}^\infty \,v_2 + c^\infty \,  (\f{|v|^2}{2}-\f32)\Big\}\sqrt{\mu(v)},
\end{align*}
with
\begin{equation}\nonumber
	b_{1}^\infty=-\phi_1^\infty,\qquad b_{2}^\infty=-\phi_2^\infty,\qquad
	c^\infty=-\phi_3^\infty.
\end{equation}
Then it is direct to know that $f, f^\infty$ constructed above satisfies \eqref{1.7-3}. Therefore the proof of Theorem \ref{thm3.1} is complete. $\hfill\Box$

%%%%%%%%%%%%%%%%%%%%%%%%%%%%%%%%%%%%%%%%%%%%%%%%%%%%%%%%%%%%%%%%%%%%%%%%%%%%

\section{Proof of Theorem \ref{thm1.1}} \label{M}
We consider the following iterative sequence
\begin{align}\label{5.1}
	\begin{cases}
		v_3\cdot \partial_x f_{j+1}+p_E^0 \FL f_{j+1}=\Gamma(f_j,f_j) +\mathfrak{S},\\[1.5mm]
		f_{j+1}(0,v)|_{v_3>0}=P_\g f_{j+1}(0,v)-(I-P_\g)f^\infty_{j+1}+\mathfrak{R},\\[1.5mm]
		\lim_{x\rightarrow\infty} f^{j+1}(x,v)=0,
	\end{cases}
\end{align}
for $j=0,1,2\cdots$ with $f_0\equiv0$.  It follows from \eqref{2.16} that  $\Gamma(f_j,f_j)\in \mathcal{N}^{\perp}$
and
\begin{align}\label{5.2}
	\|\nu^{-1} w \Gamma(f_j,f_j)\|_{L^\infty_v}\leq C \|wf_j\|^2_{L^\infty_v}.
\end{align}
Taking $\frac12\s_0\leq\s<\s_0$, using \eqref{5.2} and Theorem  \ref{thm3.1}, we can find solution $(f_{j+1},f_{j+1}^\infty)$  of  \eqref{5.1} inductively for $j=0,1,2,\cdots$, where
\begin{align}
f_{j+1}^\infty&=\mathcal{G}\big(\Gamma(f_j,f_j) +\mathfrak{S},\, \mathfrak{R} \big)
:=\Big\{ b_{j+1,1}^\infty\, v_1+   b_{j+1,2}^\infty \,v_2 + c^\infty_{j+1} \,  (\f{|v|^2}{2}-\f32)\Big\}\sqrt{\mu(v)},\nonumber
\end{align}
It follows from \eqref{5.2}, \eqref{2.21} and \eqref{1.8-0} that
\begin{align}\label{5.3-1}
&|(b^\infty_{j+1,1},\, b^\infty_{j+1,2},\, c^\infty_{j+1})| \leq \hat{C}_1 \Big(|w\mathfrak{R}|_{L^\infty(\gamma_-)} +\|e^{\s_0 x}\nu^{-1}w\mathfrak{S}\|_{L^\infty_{x,v}} + \|e^{\f12\s_0 x}wf_j\|^2_{L^\infty_{x,v}}\Big)
\end{align}
and
\begin{align}\label{5.3}
&\|e^{\s x}wf_{j+1}\|_{L^\infty_{x,v}}+|wf_{j+1}(0)|_{L^\infty(\gamma)}
\nonumber\\
&\leq \hat{C}_1 |w\mathfrak{R}|_{L^\infty(\gamma_-)}+ \frac{\hat{C}_1}{\s_0-\s} \Big( \|e^{\s_0 x}\nu^{-1}w\mathfrak{S}\|_{L^\infty_{x,v}} + \|e^{\f12\s_0 x}wf_j\|^2_{L^\infty_{x,v}} \Big).
\end{align}

\smallskip

We denote
\begin{equation*}
\delta:=|w\mathfrak{R}|_{L^\infty(\gamma_-)} + \|e^{\s_0 x}\nu^{-1}w\mathfrak{S}\|_{L^\infty_{x,v}}
\end{equation*}
By induction, we shall prove  that
\begin{align}\label{5.4}
\begin{split}
|(b^\infty_{j,1},\, b^\infty_{j,2},\, c^\infty_{j})|&\leq 2\hat{C}_1\d,\\
\|e^{\s x} wf_{j}\|_{L^\infty_{x,v}}+|wf_{j}(0)|_{L^\infty(\gamma)}
&\leq \frac{2\hat{C}_1\delta}{\s_0-\s},
\end{split}
\quad\mbox{for} \  j=1,2,\cdots.
\end{align}
Indeed, for $j=0$, it follows from $f_0\equiv0$ and \eqref{5.3-1}-\eqref{5.3} that
\begin{equation*}
\begin{split}
|(b^\infty_{1,1},\, b^\infty_{1,2},\, c^\infty_{1})|&\leq \hat{C}_1 \delta, \\
\|e^{\s x}wf_{1}\|_{L^\infty_{x,v}}+|wf_{1}(0)|_{L^\infty(\gamma)}
&\leq \frac{\hat{C}_1\delta}{\s_0-\s}.
\end{split}
\end{equation*}
Now we assume that \eqref{5.4} holds for $j=1,2\cdots, l$,  then we consider the case for $j=l+1$. Indeed it follows from \eqref{5.3-1}-\eqref{5.3}  that
\begin{align*}
|(b^\infty_{l+1,1},\, b^\infty_{l+1,2},\, c^\infty_{l+1})|
&\leq \hat{C}_1\delta + \hat{C}_1 \|e^{\f12\s_0 x}wf_j\|^2_{L^\infty_{x,v}}\nonumber\\
&\leq \hat{C}_1\delta \Big(1+\d (\f{4\hat{C}_1}{\s_0})^2\Big)\leq \f32 \hat{C}_1 \delta,
\end{align*}
and
\begin{align*}%\label{S3.130}
\|e^{\s x} wf_{l+1}\|_{L^\infty_{x,v}}+|wf_{l+1}(0)|_{L^\infty(\gamma)}
&\leq \frac{\hat{C}_1\delta}{\s_0-\s}+\frac{\hat{C}_1}{\s_0-\s} \|e^{\f12\s_0 x}wf_l\|^2_{L^\infty_{x,v}}\nonumber\\
&\leq \frac{\hat{C}_1\delta}{\s_0-\s}\Big(1+\delta(\frac{4\hat{C}_1}{\s_0})^2\Big)\leq  \frac{\f{3}{2}\hat{C}_1\delta}{\s_0-\s},
\end{align*}
where we have used  $\eqref{5.4}$ with $j=l$, and chosen $\delta\leq \d_0$ with $\delta_0$ small enough such that  $(\frac{4\hat{C}_1}{\s_0})^2 \delta_0\leq 1/2$.
Therefore we have proved \eqref{5.4} by induction.

Finally we consider the convergence of sequence  $f_j$. For the difference $f_{j+1}-f_j$,  we apply \eqref{1.22} to have
\begin{align}\label{5.5}
&\|e^{\f12\s_0 x}w\{f_{j+1}-f_j\}\|_{L^\infty_{x,v}}+|w\{f_{j+1}-f_j\}(0)|_{L^\infty(\gamma)}\nonumber\\
&\quad + |(b^\infty_{j+1,1}-b^\infty_{j,1},\, b^\infty_{j+1,2}-b^\infty_{j,2},\, c^\infty_{j+1}-c^\infty_{j})| \nonumber\\
&\leq \frac{4\hat{C}_1}{\s_0}\Big\{ \|e^{\s_0 x}\nu^{-1}w\Gamma(f_{j}-f_{j-1},f_{j})\|_{L^\infty_{x,v}} +\|e^{\s_0 x}\nu^{-1}w\Gamma(f_{j-1},f_{j}-f_{j-1})\|_{L^\infty_{x,v}}\Big\}\nonumber\\
&\leq \frac{4\hat{C}_1}{\s_0}[\|e^{\f12\s_0 x}wf_{j}\|_{L^\infty_{x,v}}+\|e^{\f12\s_0 x}wf_{j-1}\|_{L^\infty_{x,v}}]\cdot \|e^{\f12\s_0 x}w(f_{j}-f_{j-1})\|_{L^\infty_{x,v}}\nonumber\\
&\leq \delta 2(\frac{4\hat{C}_1}{\s_0})^2  \|e^{\s x}w(f_{j}-f_{j-1})\|_{L^\infty_{x,v}} %\nonumber\\
\leq \frac12 \|e^{\s x}w(f_{j}-f_{j-1})\|_{L^\infty_{x,v}},
\end{align}
where we have used  that fact $2(\frac{4\hat{C}_1}{\s_0})^2 \delta_0\leq 1/2$. Hence both $f_j$ and $ f^\infty_j$ are  Cauchy sequences, then we obtain the solution by taking the limits
\begin{equation}\label{5.8}
\mathfrak{f}=\lim_{j\rightarrow\infty} f_j\qquad\mbox{and}\qquad \mathbb{G}(\mathfrak{S},\mathfrak{R})\equiv\mathfrak{f}^\infty:=\lim_{j\to\infty} \mathfrak{f}^\infty_j.
\end{equation}
The estimates $\eqref{2.6-3}_2$ and \eqref{1.14} follow from \eqref{5.8} and \eqref{5.4}. The uniqueness can also be obtained by using the inequality as \eqref{5.5}. The continuity of $f$ is   a direct consequence of $L^\infty_{x,v}$-convergence. It follows from  the uniqueness that $\mathbb{G}(0,0)=0$.  %The positivity of $F_*:=\mu+\sqrt{\mu} f_*$ will be proved in {Section 4}.

\smallskip

Finally we prove that the solution $\mathfrak{f}, \mathfrak{f}^\infty$ depend continuously on $\mathfrak{S}, \mathfrak{R}$. Let $\mathfrak{S}_i\in\mathcal{N}^{\perp}$ and $\int_{v_3>0} v_3 \sqrt{\mu(v)} \mathfrak{R}_i(v) dv=0$, $i=1,2$ satisfying \eqref{2.5}. Let $\mathfrak{f}_i, \mathfrak{f}^\infty_i$ be the solutions obtain obtained above  by replacing $\mathfrak{S},\, \mathfrak{R}$ by $\mathfrak{S}_i,\, \mathfrak{R}_i$, and  denote
\begin{equation*}
	\mathfrak{f}_i^\infty(v)= \mathbb{G}(\mathfrak{S}_i,\mathfrak{R}_i)=\Big\{ \mathfrak{b}_{i,1}^\infty\, v_1+   \mathfrak{b}_{i,2}^\infty \,v_2 + \mathfrak{c}_i^\infty \,  (\f{|v|^2}{2}-\f32)\Big\}\sqrt{\mu(v)}\,\,\, i=1,2.
\end{equation*}
It is direct to check that $\mathfrak{f}_1-\mathfrak{f}_2$ satisfies
\begin{equation}\label{5.9}
\begin{cases}
\dis v_3\cdot\partial_x (\mathfrak{f}_1-\mathfrak{f}_2)+p_E^0 \FL (\mathfrak{f}_1-\mathfrak{f}_2)=\Gamma(\mathfrak{f}_1-\mathfrak{f}_2,\mathfrak{f}_1)+\Gamma(\mathfrak{f}_2,\mathfrak{f}_1-\mathfrak{f}_2)+\mathfrak{S}_1-\mathfrak{S}_2,\\[2mm]
\dis (\mathfrak{f}_1-\mathfrak{f}_2)(0,v)|_{v_3>0}= P_\g (\mathfrak{f}_1-\mathfrak{f}_2)(0,v)
-(I-P_\g)(\mathfrak{f}_1^\infty-\mathfrak{f}_2^\infty) +\mathfrak{R}_1-\mathfrak{R}_2,\\[2mm]
\lim_{x\rightarrow\infty}(\mathfrak{f}_1-\mathfrak{f}_2)(x,v)=0.
\end{cases}
\end{equation}
We can look at \eqref{5.9} as a linear problem. Taking $\frac12\s_0\leq\s<\s_0$, then we apply  \eqref{1.22} to \eqref{5.9} to have
\begin{align}\label{5.11}
&\|e^{\s x}w(\mathfrak{f}_1-\mathfrak{f}_2)\|_{L^\infty_{x,v}}+|e^{\s x}(\mathfrak{f}_1-\mathfrak{f}_2)(0)|_{L^\infty(\g)} +|(\mathfrak{b}_{1,1}^\infty-\mathfrak{b}_{2,1}^\infty, \, \mathfrak{b}_{1,2}^\infty-\mathfrak{b}_{2,2}^\infty, \, \mathfrak{c}_{1}^\infty-\mathfrak{c}_{2}^\infty)| \nonumber\\
&\leq \frac{2\hat{C}_1}{\s_0-\s}\Big\{ \|e^{\s_0 x}\nu^{-1}w\Gamma(\mathfrak{f}_1-\mathfrak{f}_2,\mathfrak{f}_1)\|_{L^\infty_{x,v}} +\|e^{\s_0 x}\nu^{-1}w\Gamma(\mathfrak{f}_2,\mathfrak{f}_1-\mathfrak{f}_2)\|_{L^\infty_{x,v}}\nonumber\\
&\qquad\qquad+\|e^{\s_0 x}\nu^{-1}w(\mathfrak{S}_1-\mathfrak{S}_2)\|_{L^\infty_{x,v}}+|w(\mathfrak{R}_1-\mathfrak{R}_2)|_{L^\infty_v(\g_-)} \Big\}\nonumber\\
&\leq \frac{2\hat{C}_1}{\s_0-\s}\Big\{ \big[\|e^{\f12\s_0 x}w\mathfrak{f}_1\|_{L^\infty_{x,v}}+ \|e^{\f12\s_0 x}w\mathfrak{f}_2\|_{L^\infty_{x,v}}\big]\cdot \|e^{\f12\s_0 x}w(\mathfrak{f}_1-\mathfrak{f}_2)\|_{L^\infty_{x,v}}\nonumber\\
&\qquad\qquad +\|e^{\s_0 x}\nu^{-1}w(\mathfrak{S}_1-\mathfrak{S}_2)\|_{L^\infty_{x,v}}+|w(\mathfrak{R}_1-\mathfrak{R}_2)|_{L^\infty_v(\g_-)} \Big\}.
\end{align}
Taking $\s=\f12\s_0$, by a direct calculation,  we have from \eqref{5.11} that
\begin{align}\label{5.12}
\|e^{\f12\s_0 x}w(\mathfrak{f}_1-\mathfrak{f}_2)\|_{L^\infty_{x,v}}\leq \f{8\hat{C}_1}{\s_0} \Big\{ \|e^{\s_0 x}\nu^{-1}w(\mathfrak{S}_1-\mathfrak{S}_2)\|_{L^\infty_{x,v}}+|w(\mathfrak{R}_1-\mathfrak{R}_2)|_{L^\infty_v(\g_-)} \Big\}.
\end{align}
Substituting \eqref{5.12} into \eqref{5.11}, one obtains
\begin{align*}
&\|e^{\s x}w(\mathfrak{f}_1-\mathfrak{f}_2)\|_{L^\infty_{x,v}}+|e^{\s x}(\mathfrak{f}_1-\mathfrak{f}_2)(0)|_{L^\infty_v(\R^3)} +|(\mathfrak{b}_{1,1}^\infty-\mathfrak{b}_{2,1}^\infty, \, \mathfrak{b}_{1,2}^\infty-\mathfrak{b}_{2,2}^\infty, \, \mathfrak{c}_{1}^\infty-\mathfrak{c}_{2}^\infty)| \nonumber\\
&\leq \frac{\hat{C}_1}{\s_0-\s} \Big[1+\f{64\hat{C}_1^2\d}{\s_0^2}\Big]\Big\{\|e^{\s_0 x}\nu^{-1}w(\mathfrak{S}_1-\mathfrak{S}_2)\|_{L^\infty_{x,v}} +  |w(\mathfrak{R}_1-\mathfrak{R}_2)|_{L^\infty_v(\R^3_+)}  \Big\},
\end{align*}
which proves \eqref{1.30}. Therefore we complete the proof of Theorem \ref{thm1.1}. $\hfill\Box$

%%%%%%%%%%%%%%%%%%%%%%%%%%%%%%%%%%%%%%%%%%%%%%%%%%%%%%%%%%%%%%%%%%%%%%%%%%%%%%%%%%%%%
\

\appendix
\section{ Some useful known results }\label{APP}

The operator $K$ satisfies the following Grad's estimates
\begin{align*}
Kf(v)=\int_{\mathbb{R}^3}k(v,\eta)f(\eta)\,d\eta,\quad
%	K_2f(v)=\int_{\mathbb{R}^3}k_2(v,\eta)f(\eta)\,dv,\notag
\end{align*}
where  $k(v,\eta)$ %and $k_2(v,\eta)$ satisfy
	%	\begin{align}
	%	0\leq k_1(v,\eta)=c_1|v-\eta| e^{-\f{|v|^2}{4}}e^{-\f{|\eta|^2}{4}},\notag
	%	\end{align}
	%	and
\begin{align}\label{2.15}
0\leq |k(v,\eta)|\leq \f{C}{|v-\eta|}e^{-\f{|v-\eta|^2}{8}}e^{-\f{||v|^2-|\eta|^2|^2}{8|v-\eta|^2}}+C|v-\eta| e^{-\f{|v|^2}{4}}e^{-\f{|\eta|^2}{4}},
\end{align}
where $C>0$ is a given constant. Following \eqref{2.15}, it is direct to have
\begin{align*}%\label{2.16}
\int_{\mathbb{R}^3}\Big|k(v,\eta)\cdot\f{(1+|v|)^{\a}}{(1+|\eta|)^{\a}}\Big|d\eta&\leq  C_{\alpha}(1+|v|)^{-1}.
\end{align*}
It is known \cite{Gl} that $\FL$ satisfies
\begin{equation*}
\int_{\mathbb{R}^3}g\FL g dv \geq c_0\|(\mathbf{I-P}) g \|_{\nu}^2,\\
%	\int_{\mathbb{R}^3}\nu |\FL^{-1} h|^2 dv &\leq \bar{c}_0 \int_{\mathbb{R}^3}\nu^{-1} | h|^2 dv,
\end{equation*}
From Guo \cite{Guo2}, it holds
\begin{align}\label{2.16}
\Gamma(f,f)\in \mathcal{N}^{\perp},\qquad \|\nu^{-1} w \Gamma(f,f)\|_{L^\infty_v}\leq C \|w f \|^2_{L^\infty_v}.
\end{align}

\

We introduce a lemma which will be used to obtain the uniform $L^\infty_{x,v}$ of approximate solutions.
\begin{lemma}\cite{DHWZ} \label{lemA.1}
	Consider  a sequence $\{a_i\}_{i=0}^\infty $  with each $a_i\geq0$. % for $i=0,1,\cdots$.
	For any fixed $k\in\mathbb{N}_+$, we denote $$A_i^k=\max\{a_i, a_{i+1},\cdots, a_{i+k}\}.$$
	
	\noindent{(1)} Assume $D\geq0$.  If $a_{i+1+k}\leq \f18 A_i^{k}+D$ for $i=0,1,\cdots$, then it holds that
	\begin{equation*}%\label{A.1}
	A_i^k\leq \left(\f18\right)^{\left[\frac{i}{k+1}\right]}\cdot\max\{A_0^k, \ A_1^k, \cdots, \ A_k^k \}+\f{8+k}{7} D,\quad\mbox{for}\quad i\geq k+1.
	\end{equation*}
	
	\noindent{(2)} %Let  $\{d_i\}_{i=0}^\infty $ be a decreasing sequence, i.e.,  $d_{i+1}\leq d_i $.
	Let $0\leq \eta<1$ with $\eta^{k+1}\geq\frac14$.  If $a_{i+1+k}\leq \f18 A_i^{k}+C_k \cdot \eta^{i+k+1}$ for $i=0,1,\cdots$, then it holds that
	\begin{align*}%\label{A.1-1}
	A_i^k\leq \left(\f18\right)^{\left[\frac{i}{k+1}\right]}\cdot\max\{A_0^k, \ A_1^k, \cdots, \ A_k^k \}+2C_k\f{8+k}{7} \eta^{i+k},\quad\mbox{for}\quad i\geq k+1.
	\end{align*}
\end{lemma}

\

\noindent{\bf Acknowledgments.} %Y. Guo is supported by an NSF Grant (DMS \#1810868).
Feimin Huang's research is partially supported by National Natural Sciences Foundation of China No. 11688101. Yong Wang's research  is partially supported by National Natural Sciences Foundation of China No. 11688101, 11771429, 12022114, and Youth Innovation Promotion Association of Chinese Academy of Sciences No. 2019002.
% The authors would thank the anonymous referees for the valuable and helpful comments on the paper.

\

\end{document}